\pdfoutput=1
\PassOptionsToPackage{numbers,sort&compress}{natbib}
\documentclass[11pt]{article}

\usepackage{arxiv}
\hypersetup{hidelinks}
\setcitestyle{numbers,square,comma,sort&compress}

\usepackage{amsmath,amsfonts,latexsym,fullpage,graphicx,vmargin, mathrsfs}
\usepackage{comment}
\usepackage[inline]{enumitem}
\usepackage{cleveref}

\usepackage{tikz-cd}
\usepackage{tikz}
\usetikzlibrary{graphs}
\usetikzlibrary{cd}

\usepackage{xcolor}
\usepackage{url}

\usepackage{amssymb}
\usepackage{mathtools}
\usepackage{amsmath,amsfonts}

\usepackage{subcaption}
\usepackage{scalerel}
\usepackage{url}
\usepackage{graphicx}
\usepackage{booktabs}
\usepackage{cleveref}
\usepackage{float}
\floatstyle{plaintop}
\restylefloat{table}

\usepackage{multirow}

\newtheorem{teo}{Theorem}
\newtheorem{prop}{Proposition}
\newtheorem{defi}{Definition}
\newtheorem{rmk}{Remark}
\newtheorem{lem}{Lemma}
\newtheorem{cor}{Corollary}
\newtheorem{ex}{Example}

\newtheorem{assumption}{Assumption}
\newtheorem{propfromcor}[prop]{Proposition}

\newtheorem{remark}[rmk]{Remark}

\crefname{teo}{theorem}{theorems}
\Crefname{teo}{Theorem}{Theorems}
\crefname{prop}{proposition}{propositions}
\Crefname{prop}{Proposition}{Propositions}
\crefname{propfromcor}{proposition}{propositions}
\Crefname{propfromcor}{Proposition}{Propositions}
\crefname{lem}{lemma}{lemmas}
\Crefname{lem}{Lemma}{Lemmas}
\crefname{cor}{corollary}{corollaries}
\Crefname{cor}{Corollary}{Corollaries}
\crefname{conj}{conjecture}{conjectures}
\Crefname{conj}{Conjecture}{Conjectures}
\crefname{rmk}{remark}{remarks}
\Crefname{rmk}{Remark}{Remarks}
\crefname{remark}{remark}{remarks}
\Crefname{remark}{Remark}{Remarks}
\crefname{ex}{example}{examples}
\Crefname{ex}{Example}{Examples}
\crefname{defi}{definition}{definitions}
\Crefname{defi}{Definition}{Definitions}
\crefname{assumption}{assumption}{assumptions}
\Crefname{assumption}{Assumption}{Assumptions}

\usepackage{caption}
\usepackage{subcaption}
\usepackage[section]{placeins}

\DeclareMathOperator{\Top}{\mathbf{Top}}

\DeclareMathOperator{\dTop}{\mathbf{dTop}}

\DeclareMathOperator{\f}{\mathbf{F}}
\DeclareMathOperator{\ff}{\mathbf{f}}

\DeclareMathOperator{\h}{H}

\DeclareMathOperator{\DCh}{DCh}

\DeclareMathOperator{\Fop}{F}

\DeclareMathOperator{\Crit}{Crit}

\DeclareMathOperator{\supp}{supp}

\DeclareMathOperator{\cov}{cov}

\DeclareMathOperator{\diam}{diam}

\DeclareMathOperator{\In}{In}
\DeclareMathOperator{\Out}{Out}
\DeclareMathOperator{\Adm}{Adm}
\DeclareMathOperator{\Hom}{Hom}

\newcommand{\virgolette}[1]{``#1''}

\newcommand{\R}{\mathbb{R}}
\newcommand{\Q}{\mathbb{Q}}
\newcommand{\Z}{\mathbb{Z}}
\newcommand{\N}{\mathbb{N}}

\newcommand{\U}{\mathcal{U}}

\newcommand{\Rcal}{\mathcal{R}}

\newcommand{\Hcal}{\mathcal{H}}
\renewcommand{\O}{\mathcal{O}}
\newcommand{\Ss}{\mathbb{S}^1}

\newcommand{\Dcal}{\mathcal{D}}

\usepackage{scalerel}

\firstpageno{1}

\begin{document}

\title{Periodic Tunnel Systems: Reeb Graph Representations and Circular Max-Flow}

\author{Matteo Pegoraro\thanks{Institute of Computing, Faculty of Informatics,
       Universit\`a della Svizzera italiana, Lugano, Switzerland.\protect\newline
       E-mail: \texttt{matteo.pegoraro@usi.ch}}
  \and
  Lisbeth Fajstrup\thanks{Department of Mathematical Sciences, Aalborg University, Aalborg, Denmark.\protect\newline
       E-mail: \texttt{fajstrup@math.aau.dk}}}

\maketitle

\begin{abstract}
Inspired by the study of transport properties in materials science, we investigate the Reeb graph of a map \(f:X\to\Ss\), interpreting \(X\) as a periodic tunnel system in which transport takes place. The map \(f\) encodes a chosen periodic direction and induces a directed structure on the graph. Under constructibility, we establish a lifting result for directed paths that makes precise how the graph represents tunnel segments in \(X\). To turn this representation into a numerical transport descriptor, we equip the graph with cross-sectional capacities defined, in a Riemannian setting, by weighted Hausdorff-measure integrals. Under boundary regularity and Morse assumptions in dimension \(d\ge2\), we establish continuity over regular Reeb values and continuity of the total capacity across fibers. Taking edgewise infima gives the bottleneck capacities for a circular max-flow formulation, which measures flow through fibers of the map to the circle without choosing an arbitrary inlet or outlet. The resulting descriptor is computable by linear optimization and is related to minimum-cost circulations. For nonnegative real capacities, we prove that the max-flow value is invariant under periodic unrolling: the same descriptor is obtained whether the medium is modeled by a toroidal cell or by its periodic replication along the chosen direction. Simulations on periodic point clouds illustrate the construction and the effects of imposing source--target cuts.
\end{abstract}


\section{Introduction}

In recent years, applied topology has brought fundamental contributions to materials science by
devising principled methods for analyzing 2D and 3D samples that represent portions of a medium of
interest, sometimes encoded as force networks or other geometric complexes
\citep{materials_1, materials_2, materials_3, materials_4}.
The versatility of topological techniques has enabled successful applications across a broad range
of material structures, including amorphous materials \citep{materials_5, sorensen2020revealing},
quasi-crystalline materials \citep{cramer2020evolution}, and crystalline systems
\citep{skraba2012topological, hyde2014unification, edelsbrunner2021density, kurlin2022mathematics, bright2023geographic}.

\medskip

A particularly active direction concerns transport and percolation phenomena in complex
microstructures \citep{patil2021conductive, minamitani2022topological, yang2024topological}. In porous
media this includes permeability, while in electrolyte materials it includes ionic conductivity.
In both cases, the goal is to relate an effective macroscopic behavior to geometric and topological
features of the underlying structure. A fruitful approach, pioneered in the porous-media setting,
combines a graph representation of the medium (viewed as a network of interconnected channels) with
max-flow computations, yielding interpretable transport proxies
\citep{ushizima2012augmented, armstrong2021correspondence}. In \cite{armstrong2021correspondence}, this
kind of max-flow pipeline is also linked formally to a transport quantity under modeling
assumptions.

Max-flow originates in optimization theory \citep{ford1956maximal, schrijver2002history, schrijver2003combinatorial}
and has since developed through algorithms \citep{leighton1999multicommodity, cormen2022introduction, orlin2013max, madry2016computing},
applications \citep{garg1996approximate, boykov2004experimental, dwivedi2011maximum}, and variants adapted to different constraints and domains
\citep{couprie2011combinatorial, yuan2010study, han2014maximum, cui2016quantum, gayathri2024max}.
In our setting, however, the graph on which max-flow is computed is not given a priori: as in \cite{ushizima2012augmented}, it is obtained from the ambient space by a Reeb-graph dimensionality reduction, and its capacity constraints are induced by cross-sectional areas of the corresponding level-set components. The periodic boundary conditions, however, require a circular formulation of flow: the chosen transport direction closes on itself and supplies no distinguished inlet or outlet. This brings the problem into contact with the interaction between flow and the topology of graphs and cell complexes
\citep{chambers2012homology, ghrist2013topological, krishnan2014flow}.

\subsection*{Motivating application: periodic transport in materials}

We briefly outline the atomistic setting that motivates the pipeline developed in this work.
The input consists of atomic positions forming a finite point cloud \(B\), often called the
backbone of the material (in particular for amorphous solids such as glasses), together with
ion trajectories \(I_t\) evolving in time. Chemical interactions strongly constrain ionic motion:
ions rarely approach atoms closer than a material-dependent exclusion distance. This can be modeled by
thickening the backbone with metric balls,
\(\widetilde{B}=\bigcup_{b\in B} B_{\varepsilon_b}(b)\), so that ionic trajectories are effectively
confined to the complementary void space
\(X=(\Ss\times\Ss\times\Ss)\setminus\widetilde{B}\).
The flat torus \(\Ss\times\Ss\times\Ss\) models periodic boundary conditions, which are routinely used to
suppress boundary artifacts in simulations \citep{jensen2017introduction}.

Our aim, like that of the augmented-network approaches of
\cite{ushizima2012augmented, armstrong2021correspondence}, is to extract from the void space a
geometric description of the constraints imposed on directional transport. In this framework,
ionic motion is mediated by the tunnel structure of \(X\): ions can move through connected void
channels, while narrow cross-sections act as bottlenecks. We pursue this aim in two stages of
independent interest. The first is a compact and interpretable directed Reeb-graph representation
of the periodic tunnel system. The second is a real-valued descriptor of that representation,
obtained by augmenting its tunnels with geometric capacities and optimizing a circular flow.
The resulting area-normalized max-flow descriptor has now been used in a materials-science study
of sulfide glass solid-state electrolytes, where it was compared with simulated ionic
conductivity \citep{pegoraro2026ionic}.

The earlier approaches provide an important foundation for both stages, but their source--target
formulations are not directly applicable under periodic boundary conditions. A direct translation
of the porous-media paradigm meets a basic obstruction: periodicity removes any canonical notion
of inlet and outlet. A standard source--sink max-flow therefore requires an arbitrary cut that
declares an ``upper'' and a ``lower'' side of the domain; on a torus, different choices can produce
inconsistent results (cf.\ \Cref{sec:sim_2}).

\subsection*{Detailed Pipeline}

\begin{figure}[!t]
    \centering
    \includegraphics[width=\textwidth]{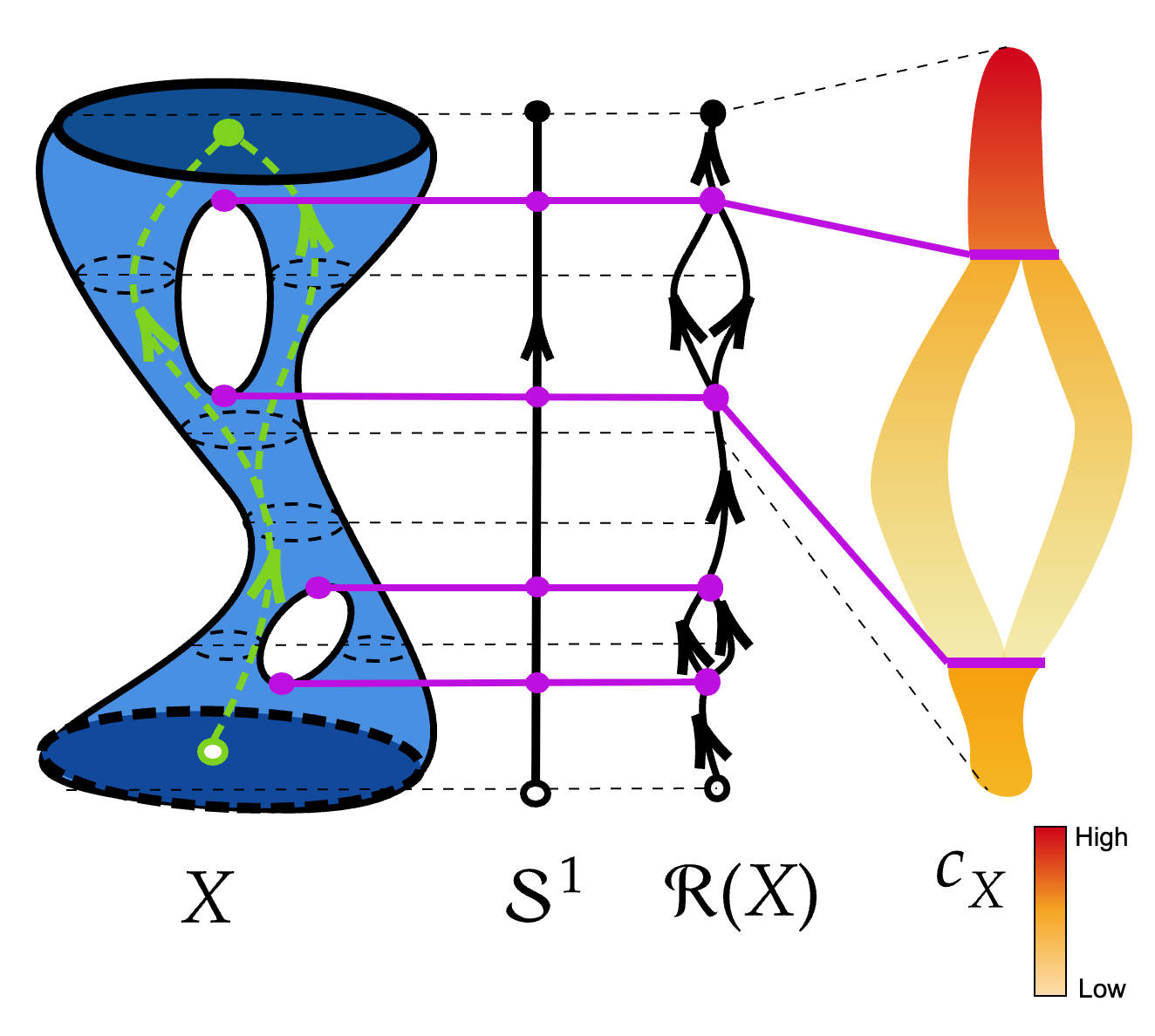}
\caption{Visual summary of the two contributions. Starting from a space \(X\) equipped with a map \(f:X\to\Ss\), represented by horizontal level sets, we first form the directed Reeb graph \(\Rcal(X)\). Its points represent path-connected components of fibers, its arcs represent tunnels, and its direction is induced by the chosen direction on \(\Ss\). The green dashed path indicates how a directed tunnel segment is represented by a directed path in \(X\); the purple points indicate critical points or values. We then augment the representation by the pointwise capacity \(c_X\), shown by the coloring on the right, and summarize the resulting network by circular max-flow. The top and bottom of the displayed objects are identified throughout.}
\label{fig:augmented_reeb}
\end{figure}

We summarize the two principal stages of our construction in \Cref{fig:augmented_reeb}. Some
notation in the figure, including the directed structures and the pointwise capacity function, is
introduced in the technical sections.

\begin{itemize}
\item \textbf{Stage I: a directed Reeb representation of a periodic tunnel system.}
Starting from \(f:X\to\Ss\), we form the Reeb graph \(\Rcal(X)\), which records the evolution of
path-connected components of the level sets of \(f\) around the circular coordinate. In the
materials-science setting these components are tunnel cross-sections, and the graph edges record
the tunnels traced by them. Pulling back the clockwise or counterclockwise direction on \(\Ss\)
turns this graph into a directed representation. We make the tunnel interpretation precise by
showing that directed paths on the graph have representatives given by actual directed paths in
\(X\). In doing this, we use the universal cover
\(\R\to\Ss\) to relate the circular representation to its periodically unrolled counterpart.

\item \textbf{Stage II: a real-valued descriptor of the representation.}
We first attach geometric information to the Reeb graph. A point of \(\Rcal(X)\) represents one
path-connected component of a level set. We can therefore augment each such point with the weighted
\((d-1)\)-dimensional Hausdorff measure
\[
\int_{\pi_f^{-1}(y)\cap\Omega}\mu(p)\,d\Hcal^{d-1}(p)
\]
thereby obtaining a pointwise capacity function on \(\Rcal(X)\).

Continuity at differential regular values and at Morse critical values justifies assigning to
each tunnel edge its bottleneck, namely the infimum of the pointwise capacities along that edge.
We then define circular max-flow for the resulting capacity-augmented directed graph, using bottlenecks to constrain the value of flow in each edge. Because circular max-flow 
is measured through the fibers of a map to \(\Ss\), no source or target is chosen. The resulting value is thus a
single intrinsic throughput descriptor, and we show that it is computable by linear optimization or as a minimum-cost
circulation. We further prove that the same max-flow value
is obtained whether the medium is modeled by a toroidal cell or by its periodic
replication along the chosen direction. This agreement captures the modeling
meaning of periodic boundary conditions: identifying opposite boundaries of a
finite cell represents the periodic continuation of the medium.
\end{itemize}

The first construction is developed in \Cref{sec:reeb}; the geometric capacity construction is
established in \Cref{sec:mapper_max_flow}; the abstract scalar descriptor is developed in
\Cref{sec:circ_flow}; and \Cref{sec:reeb_pipeline} combines the ingredients for periodic tunnel
systems.

\subsection*{Outline}
\Cref{sec:reeb} establishes the directed Reeb-graph representation of periodic tunnel systems,
including pullbacks, directed structures, constructibility, and the lifting of directed tunnel
segments to paths in the original space. \Cref{sec:mapper_max_flow} constructs geometric weights
from fiberwise Hausdorff-measure integrals and proves their compatibility with periodic unrolling.
\Cref{sec:circ_flow} introduces circular max-flow as a real-valued descriptor of a
capacity-augmented circular graph, gives its linear-programming formulation, and proves agreement
with the corresponding unrolled graph. \Cref{sec:reeb_pipeline} finally transfers pointwise capacities to graph arcs and assembles the representation
and the descriptor. \Cref{sec:discussion} discusses the two
contributions and future directions; the appendices contain simulation studies and technical
proofs in the order in which their results are used.

\section{Directed Reeb-Graph Representations of Periodic Tunnel Systems}
\label{sec:reeb}

This section establishes the first contribution: a directed Reeb-graph representation of the
tunnel structure of a space with periodic boundary conditions.

Although the motivating setting is typically a coordinate projection restricted to the void space
\(X\subset\Ss\times\Ss\times\Ss\), the construction requires only a
topological space \(X\) and a continuous map \(f:X\to\Ss\). We therefore work in this generality
and specialize to coordinate projections in the examples and simulations.

\subsection{Reeb graphs and Reeb Spaces}\label{sec:reeb_graph_spaces}

To start off, we introduce Reeb graphs as quotient spaces. In the following, we consider $Y-\Top$, the slice category of $\Top$ over $Y$, which collects the maps from objects in $\Top$ to $Y$.
We refer to an object in $Y-\Top$ as a $Y$-space. We also make use of the following notation: given a topological space $X$, $\O(X)$ is the  set of open sets of $X$.

We need some general facts and definitions about Reeb graphs/spaces and their topology. Whenever there is no practical advantage in stating them for our particular scenario, we state them for a general continuous map $f:X\rightarrow Y$.

\begin{defi}[Reeb space / Reeb graph]\label{def:reeb}
Let \(f:X\to Y\) be a \(Y\)-space. Define an equivalence relation on \(X\) by
\[
x\sim x'
\quad\Longleftrightarrow\quad
\exists\,t\in Y\ \text{such that}\ x,x'\ \text{lie in the same path-connected component of }f^{-1}(t).
\]
The Reeb space of \(f\) is the quotient \(\Rcal(X):=X/\!\sim\), with
quotient map \(\pi_f:X\to \Rcal(X)\). Since \(x\sim x'\) implies \(f(x)=f(x')\), the map \(f\)
factors uniquely through \(\pi_f\): there is an induced continuous map \(\Rcal(f):\Rcal(X)\to Y\)
defined by \(\Rcal(f)([x]):=f(x)\).
Equivalently, we have the commutative diagram
\begin{equation}\label{eq:reeb_quotient}
\begin{tikzcd}
X \ar[r,"f"] \ar[d,"\pi_f"'] & Y \\
\Rcal(X) \ar[ur,"\Rcal(f)"'] &
\end{tikzcd}
\end{equation}
and we refer to \(\Rcal\) as the Reeb functor.
\end{defi}

\begin{remark}\label{rem:reeb-graph-terminology}
In general \(\Rcal(X)\) need not be a graph. In the following we will nevertheless use the term
Reeb graph for \(\Rcal(X)\) as a convenient shorthand, postponing the
precise hypotheses under which \(\Rcal(X)\) admits a graph stratification to later sections.
\end{remark}

\begin{prop}\label{prop:reeb_topology}
Consider a continuous map $f:X\rightarrow Y$ and assume $X$ is locally path-connected.
    Take $V\in \O(Y)$, and $U\in\pi_0(f^{-1}(V))$. We have $\pi_f^{-1}(\pi_f(U))=U$ and so $\pi_f(U)\in \O(\Rcal(X))$. 
\end{prop}

\begin{proof}
Clearly \(U\subseteq \pi_f^{-1}(\pi_f(U))\). Conversely, let
\(y\in \pi_f^{-1}(\pi_f(U))\). Then there exists \(x\in U\) such that
\(\pi_f(y)=\pi_f(x)\). Hence \(x\) and \(y\) lie in the same path-connected
component of the fiber \(f^{-1}(f(x))\). Since \(x\in U\subset f^{-1}(V)\),
we have \(f(x)\in V\), and this fiber is contained in \(f^{-1}(V)\). Therefore
\(x\) and \(y\) lie in the same path-connected component of \(f^{-1}(V)\).
Since \(U\) is the component containing \(x\), we get \(y\in U\). Thus
\(\pi_f^{-1}(\pi_f(U))=U\). Since \(U\) is open in \(X\), the definition of
the quotient topology implies that \(\pi_f(U)\) is open in \(\Rcal(X)\).
\end{proof}

\begin{prop}\label{prop:reeb_lpc}
Assume that $X$ is locally path-connected.

\begin{enumerate}
    \item If $W\subset \Rcal(X)$ is open and path-connected, then $\pi_f^{-1}(W)$ is path-connected.

    \item $\Rcal(X)$ is locally path-connected.
\end{enumerate}
\end{prop}

\begin{proof}
We first prove (1).
Let $W\subset \Rcal(X)$ be open and path-connected.
Since $\pi_f^{-1}(W)$ is open in $X$ and $X$ is locally path-connected, it is enough to prove that
$\pi_f^{-1}(W)$ is connected.

Suppose by contradiction that \(\pi_f^{-1}(W)=A\sqcup B\)
with $A,B$ nonempty, disjoint, open-and-closed in $\pi_f^{-1}(W)$.
Every fiber of $\pi_f$ is a Reeb equivalence class, hence path-connected, so no fiber can meet both $A$ and $B$.
Thus $A$ and $B$ are saturated for the quotient map
\(\pi_f|_{\pi_f^{-1}(W)}:\pi_f^{-1}(W)\to W\).
Therefore $\pi_f(A)$ and $\pi_f(B)$ are disjoint, nonempty, open subsets of $W$ whose union is $W$,
contradicting the connectedness of $W$.
Hence $\pi_f^{-1}(W)$ is connected, and therefore path-connected.

We now prove (2).
Fix $p\in \Rcal(X)$ and let $W$ be an open neighborhood of $p$.
Choose $x\in \pi_f^{-1}(p)$, and let $U$ be the path-connected component of $\pi_f^{-1}(W)$ containing $x$.
Moreover, \(U\) is saturated for the quotient map \(\pi_f\): every fiber of
\(\pi_f\) is path-connected, hence if it
meets the component \(U\), it is contained in \(U\). Since \(U\) is open and
saturated, \(\pi_f(U)\) is open in \(\Rcal(X)\) by the quotient topology.
Thus $\pi_f(U)$ is an open path-connected neighborhood of $p$, proving that $\Rcal(X)$ is locally path-connected.
\end{proof}

\subsection{Pullbacks of Reeb Spaces}
\label{sec:pull_reeb}

To compare the circular tunnel representation with its periodically unrolled counterpart, we need the following definitions and results on pullbacks. Some of these facts may be known, but we could not find direct references in the form needed here, so we give brief proofs.

\begin{defi}
    The pullback of a continuous map $f:X\rightarrow Y$ along a map $p:E\rightarrow Y$ is the set
    \(X\times_Y E=\{(x,e)\in X\times E\mid f(x)=p(e)\}\), along with the
    projections to $X$ and $E$ and the subspace topology.
\[
\begin{tikzcd}
X\times_Y E\ar[r,"f'"]\ar[d,"p'"]&E\ar[d,"p"]\\
X\ar[r,"f"]&Y.
\end{tikzcd}
\]
The above construction gives the pullback functor $p^*:Y-\Top\rightarrow E-\Top$ induced by
\((f:X\rightarrow Y)\mapsto f':X\times_Y E\rightarrow E\).
\end{defi}

\begin{defi}
    A topological space is Hausdorff if every pair of points has disjoint open neighborhoods. A topological space is said to be locally Hausdorff if every point has a neighbourhood that is a Hausdorff space under the subspace topology.
    A topological space is locally compact if for every point there is a compact neighborhood. 
\end{defi}

For the sake of simplicity, in the upcoming lemma we write
\(p^*X:=p^*(f:X\rightarrow Y)=f':X\times_Y E\rightarrow E\).

\begin{lem}\label{lem:pullb_pi0}

Consider \(p:E\rightarrow Y\) in \(\Top\), with \(E\) locally compact and \(Y\) locally Hausdorff. 

Then \(p^*\) preserves colimits.
\end{lem}
\begin{proof}
    In \cite{niefield1982cartesianness} it is proven that if $E$ is locally compact and $Y$ locally Hausdorff, then the functor:
\[
(A\rightarrow Y) \mapsto p^*A=(A\times_Y E\rightarrow E)
\]
has a right adjoint. Denoting its value at \(B\to E\) by
\(B^p\to Y\), the adjunction gives
\[
\Hom_{Y-\Top}(A,B^p)\cong \Hom_{E-\Top}(p^*A,B),
\]
naturally in \(A\) and \(B\).
Thus, \(p^*\), being a left adjoint, preserves colimits.
\end{proof}

In addition, we will also use the standard stability of covering maps under pullback.

\begin{lem}[\cite{calcut2012topological}]\label{lem:pull_cover}
Consider  $f:X\rightarrow Y$ in $\Top$: if $p:E\rightarrow Y$ is a covering map, then
\(X\times_Y E\to X\) is a covering map.
\end{lem}

Now we obtain a result which allows us to explicitly characterize $\Rcal(X\times_Y E)$, which will be used to unroll the graphs as in \Cref{sec:int_flow}.

\begin{prop}\label{lem:pull_back}
In the setting of Lemma~\ref{lem:pullb_pi0}, consider the following pullback diagram in \(\Top\):
\[
\begin{tikzcd}
X\times_Y E\ar[r,"f'"]\ar[d,"p'"]&E\ar[d,"p"]\\
X\ar[r,"f"]&Y .
\end{tikzcd}
\]
Then there is a natural homeomorphism
\[
\Rcal(X\times_Y E)\cong \Rcal(X)\times_Y E .
\]
\end{prop}

\begin{proof}
  Consider the fiber product
\[
X\times_Y X
:=
\{(x,x')\in X\times X : f(x)=f(x')\}.
\]
The Reeb equivalence relation on $X$ is represented by 
\[
R_f
:=
\{(x,x')\in X\times X : x\sim x'\},
\]
endowed with the subspace topology inherited from \(X\times_Y X\).

Let
\(
r_1,r_2:R_f\rightrightarrows X
\)
be the restrictions of the two coordinate projections. Since
\(f\circ r_1=f\circ r_2\), these are morphisms in \(Y\text{-Top}\).
By definition of the Reeb quotient, \(\Rcal(X)\), together with the
quotient map \(\pi_f:X\to\Rcal(X)\), is the coequalizer in
\(Y\text{-Top}\) of
\(
R_f\rightrightarrows X.
\)

By \Cref{lem:pullb_pi0}, the pullback functor
\[
p^*:Y\text{-Top}\longrightarrow E\text{-Top}
\]
preserves colimits. Hence \(\Rcal(X)\times_YE\) is the coequalizer in
\(E\text{-Top}\) of the pulled-back pair
\[
p^*R_f \rightrightarrows X\times_YE.
\]

We now identify \(p^*R_f\) with the Reeb equivalence relation
\(R_{f'}\) associated with
\[
f':X\times_YE\longrightarrow E.
\]
For every \(e\in E\), projection onto the first factor restricts to a
homeomorphism
\[
(f')^{-1}(e)
\xrightarrow{\sim}
f^{-1}(p(e)),
\qquad
(x,e)\longmapsto x.
\]
Therefore \((x,e)\) and \((x',e)\) lie in the same path-connected
component of \((f')^{-1}(e)\) if and only if \(x\) and \(x'\) lie in
the same path-connected component of \(f^{-1}(p(e))\).

More explicitly,
\[
p^*R_f
=
R_f\times_YE
=
\{(x,x',e):x\sim x',\ f(x)=f(x')=p(e)\}.
\]
The map
\[
(x,x',e)
\longmapsto
\bigl((x,e),(x',e)\bigr)
\]
gives the homeomorphism
\((X\times X)\times E
\;\cong\;
(X\times E)\times_E(X\times E),\)
and, by restriction,
identifies \(p^*R_f\) homeomorphically with the subspace
\[
R_{f'}
\subset
(X\times_YE)\times_E(X\times_YE)
\]
consisting of pairs that lie in the same path-connected component of
a fiber of \(f'\). Hence
\(
p^*R_f\cong R_{f'}
\)
as \(E\)-spaces, compatibly with the two projections to
\(X\times_YE\).

By definition, \(\Rcal(X\times_YE)\) is the coequalizer in
\(E\text{-Top}\) of
\[
R_{f'}\rightrightarrows X\times_YE.
\]
Thus \(\Rcal(X\times_YE)\) and \(\Rcal(X)\times_YE\) are coequalizers
of the same pair in \(E\text{-Top}\), and hence are canonically
isomorphic:
\[
\Rcal(X\times_YE)\cong \Rcal(X)\times_YE.
\]
Explicitly, the resulting homeomorphism is
\(
\pi_{f'}(x,e)\longmapsto (\pi_f(x),e).
\)
\end{proof}

As a consequence of the results in this section, we have proven that the following are both  pullback diagrams in $\Top$. Note that the right diagram is obtained by applying the Reeb functor to the left diagram, plus, $\Rcal(p')$ is a covering map.

\begin{equation*}
\begin{tikzcd}
X\times_{\Ss} \R\ar[r,"f'"]\ar[d,"p'"]&\R\ar[d,"p"]\\
X\ar[r,"f"]&\Ss.
\end{tikzcd}
\qquad
\begin{tikzcd}
\Rcal(X\times_{\Ss} \R)\ar[r,"\Rcal(f')"]\ar[d,"\Rcal(p')"]&\R\ar[d,"p"]\\
\Rcal(X)\ar[r,"\Rcal(f)"]&\Ss.
\end{tikzcd}  
\end{equation*}

\subsection{Directed Spaces and Directed Reeb Graphs}\label{sec:dir_nerve}

To encode coherently directions on the edges of the Reeb representation, we induce their orientations from the maps $f:X\rightarrow \Ss$ and $f':X\times_{\Ss} \R\rightarrow \R$. We use definitions and ideas from directed topology, following \cite{fajstrup2016directed}. The reader may refer to \Cref{fig:augmented_reeb} for a visual summary.

\begin{defi}[\cite{fajstrup2016directed}]
    A d-space $(X,dX)$ consists of a topological space $X$ and a subset of paths $dX\subset X^{[0,1]}$, called directed paths, dipaths or d-paths, such that:
    \begin{itemize}
        \item $dX$ contains all constant paths;
        \item $dX$ is closed under reparametrization with continuous nondecreasing functions $f:[0,1]\rightarrow [0,1]$ (not necessarily strictly increasing and not necessarily surjective);
        \item given two paths $\alpha,\beta \in dX$ with $\alpha(1)=\beta(0)$, the composition $\alpha \ast\beta$ is defined as: 
$\alpha \ast\beta(t)=\alpha(2t)$ for $t\in [0,1/2]$ and $\alpha \ast\beta(t)=\beta(2t-1)$ otherwise. The set $dX$ is closed under composition.
    \end{itemize}
    A continuous map $f:X\rightarrow Y$ induces a push-forward $d_f:X^{[0,1]}\rightarrow Y^{[0,1]}$ between the sets of paths: $d_f(\alpha)=f\circ \alpha$. Accordingly, a d-map $f:(X,dX)\rightarrow (Y,dY)$ is a continuous map so that $d_f(dX)\subset dY$. We call $\dTop$ the category of d-spaces with d-maps. 
\end{defi}

\begin{ex}\label{ex:directed_R}
Given $\R$ we can consider the directed structure given by non-decreasing paths.
\end{ex}

\begin{ex}[\cite{fajstrup2016directed}]\label{ex:directed_S1}
Consider now $\Ss\subset \mathbb{C}$ as the unit circle centered at $0$. We can consider the directed structure $d\Ss$ generated by all paths $\alpha:I\rightarrow\Ss\hookrightarrow\mathbb{C}$ so that $t\mapsto  \alpha(t) $ goes counterclockwise. 
\end{ex}

While $d_f$ is the map that acts on paths, $\Dcal_f$ and $\Dcal^f$ are going to be the operators which act on directed structures.

Given a d-space $(X,dX)$ and a continuous function $f:X\rightarrow Y$, the pushforward directed structure $\Dcal_f (dX)$ is defined as the directed structure generated by $d_f (dX)$ (i.e. the smallest one that contains $d_f (dX)$). Conversely, if $(Y,dY)$ is a directed space, the pullback structure is defined as $\Dcal^f(dY):=\{\alpha\in X^I \mid d_f(\alpha)\in dY\}$.
These are all directed structures. In particular, $\Dcal_f(dX)$ is the smallest directed structure so that $f$ is a d-map, while $\Dcal^f(dY)$ is the biggest directed structure so that $f$ is a d-map. 
Note that, in general, $\Dcal_f (\Dcal^f(dY))\subset dY$. 
In this work we will often employ pullback directed structures.

\begin{ex}
The restriction on $\Ss$ of the directed structure in \Cref{ex:directed_S1} can  also be built as the pushforward via $p:\R\rightarrow\Ss$ of the directed structure induced on $\R$ by monotone non decreasing paths w.r.t. its total ordering. Similarly, the pullback of the counterclockwise structure on $\Ss$, gives us the original d-structure on $\R$.
\end{ex}

The following proposition states that, if \(Z\) carries a directed structure and we pull this directed structure back along the objects of \(Z\text{-}\Top\), then every morphism in the slice category is a directed map.

\begin{prop}\label{prop:pull_directed}
    Let $f:X\rightarrow Z$ and $g:Y\rightarrow Z$ be continuous maps and let $(Z,dZ)$ be a directed space. Then any $h:X\rightarrow Y$ such that $f=g\circ h$ is a directed map w.r.t. $\Dcal^f(dZ)$ and $\Dcal^g(dZ)$.
    \begin{proof}
        Let $\alpha \in \Dcal^f(dZ)$. Then $f\circ \alpha = g\circ h\circ \alpha\in dZ$. Thus, $h\circ \alpha \in \Dcal^g(dZ)$.
    \end{proof}
\end{prop}

Using these tools we can obtain limits and colimits in $\dTop$ using the ones in $\Top$. See \cite{fajstrup2016directed} for more details.
Consider for example the following pullback diagram:
\[
\begin{tikzcd}
X\times_Y E\ar[r,"f'"]\ar[d,"p'"]&E\ar[d,"p"]\\
X\ar[r,"f"]&Y.
\end{tikzcd}
\]

Then $d(X\times_Y E)$ is given by all $(\alpha,\beta)\in dX\times dE$ such that, for every $t\in [0,1]$, $f(\alpha(t))=p(\beta(t))$.
One can also check that $d(X\times_Y E)=\Dcal^{f'}(dE)\cap \Dcal^{p'}(dX)$.

\subsection{Countably Constructible Reeb graphs}
\label{sec:construtcible}

In the following, we ask that the maps $f:X\rightarrow \Ss$ and $f':X\times_{\Ss} \R\rightarrow \R$ satisfy certain finiteness conditions, which also imply that the image of the Reeb functor is a Reeb graph.

 To simplify the notation, we may write $f:X\rightarrow M$ to prove results or state definitions holding both for $M$ equal to $\R$ or $\Ss$.
First we report the definition of constructible $\R$-spaces from  \cite{de2016categorified}, and then we extend it.

\begin{defi}[\cite{de2016categorified}]\label{def:de_constr}
 An $\R$-space $f:X\rightarrow \R$ is constructible if 
 there is a finite set $S=\{t_0<t_1<\ldots<t_n\}\subset \R$  satisfying:
 \begin{itemize}
    \item for all $i=0,\ldots,n$  we have a locally path-connected compact space $V_i$;
    \item for all $i=0,\ldots,n-1$ we have  locally path-connected compact space $E_i$;
    \item for all $i=0,\ldots,n-1$ we have continuous maps $l_i:E_i\rightarrow V_i$ and $r_i: E_i\rightarrow V_{i+1}$ 
\end{itemize}   

and such that  $f:X\rightarrow \R$ is homeomorphic to:
\[
\left(\coprod_{i=0}^n (V_i \times \{t_i\})\right) \coprod \left(\coprod_{i=0}^{n-1} (E_i \times [t_i,t_{i+1}]) \right) \big/\sim
\]
with $\sim$ generated by $(l_i(x),t_i)\sim (x,t_i)$ and 
$(r_i(x),t_{i+1})\sim (x,t_{i+1})$ with the map into $\R$ induced by the projection on the second factor of the products. The set $S$ is usually referred to as a critical set.
\end{defi}

We now relax \Cref{def:de_constr} by allowing the gluing procedure to be iterated countably many times, provided that the set of critical values is closed and discrete (in particular, it has no accumulation points).

\begin{defi}\label{def:cont_constr}
An \(\mathbb R\)-space \(f:X\to\mathbb R\) is countably constructible if there
exists a closed discrete set \(S\subset\mathbb R\) such that, for every compact interval
\([a,b]\subset\mathbb R\), the restricted \([a,b]\)-space
\(f^{-1}([a,b])\longrightarrow [a,b]\)
is constructible in the sense of \Cref{def:de_constr}, with critical set contained in
\(\{a,b\}\cup\bigl(S\cap(a,b)\bigr)\).

An \(\Ss\)-space \(f:X\to\Ss\) is constructible if the pullback
\(f':X\times_{\Ss}\mathbb R\longrightarrow\mathbb R\) along the covering map
\(p:\mathbb R\to\Ss\) is countably constructible.

We call (countably) constructible \(M\)-Reeb graphs the Reeb graphs of
(countably) constructible \(M\)-spaces.
\end{defi}

In the sequel we will write ``(countably) constructible'' to mean:
countably constructible when \(M=\R\), and constructible when \(M=\Ss\).
Note that, by definition, if \(f:X\to \Ss\) is constructible then its pullback along \(p:\R\to \Ss\)
is a countably constructible \(\R\)-space.

\begin{rmk}\label{rmk:constructible_reeb}
The point of \Cref{def:cont_constr} is that (countably) constructible \(M\)-spaces yield Reeb spaces
that admit a graph stratification, hence genuine Reeb graphs. This follows by iterating, in a
countable setting, the construction of \(\R\)-Reeb graphs from \(\R\)-spaces in
\cite{de2016categorified}.
\end{rmk}

\begin{defi}
    A continuous map $f:X\rightarrow Y$ is quasi-finite if its fibers have finite cardinality.
\end{defi}

\begin{rmk}
    The map $\Rcal(f):\Rcal(X)\rightarrow M$ of a (countably) constructible $M$-Reeb graph is quasi-finite.
\end{rmk}

We now introduce two standard notions of critical values: one for persistent sets/persistence modules, phrased directly in the setting of Reeb graphs, and another one for smooth maps into $M$. 

\begin{defi}[Topological critical values for a Reeb graph]\label{def:crit_top}
Let \(\Rcal(f):\Rcal(X)\to M\) be a Reeb graph. A value \(t\in M\) is regular for \(\Rcal(f)\) if
there exists an open neighborhood \(U\subset M\) of \(t\) such that, for every \(t'\in U\),
the natural map
\[
\pi_0\bigl(\Rcal(f)^{-1}(t')\bigr)
\longrightarrow
\pi_0\bigl(\Rcal(f)^{-1}(U)\bigr),
\]
induced by the inclusion
\(\Rcal(f)^{-1}(t')\hookrightarrow \Rcal(f)^{-1}(U)\), is a bijection.
A value is critical if it is not regular. We denote by \(\Crit(\Rcal(f))\) the set of critical
values of \(\Rcal(f)\).
\end{defi}

\begin{defi}[Differential critical points and values]
\label{def:crit_diff}
Let $X$ be a $C^2$ $d$-manifold (with or without boundary) and let
$f:X\to M$ ($M=\R$ or $M=\Ss$) be a $C^2$ map.
A point $p\in X$ is a (differential) critical point of $f$ if \(df_p=0\).
A value $t\in M$ is a (differential) critical value of $f$ if $t=f(p)$ for some
critical point $p$. We call values in $M$ differentially
regular if they are not critical. Lastly, a critical point is called Morse if it is isolated and its Hessian is non-degenerate.
\end{defi}

\begin{rmk}\label{rmk:critical}
In the smooth setting, critical values are naturally defined through the
differential notion in \Cref{def:crit_diff}. Morse theory
\citep{milnor1963morse} relates this notion to the topology of level and
sublevel sets. If \(f:X\to\R\) is a Morse function on a compact manifold
(i.e.\ all its critical points are nondegenerate), then the homotopy type,
and in particular \(\pi_0\), of the sublevel sets can change only when a
differential critical value is crossed; the change is determined by the
index of the corresponding critical point(s), namely the number of negative
eigenvalues of the Hessian. Consequently, when \(\Rcal(f)\) is the Reeb
graph of a Morse function, every topological critical value in
\(\Crit(\Rcal(f))\) arises from a differential critical value of \(f\).

Sard's theorem already implies that the set of differential critical values has
measure zero. In our \(\Ss\)-constructible setting, however, we impose the stronger
finiteness condition that the set of topological critical values
\(\Crit(\Rcal(f))\) is finite. In what follows, we use the term
``critical'' whenever the intended notion is clear from the context:
for smooth maps it refers to differential critical points or values, while
for Reeb graphs it refers to
\Cref{def:crit_top}.
\end{rmk}

\begin{rmk}[Local homeomorphism]\label{rmk:height_local_homeo}
Let \(f:X\to M\) be a (countably) constructible \(M\)-space, and let
\(t_0\in M\setminus \Crit(\Rcal(f))\). By constructibility, after possibly
shrinking to an open interval/arc \(J\subset M\) containing \(t_0\) and satisfying
\(J\cap \Crit(\Rcal(f))=\emptyset\), the restriction of \(\Rcal(f)\) over \(J\)
is a finite disjoint union of open graph edges. On each such edge the height map
has no critical value and is monotone with image \(J\). Hence every
path-connected component \(V\) of \(\Rcal(f)^{-1}(J)\) is mapped
homeomorphically onto \(J\) by \(\Rcal(f)\).
Equivalently, for every \(y_0\in \Rcal(f)^{-1}(t_0)\), if \(V\) denotes the
path-connected component of \(\Rcal(f)^{-1}(J)\) containing \(y_0\), then
\(\Rcal(f)|_V:V\longrightarrow J\) is a homeomorphism. In particular, locally over a regular value, the restriction of the height map to each component of its preimage admits a well-defined continuous inverse.
\end{rmk}

Note that for an $M$-space as in \Cref{def:cont_constr}, the set $S$ always contains the critical values of the associated Reeb graph.

\subsection{Directed structure and tunnels}

We now apply the preceding machinery to the spaces in the Reeb-graph framework and make precise how directed graph segments represent tunnels in the original space.

Consider the pullback diagrams
\begin{equation}\label{eq:pullbacks}
\begin{tikzcd}
X\times_{\Ss} \R\ar[r,"f'"]\ar[d,"p'"]&\R\ar[d,"p"]\\
X\ar[r,"f"]&\Ss,
\end{tikzcd}
\qquad
\begin{tikzcd}
\Rcal(X)\times_{\Ss} \R\ar[r,"\Rcal(f')"]\ar[d,"\Rcal(p')"]&\R\ar[d,"p"]\\
\Rcal(X)\ar[r,"\Rcal(f)"]&\Ss .
\end{tikzcd}
\end{equation}
By \Cref{lem:pull_back}, the diagram on the right is obtained from the one on
the left by applying the Reeb functor.

Using \Cref{ex:directed_S1}, we endow \(\Ss\) with its clockwise or counterclockwise directed structure and pull that structure back through all the maps in
the diagrams in \ref{eq:pullbacks}. We denote the resulting directed structures by
\[
d\Ss,\qquad dX,\qquad d\Rcal(X),\qquad d(X\times_{\Ss}\R),\qquad
d(\Rcal(X)\times_{\Ss}\R).
\]
By construction, these directed structures satisfy
\[
d(X\times_{\Ss}\R)
=
\Dcal^{f'}(\Dcal^{p}(d\Ss))
=
\Dcal^{p'}(\Dcal^{f}(d\Ss)),
\]
\[
d(\Rcal(X)\times_{\Ss}\R)
=
\Dcal^{\Rcal(f')}(\Dcal^{p}(d\Ss))
=
\Dcal^{\Rcal(p')}(\Dcal^{\Rcal(f)}(d\Ss)),
\]
and, using the commutativity of the Reeb quotient diagram,
\[
dX=\Dcal^{\pi_f}(d\Rcal(X)),
\qquad
d(X\times_{\Ss}\R)=\Dcal^{\pi_{f'}}\bigl(d(\Rcal(X)\times_{\Ss}\R)\bigr).
\]

We now state the path-lifting result showing how directed paths in the Reeb graph
represent directed paths in the original space.

\begin{prop}\label{prop:path_conn_regular}
Let \(f:X\to \Ss\) be a constructible \(\Ss\)-space, and let
\(\Rcal(f):\Rcal(X)\to\Ss\) be the associated constructible Reeb graph. Let \(\alpha:[0,1]\to \Rcal(X)\) be an
embedding.

Then there exists a path \(\alpha':[0,1]\to X\) such that
\[
\pi_f(\alpha'([0,1]))=\alpha([0,1]),
\qquad
\pi_f(\alpha'(0))=\alpha(0),\qquad
\pi_f(\alpha'(1))=\alpha(1).
\]

Moreover, if \(\alpha\in d\Rcal(X)\), then \(\alpha'\) may be chosen in \(dX\).
\end{prop}

The proof is given in Appendix~\ref{app:reeb_pipeline_proofs}. It subdivides the
arc into edgewise pieces, lifts them through constructible presentations, and
glues the lifts inside path-connected Reeb fibers.

Consequently, every directed tunnel segment in \(\Rcal(X)\)
is the \(\pi_f\)-image of an actual directed path in \(X\). The same conclusion
holds for the pullback objects \(X\times_{\Ss}\R\) and
\(\Rcal(X)\times_{\Ss}\R\). Thus the induced directed structures are compatible
with both the Reeb quotient maps and the maps into \(\Ss\) and \(\R\); see
\Cref{fig:augmented_reeb}.

\section{Capacity Functions}
\label{sec:mapper_max_flow}

Building on the representation in \Cref{sec:reeb}, we now augment it with geometric information measuring the cross-sections of the represented tunnels. A point of the Reeb graph represents a path-connected component of a level set of
\(f\), so its natural cross-sectional weight is obtained by integrating a continuous function over
that component with respect to \((d-1)\)-dimensional Hausdorff measure. The purpose of this
section is to show that these weights vary continuously under mild regularity assumptions and
therefore define stable pointwise capacities.
Notably, our setting covers spaces of interest like the offset geometries arising in the materials-science model (unions of metric balls
and their complementary voids) and smooth Morse-type functions on compact manifolds. 

The long
technical proofs are deferred to Appendix~\ref{app:capacity_proofs}: at regular values they use
submersion charts and Sard-type control of exceptional boundary slices, while at Morse critical values
they use the explicit local models supplied by the Morse lemma.

We begin by recalling the definition  of the Hausdorff
measure.

\begin{defi}
    Let $(X,d_X)$ be a metric space and consider $S\subset X$. For any $\delta>0$ and $d\in \N$, $d\ge 1$, define:
    \[
    \Hcal_\delta^d(S) := \frac{\omega_d}{2^d} \inf_{\U \in \cov_\delta(S)}\sum_{U\in \U} \diam(U)^d,
    \]
    where \(\omega_d\) denotes the \(d\)-dimensional Lebesgue measure of
the unit ball in \(\R^d\) and $\cov_\delta(S)$ collects all countable covers of $S$ by sets of diameter less than $\delta$. Lastly, we define the $d$-dimensional Hausdorff measure of $S$ as:
    \[
    \Hcal^d(S)= \sup_{\delta>0}\Hcal_\delta^d(S) = \lim_{\delta \rightarrow 0^+}\Hcal_\delta^d(S).
    \]
\end{defi}

We recall that on a \(d\)-dimensional \(C^1\) Riemannian manifold, the
\(d\)-dimensional Hausdorff measure associated with the induced distance
coincides with the Riemannian volume measure (see, e.g., \cite{taylor2006measure}).
Throughout, Hausdorff measures are defined using the ambient metric.

We now establish the first key continuity statements needed in this section: the
fiberwise integral of a continuous weight over the \((d-1)\)-dimensional level
sets of \(f\) varies continuously as the level crosses a regular value of
a \(C^1\) extension of \(f\). We start by considering non-negative weight functions in the Euclidean setting, then move to the manifold setting, and, lastly, extend the result to general continuous functions.

\begin{lem}\label{prop:hausd_continuous}
Let \(\Omega\subset\mathbb R^d\) be a bounded open set, let
\(\mu:\overline\Omega\to[0,\infty)\) be continuous, and let
\(f:U\to\mathbb R\) be a \(C^1\) function on an open neighborhood
\(U\) of \(\overline\Omega\).
Assume that
\[
\partial\Omega\subset\bigcup_{j=1}^N\Sigma_j,
\]
where each \(\Sigma_j\) is a compact \(C^1\) submanifold of \(U\)
of dimension at most \(d-1\), with \(C^1\) boundary.
Define
\[
c(t):=\int_{\Omega\cap f^{-1}(t)}\mu(x)\,d\mathcal H^{d-1}(x),
\]
with value zero if \(\Omega\cap f^{-1}(t)=\emptyset\).
Fix \(t_0\in\mathbb R\) and assume that:
\begin{enumerate}
\item[\textup{(i)}] \(\nabla f(x)\ne0\) for every
\(x\in\overline\Omega\cap f^{-1}(t_0)\);
\item[\textup{(ii)}] for every \(j\) with \(\dim\Sigma_j=d-1\),
\[
\mathcal H^{d-1}\Bigl(\bigl\{x\in\Sigma_j^\circ\cap\overline\Omega\cap f^{-1}(t_0):
d(f|_{\Sigma_j^\circ})(x)=0\bigr\}\Bigr)=0,
\qquad \Sigma_j^\circ:=\Sigma_j\setminus\partial\Sigma_j.
\]
\end{enumerate}
Then \(c\) is continuous at \(t_0\).
\end{lem}

The proof is given in Appendix~\ref{app:capacity_proofs}.
The same argument gives the following localized version, in which the
hypotheses are imposed only near the relevant part of the weight's support.

\begin{cor}[Localization]\label{cor:hausd_localized}
Let \(\Omega\subset\mathbb R^d\) be a bounded open set, let
\(\mu:\overline\Omega\to[0,\infty)\) be continuous, and let
\(f:U\to\mathbb R\) be a \(C^1\) function on an open neighborhood
\(U\) of \(\overline\Omega\).
Let \(c\) be the slice integral defined in \Cref{prop:hausd_continuous}.
Fix \(t_0\in\mathbb R\). Suppose there exist an open interval
\(I\ni t_0\) and a compact set \(K\subset\overline\Omega\) such that
\[
\operatorname{supp}_{\overline\Omega}(\mu)\cap f^{-1}(I)\subset K,
\]
where the support is taken relative to \(\overline\Omega\).
Assume that there is an open neighborhood \(W\subset U\) of \(K\) such that
\[
\partial\Omega\cap W\subset\bigcup_{j=1}^N\Sigma_j,
\]
where each \(\Sigma_j\) is a compact \(C^1\) submanifold of \(U\)
of dimension at most \(d-1\), with \(C^1\) boundary.
Assume furthermore that:
\begin{enumerate}
\item[\textup{(i)}] \(\nabla f(x)\ne0\) for every
\(x\in K\cap f^{-1}(t_0)\);
\item[\textup{(ii)}] for every \(j\) with \(\dim\Sigma_j=d-1\),
\[
\mathcal H^{d-1}\Bigl(\bigl\{x\in\Sigma_j^\circ\cap K\cap f^{-1}(t_0):
d(f|_{\Sigma_j^\circ})(x)=0\bigr\}\Bigr)=0,
\qquad \Sigma_j^\circ:=\Sigma_j\setminus\partial\Sigma_j.
\]
\end{enumerate}
Then \(c\) is continuous at \(t_0\).
Taking \(K=\overline\Omega\) and \(W=U\) recovers
\Cref{prop:hausd_continuous}.
\end{cor}

The proof is given in Appendix~\ref{app:capacity_proofs}.

Informally, the different sets involved in these two statements serve the following purposes (see also \Cref{fig:local_charts}): \(\Omega\) is the integration domain, while \(\overline\Omega\)
is the compact domain of the weight. The open set \(U\) provides the ambient domain needed
to construct submersion charts, including at boundary points.
The compact set \(K\) localizes the support of \(\mu\) at levels near
\(t_0\), while \(W\) localizes the covering of \(\partial\Omega\)
around \(K\). Note that no regularity of \(\partial K\) is assumed.

\begin{rmk}\label{rmk:hausd_continuous_transversality}
Assumption \textup{(ii)} is necessary in general: if $f$ is regular on $\overline\Omega$ but is tangent to a
$(d\!-\!1)$-dimensional boundary stratum $\Sigma_j$ at level $t_0$, then the slices
$\Omega\cap f^{-1}(t)$ can gain/lose $(d\!-\!1)$-dimensional pieces as $t$ crosses $t_0$, producing
a jump discontinuity in $c(t)$. For instance, if we consider $\Omega = (0,3)\times(0,3)\setminus [1,2]\times[1,2]$ with $f(x,y)=y$ we obtain discontinuities for $t_0\in \{1,2\}$.
\end{rmk}

\begin{teo}\label{cor:manifold}
Let $T$ be a compact $C^2$ $d$-dimensional Riemannian manifold without boundary.

Let $\Omega\subset T$ be an open subset with compact closure $\overline{\Omega}\subset T$, and assume that
\(\partial\Omega\subset\bigcup_{j=1}^J\Sigma_j\),
where each $\Sigma_j$ is a compact $C^1$ submanifold of $T$ of dimension at most $d-1$,
with $C^1$ boundary.

Let $\mu:\overline{\Omega}\to[0,\infty)$ be continuous, and let $f$ be a $C^1$ function on an open
neighborhood $U\subset T$ of $\overline{\Omega}$. Let $t_0\in\mathbb{R}$ and assume:
\begin{enumerate}
\item[\textup{(i)}] regularity on $\overline\Omega$: $df_x\neq 0$ for all
$x\in \overline{\Omega}\cap f^{-1}(t_0)$;
\item[\textup{(ii)}] a.e.\ boundary regularity at $t_0$: for every $j$ with $\dim\Sigma_j=d-1$,
\[
\mathcal{H}^{d-1}\Bigl(\bigl\{x\in \Sigma_j^\circ\cap \overline{\Omega}\cap f^{-1}(t_0):
d(f|_{\Sigma_j^\circ})(x)=0\bigr\}\Bigr)=0,
\qquad \Sigma_j^\circ:=\Sigma_j\setminus\partial\Sigma_j.
\]
\end{enumerate}
Then, on an open interval \(I\ni t_0\), the expression
\[
c(t):=\int_{\Omega\cap f^{-1}(t)}\mu(x)\,d\mathcal H^{d-1}(x)
\]
defines a finite real-valued function, with value zero for an empty slice,
and \(c\) is continuous at \(t_0\).
\end{teo}

The proof is given in Appendix~\ref{app:capacity_proofs}.

\begin{cor}[Continuity at regular values]\label{cor:hausd}
Within the same setting as in \Cref{cor:manifold}, let
$\mu:\overline{\Omega}\to\mathbb{R}$ be a continuous function.
For every regular value \(t_0\) of \(f\) at which the a.e.\ boundary regularity
assumption of \Cref{cor:manifold} holds, the expression
\[
t\longmapsto\int_{\Omega\cap f^{-1}(t)}\mu(x)\,d\mathcal H^{d-1}(x)
\]
defines a finite real-valued function on a neighborhood of \(t_0\), and this
function is continuous at \(t_0\).
\end{cor}

The proof is given in Appendix~\ref{app:capacity_proofs}.

We next extend the continuity statements to Morse critical values attained at interior critical points of $\Omega$.

\begin{teo}[Continuity at a Morse critical value]\label{prop:morse_cont_man}
Let \(T\) be a compact \(C^2\) \(d\)-dimensional Riemannian manifold without boundary, $d\ge 2,$
and let \(\Omega\subset T\) be an open set with compact closure \(\overline\Omega\).
Let \(\mu:\overline\Omega\to[0,\infty)\) be continuous.

Assume that the boundary \(\partial\Omega\) in \(T\) is contained in a finite union
\(\partial\Omega\subset\bigcup_{j=1}^N\Sigma_j\),
where each \(\Sigma_j\) is a compact \(C^1\) submanifold of \(T\) of dimension at most \(d-1\),
with \(C^1\) boundary.

Let \(f:U\to \R\) be a \(C^2\) map defined on an open neighborhood \(U\subset T\) of \(\overline\Omega\).
Fix \(t_0\in\R\) and assume:
\begin{enumerate}
\item[(a)] \(f^{-1}(t_0)\cap \overline\Omega\) contains finitely many critical points \(p_1,\dots,p_L\),
all contained in \(\Omega\), and each \(p_\ell\) is a Morse critical point of \(f\).
\item[(b)] There exists \(\varepsilon>0\) such that \(f\) has no critical points in
\(\overline\Omega\cap f^{-1}\bigl((t_0-\varepsilon,t_0+\varepsilon)\bigr)\) other than \(p_1,\dots,p_L\).
\item[(c)] a.e.\ boundary regularity at \(t_0\): for every \(j\) with \(\dim\Sigma_j=d-1\),
\[
\Hcal^{d-1}\Bigl(\Bigl\{x\in \Sigma_j^\circ\cap \overline\Omega\cap f^{-1}(t_0):\ d(f|_{\Sigma_j^\circ})(x)=0\Bigr\}\Bigr)=0,
\qquad \Sigma_j^\circ:=\Sigma_j\setminus\partial\Sigma_j.
\]
\end{enumerate}
Define \(c(t):=\int_{\Omega\cap f^{-1}(t)}\mu\,d\Hcal^{d-1}\),
with the convention \(c(t)=0\) if \(\Omega\cap f^{-1}(t)=\emptyset\).
Then \(c\) is continuous at \(t_0\).
\end{teo}

The proof is given in Appendix~\ref{app:capacity_proofs}.

We are now ready to state our standing assumptions on the data
$f:X\to \Ss$. In particular, we want a framework that (a)
allows us to use either \Cref{cor:hausd} or \Cref{prop:morse_cont_man} at each value in the image of $f$
and (b) is compatible with the constructibility
conditions of \Cref{sec:construtcible}.

\begin{assumption}\label{assumptions}
We fix:
\begin{itemize}
  \item a compact $C^2$ $d$-dimensional manifold $T$ without boundary, with $d\ge2$,
        equipped with a $C^1$ Riemannian metric;
  \item an open set $\Omega\subset T$ with compact closure $X=\overline{\Omega}$;
  \item the target space $\Ss$ with either the clockwise or counterclockwise directed structure;
  \item a map $f:X\to \Ss$, which gives a constructible $\Ss$-space;
  \item $X$ is a directed space, with the directed structure pulled back from $\Ss$ via $f$;

  \item $f$ extends to a $C^2$ map $f_U:U\to\Ss$ defined on an open
        neighborhood $U\subset T$ of $X$;

  \item the boundary of $\Omega$ in $T$ satisfies
  \[
     \partial\Omega \;\subset\; \bigcup_{j=1}^J \Sigma_j,
  \]
  where each $\Sigma_j$ is a compact $C^1$ submanifold of $T$ of dimension at most
  $d-1$, with $C^1$ boundary;

  \item a.e.\ boundary regularity for $f_U$: for every $t\in\Ss$ and every $j$ with $\dim\Sigma_j=d-1$,
  \[
  \Hcal^{d-1}\Bigl(\Bigl\{x\in \Sigma_j^\circ\cap X\cap f_U^{-1}(t):\
  d\!\bigl((f_U)|_{\Sigma_j^\circ}\bigr)(x)=0\Bigr\}\Bigr)=0,
  \qquad \Sigma_j^\circ:=\Sigma_j\setminus\partial\Sigma_j;
  \]

  \item all critical points of $f_U$ contained in $X$ lie in $\Omega$
        and are Morse; in particular $f_U$ has no critical points on $\partial\Omega$.
\end{itemize}
\end{assumption}

We briefly comment on \Cref{assumptions}, showing that it is general
enough to encompass several scenarios of interest.

\begin{ex}\label{ex:morse}
In \Cref{assumptions} we may take \(\Omega=T\), so that
\(X=\overline{\Omega}=T\) and \(\partial\Omega=\emptyset\).
Hence the boundary decomposition and the a.e.\ boundary regularity
requirements are vacuous.

Suppose that
\(
f:T\to\Ss
\)
is a \(C^2\) constructible \(\Ss\)-space in the sense of
\Cref{def:cont_constr}, and that all its differential critical points
are Morse in the sense of \Cref{def:crit_diff}. Taking \(U=T\) and
\(f_U=f\), and equipping \(X=T\) with the directed structure pulled
back from \(\Ss\), all the conditions in \Cref{assumptions} are
satisfied.

Ordinary Morse functions give a particular class of such examples.
Indeed, let \(g:T\to\R\) be a \(C^2\) Morse function and let
\(I\subset\R\) be a compact interval containing \(g(T)\). If
\(
\iota:I\hookrightarrow\Ss
\)
is a \(C^2\) embedding, then
\(
f:=\iota\circ g:T\to\Ss
\)
has the same critical points as \(g\), and their nondegeneracy is
preserved. Moreover, the usual Morse-theoretic decomposition of
\(g\) gives the constructibility required in
\Cref{def:cont_constr}. Thus \(f\) satisfies the preceding conditions.
\end{ex}

\begin{ex}\label{ex:atoms}
Let \(T=\Ss\times\Ss\times\Ss\) be endowed with the product Riemannian metric,
and let \(B\subset T\) be a finite set. Define
\(d_B(x):=\min_{b\in B}d(x,b)\), and let \(\Omega\) be either a sublevel or a
superlevel set of \(d_B\), for instance
\(\Omega=d_B^{-1}((-\infty,r))\) or
\(\Omega=d_B^{-1}((r,+\infty))\).
Set \(X:=\overline\Omega\), which is compact since \(T\) is compact.

Because \(B\) is finite, the boundary of \(\Omega\) is covered by finitely many
compact smooth submanifold pieces of \(T\), of dimension at most \(2\), with
possibly nonempty smooth boundary. If \(f:T\to\Ss\) is a coordinate projection,
then its critical set on each two-dimensional boundary piece at a fixed level is
\(\Hcal^2\)-null. Thus the boundary decomposition and the a.e. boundary regularity
condition in \Cref{assumptions} both hold. Taking \(f_U=f\) on \(U=T\), the map
\(f_U\) has no differential critical points, so the Morse-critical-point condition
is vacuous.

The map $f:X\to S^1$ is also constructible. As detailed in
Appendix~\ref{app:capacity_proofs}, its real-valued pullback restricted to one period
is a subanalytic function on a compact subanalytic set.
Triangulation gives a piecewise-linear model, hence a
constructible presentation, which extends periodically.

Therefore the thickened finite-backbone examples arising from sublevel or
superlevel sets of \(d_B\) satisfy \Cref{assumptions}.
\end{ex}

As a last step of this section, we show that the fiberwise integrals
\[
t\longmapsto\int_{\Omega\cap f^{-1}(t)}\mu(p)\,d\Hcal^{d-1}(p)
\]
can be decomposed along the fibers of the Reeb graph. As in \Cref{ex:morse}, embedding the image of an \(\R\)-valued map into \(\Ss\) leaves the Reeb quotient and the fiberwise capacities unchanged. We therefore formulate the construction only for \(f:X\to\Ss\).

\begin{defi}[Pointwise capacity function]\label{def:pointwise_capacity}
A function \(c:\Rcal(X)\to\R_{\ge 0}\) is called a pointwise capacity function if:
\begin{enumerate}
    \item for every \(y_0\in\Rcal(X)\) with \(\Rcal(f)(y_0)\) a regular value of \(\Rcal(f)\),
    the function \(c\) is continuous at \(y_0\);
    \item the function
    \(\Ss\longrightarrow\R_{\ge0},\ t\longmapsto\sum_{y\in\Rcal(f)^{-1}(t)}c(y)\)
    is continuous.
\end{enumerate}
\end{defi}

In other words, a pointwise capacity function assigns a nonnegative ``cross-sectional weight''
to each tunnel component of the Reeb graph, varying continuously along the portions of the graph
lying over regular values, and in such a way that the total capacity of each fiber of \(\Rcal(f)\)
depends continuously on the parameter \(t\in \Ss\).

\begin{prop}\label{cor:capacity}
Let \(f:X\to \Ss\) satisfy \Cref{assumptions} and let \(\mu:X\to\R_{\ge 0}\) be
continuous.
Define
\[
c:\Rcal(X)\longrightarrow \R_{\ge 0},\qquad
c(y)\;:=\;\int_{\pi_f^{-1}(y)\cap \Omega}\mu(p)\,d\Hcal^{d-1}(p),
\]
with the convention that the integral is \(0\) if \(\pi_f^{-1}(y)\cap\Omega=\emptyset\).
Then \(c\) is a pointwise capacity function, and, in particular, for every \(t\in \Ss\),
\[
\int_{\Omega\cap f^{-1}(t)} \mu\,d\Hcal^{d-1}
\;=\;
\sum_{y\in \Rcal(f)^{-1}(t)} c(y).
\]
\end{prop}

\begin{proof}
For every \(t\in\Ss\), the sets \(\pi_f^{-1}(y)\cap\Omega\), with
\(y\in\Rcal(f)^{-1}(t)\), are the intersections of \(\Omega\) with the
path-connected components of \(f^{-1}(t)\). The Reeb fibres are closed,
and constructibility makes \(\Rcal(f)^{-1}(t)\) finite. Additivity of the
integral therefore gives
\[
\sum_{y\in\Rcal(f)^{-1}(t)}c(y)
=\int_{\Omega\cap f^{-1}(t)}\mu\,d\Hcal^{d-1}.
\]

To prove condition \textup{(2)}, we check that the continuity results above
apply to the right-hand side. By \Cref{assumptions}, the ambient manifold,
the boundary of \(\Omega\), and the map \(f_U\) satisfy the required regularity
assumptions; in particular, the a.e.\ boundary condition holds at every level.
The weight \(\mu\) is continuous and nonnegative by hypothesis.
At a differential regular value,
\Cref{cor:hausd} applies. At a differential critical value, all critical
points in \(X\) are Morse and lie in \(\Omega\). They are finite in number,
since \(X\) is compact and these critical points are isolated. Thus a
sufficiently small target neighborhood contains no other critical values,
as required by \Cref{prop:morse_cont_man}. We use a local real coordinate
on \(\Ss\), since the proofs of these results only require a neighborhood
of the level under consideration. Consequently, the right-hand side is
finite and continuous at every \(t\in\Ss\). This also proves that each
\(c(y)\) is finite, and establishes condition \textup{(2)} of
\Cref{def:pointwise_capacity}.

For condition \textup{(1)}, fix \(y_0\) over a regular value \(t_0\) of
\(\Rcal(f)\). By \Cref{rmk:height_local_homeo}, there is an open neighborhood
\(V\) of \(y_0\) such that
\(\Rcal(f)|_V:V\to J\) is a homeomorphism onto an open arc \(J\subset\Ss\).
Write \(\theta:J\to V\) for its inverse. Choose a smaller open arc
\(J_0\ni t_0\) with \(\overline{J_0}\subset J\), and a continuous function
\(\eta:\Rcal(X)\to[0,1]\) such that
\[
\eta=1\quad\text{on }\theta(\overline{J_0}),
\qquad \supp\eta\subset V.
\]
Such a function exists since \(\Rcal(X)\) is a finite graph. Set
\(\nu:=\mu\,(\eta\circ\pi_f)\). Although \(\eta\) isolates the chosen branch,
it is defined on the whole Reeb graph, so \(\nu\) is a continuous,
nonnegative weight on all of \(X\), vanishing outside \(\pi_f^{-1}(V)\).
The domain \(\Omega\), its boundary, and the map \(f_U\) are unchanged.
Thus the same regular and Morse continuity results apply with \(\nu\)
in place of \(\mu\), exactly as for condition \textup{(2)}.

For every \(t\in J_0\), the only point of \(\Rcal(f)^{-1}(t)\) at which
\(\eta\) can be nonzero is \(\theta(t)\), where its value is one. Hence
\[
c(\theta(t))
=\int_{\Omega\cap f^{-1}(t)}\nu\,d\Hcal^{d-1}.
\]
The right-hand side is continuous at \(t_0\). Since \(\theta\) is a
homeomorphism, \(c\) is continuous at \(y_0\), proving condition
\textup{(1)}.
\end{proof}

\subsection{Compatibility of capacities under pullback.}\label{sec:capacity_pullback}

We close the section by checking that the pointwise capacities are unchanged by periodic
unrolling. 

Let \(p':X\times_{\Ss}\R\to X\) be the canonical projection, and let
\[
\Rcal(p'):\Rcal(X\times_{\Ss}\R)\longrightarrow \Rcal(X)
\]
be the induced map on Reeb graphs. We regard \(X\times_{\Ss}\R\) as a subset of
\(T\times\R\), endowed with the product Riemannian metric on the ambient manifold
\(T\times\R\) for the purpose of computing Hausdorff measures.

For every \(r\in\R\), the level set of \(f'\) over \(r\) is
\[
(f')^{-1}(r)=\{(x,r): f(x)=p(r)\}\subset T\times\{r\}.
\]
Since the projection \(T\times\{r\}\to T\) is an isometry, the restriction
\[
p':(f')^{-1}(r)\longrightarrow f^{-1}(p(r))
\]
preserves the ambient \((d-1)\)-dimensional Hausdorff measure of the corresponding
level-set components.

Set \(\Omega':=(p')^{-1}(\Omega)\) and \(\mu':=\mu\circ p'\).
Write \(c_X\) for the capacity of \Cref{cor:capacity}, and define
\(c_{X\times_{\Ss}\R}\) by the same fiberwise integral using \(\Omega'\)
and \(\mu'\).

Let \(y'\in\Rcal(X\times_{\Ss}\R)\). By the compatibility of Reeb graphs with pullbacks,
the map \(p'\) identifies the Reeb fiber component represented by \(y'\) with the Reeb fiber
component represented by \(\Rcal(p')(y')\). More precisely,
\[
p'\bigl(\pi_{f'}^{-1}(y')\cap\Omega'\bigr)
=
\pi_f^{-1}\bigl(\Rcal(p')(y')\bigr)\cap\Omega,
\]
and this identification preserves the relevant \((d-1)\)-dimensional Hausdorff measure.
Since \(\mu'=\mu\circ p'\), the fiberwise integrals defining the pointwise capacities agree:
\[
\begin{aligned}
c_{X\times_{\Ss}\R}(y')
&=
\int_{\pi_{f'}^{-1}(y')\cap\Omega'} \mu'\,d\Hcal^{d-1} \\
&=
\int_{\pi_f^{-1}(\Rcal(p')(y'))\cap\Omega} \mu\,d\Hcal^{d-1} \\
&=
c_X\bigl(\Rcal(p')(y')\bigr).
\end{aligned}
\]
Equivalently,
\[
c_{X\times_{\Ss}\R}
=
c_X\circ\Rcal(p').
\]

\section{Max-Flow for Circular Graphs}
\label{sec:circ_flow}

The preceding sections produced a directed Reeb graph representation of a periodic tunnel system and a
geometric mechanism for assigning capacities to its points. We now develop the second principal
contribution at an abstract graph level: a real-valued, max-flow-inspired descriptor for a
capacity-augmented directed graph mapped to the circle. The definition is source-free, so it is
intrinsic to the circular setting; we also give its optimization formulation and establish
its agreement with the corresponding periodically unrolled graph.

\subsection{Simplicial Complexes}
\label{sec:simpl_compl}

We are going to use the language of simplicial complexes to talk about graphs and maps between graphs. We briefly recall below the definitions needed here.

\begin{defi}
    Given a set $V$, an (abstract) simplicial complex on $V$ is a collection of subsets of $V$, that is, $K\subset 2^{V}$,  such that:
    \begin{itemize}
        \item $\{v\}\in K$ for every $v\in V$;
        \item for every $\sigma\in K$, $2^\sigma\subset K$.
    \end{itemize}
\end{defi}

\begin{defi}
A graph is a simplicial complex $(K^0,K^1)$, with $K^0$ being the vertices (identifying $\{v\}$ with $v$) and $K^1$ the edges. 
\end{defi}

\begin{defi}
Given two simplicial complexes \(K\) and \(K'\), a map
\(f:K\to K'\) is called simplicial if it is induced by a function on vertices, in the sense that, for every simplex
\(\sigma=\{x_0,\ldots,x_n\}\in K\), one has
\[
f(\sigma)
=
\{f(x_0),\ldots,f(x_n)\}
\in K',
\]
where, with a slight abuse of notation, we also write \(f\) for the restriction of \(f\) to \(K^0\). If \(f\) preserves the dimension of simplices, it is called rigid \citep{herlihy2013distributed}. If every fiber of \(f\), viewed as a map on simplices, contains only finitely many simplices, we call \(f\) quasi-finite.
\end{defi}

In order to talk about flow, we will assume that on our graph there will be a set of allowed directions, which determine the direction in which the water flows.

\begin{defi}\label{defi:dir_graph}
 A directed graph is a graph $K=(K^0,K^1)$ plus a set of directed edges $E_K\subset K^0\times K^0$.  The set of directed edges is such that there is a well defined map $E_K\rightarrow K^1$ given by $(v,w)\mapsto \{v,w\}$. A (directed) edge $(v,w)$, also called an arc, means that from vertex $v$, one can reach vertex $w$. For a vertex $v$ we define: 
 $ \In_K(v):=\{(w,v) \in E_K\}$ and $ \Out_K(v):=\{(v,w) \in E_K\}$.
When there is no possible ambiguity about the simplicial complex $K$, we omit the subscript. 
A directed simplicial map sends each directed edge either to a directed edge with the same orientation, or collapses it to a vertex. If $\Out(v)$ and $\In(v)$ are finite for all vertices $v$, the graph is called locally finite.
\end{defi}

\begin{ex}
A directed graph $K$ induces a directed topological space $(X,dX)$ with $X=\mid K \mid$, the geometric realization of $K$, and $dX$ being the result of composing, reparametrizing, etc. the paths given by the directed edges of $K$.
\end{ex}

It is well known that the combinatorial constructions obtained with simplicial complexes have a topological counterpart, that can be obtained via their geometric realization \citep{munkres2018elements}, which we denote with the standard notation $\mid K\mid $ - $K$ being a simplicial complex. In particular, simplicial maps between simplicial complexes induce continuous maps in $\Top$. In the following we will sometimes refer to the induced continuous map $|K|\to |K'|$ in $\Top$ as the topological realization of the simplicial map.

\begin{figure}
\begin{subfigure}{\textwidth} 
    \centering
    	\includegraphics[width = 0.6\textwidth]{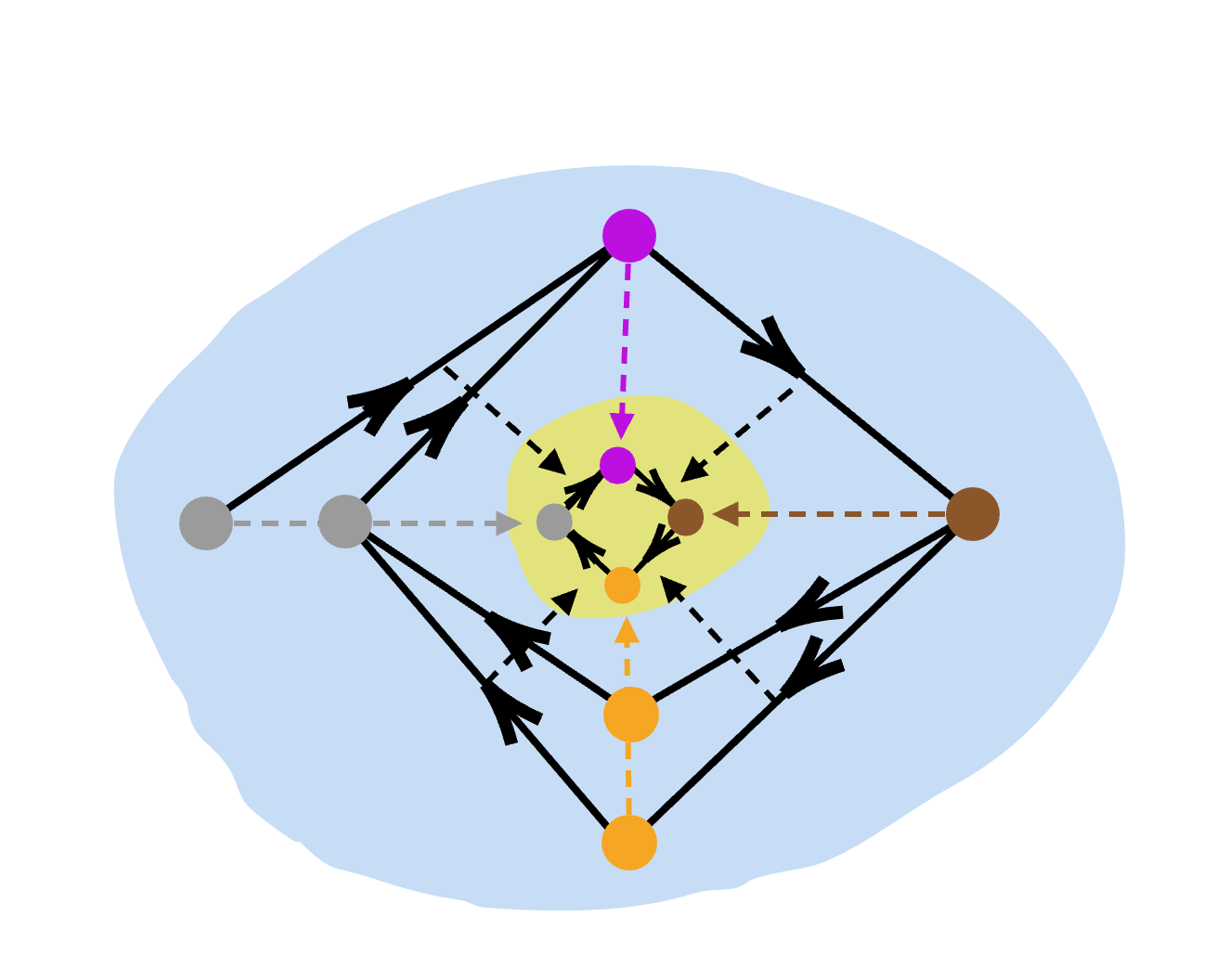}
\end{subfigure}
\caption{An example of the directed simplicial map $f:K\rightarrow Q$ which we use as starting point for our investigation: the graph on the blue background is the graph $K$, the dashed arrows represent the simplicial map $f$ on vertices and edges, and the graph on the yellow background is $Q$, whose geometric realization is homeomorphic to $\Ss$. Arrows represent arcs.}
\label{fig:reeb_map}
\end{figure}

\begin{figure}
\begin{subfigure}{\textwidth} 
    \centering
    	\includegraphics[width = 0.6\textwidth]{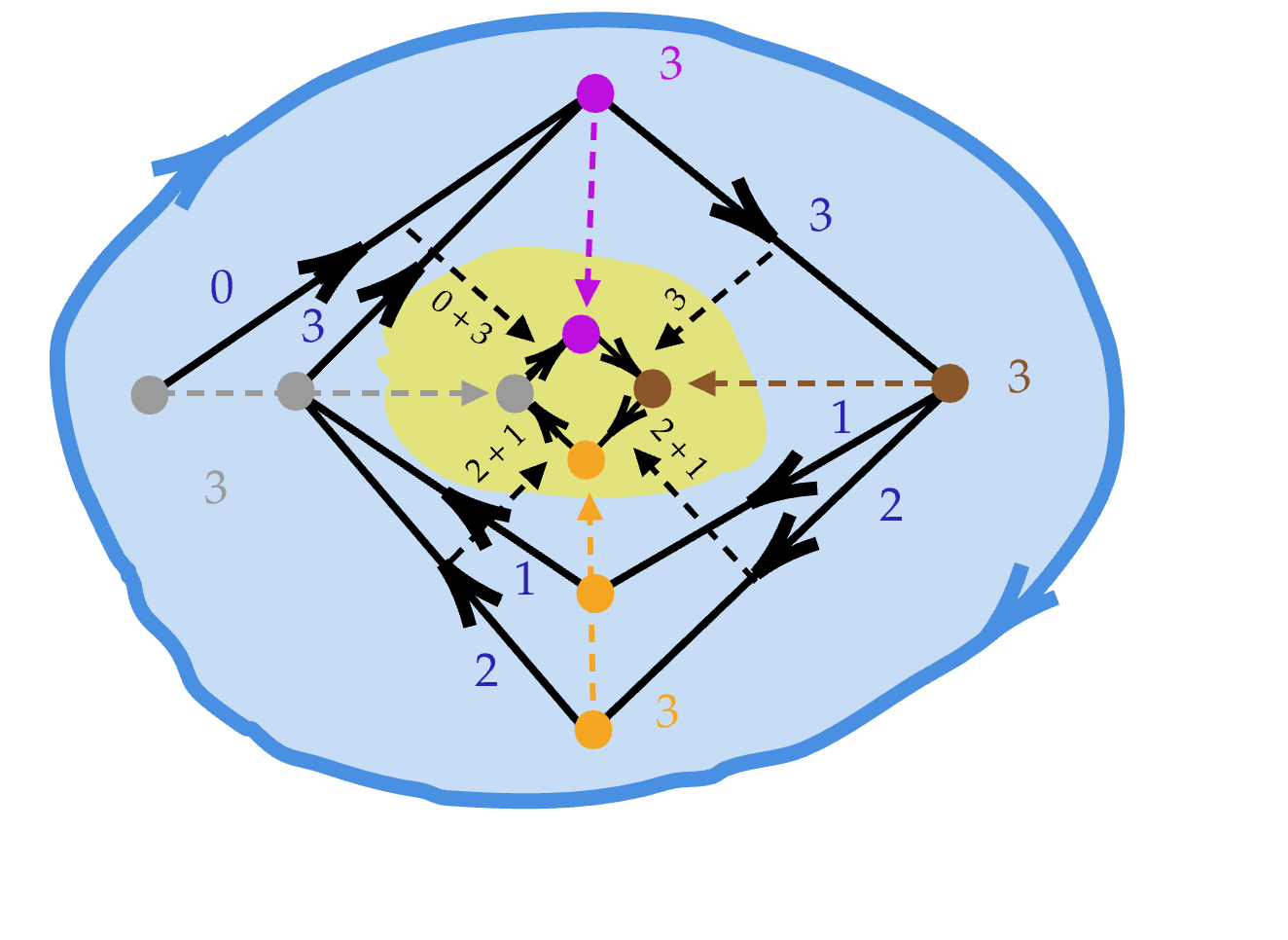}
\end{subfigure}
\caption{An example of a circulation and the map $\Fop_1(f)$, induced by the simplicial map $f$ in \Cref{fig:reeb_map}. Note that circulations on $Q$ are constant on the edges.}
\label{fig:sliced_flow}
\end{figure}

\subsection{Circulations}\label{sec:circulations}
There is a strong link between the objects used to model the flow of water in a system of pipes, called circulations, and $1$-cycles built via simplicial homology, see for instance \cite{chambers2012homology}.  
Now we model our space of circulations following the path usually taken to define simplicial homology groups, except we want our formulations to be coherent with the directed structure we have on our graph.
Also, in our definitions, we don't necessarily restrict ourselves to circulations with finite support, if the graph is not finite.

To write this down formally, we resort to functional notation, instead of usual linear combinations of edges/vertices, to better encompass the possibility of infinite graphs. Note, however, the point of view we are taking is not the one of cohomology, but rather of the linear dual of homology (at least in the case of finite graphs).

\begin{defi}
    Given a directed locally finite graph $K=(K^0,K^1)$ with $E_K$ its set of directed edges, we define the directed-chain vector spaces by
\(\DCh_0(K,\R):=\{s:K^0\rightarrow \R\}\) and
\(\DCh_1(K,\R):=\{s:E_K\rightarrow \R\}\).
We will use the following notation $s(v)=\lambda_v$, for $s\in\DCh_0(K,\R)$, and $s(e)=\lambda_e$, for $s\in\DCh_1(K,\R)$. Moreover, for an arc $e\in E_K$ we use $e^*$ to identify the characteristic function on $e$, that is $e^*(e')=1$ if $e'=e$ and $0$ otherwise, and we extend this notation linearly. Similarly for vertices. 
Then we define $\partial_1:\DCh_1(K,\R)\rightarrow \DCh_0(K,\R)$ as:
\[
\partial_1(s)(v)= \sum_{(w,v)\in \In(v)} \lambda_{(w,v)} - \sum_{(v,w)\in \Out(v)} \lambda_{(v,w)}, 
\]
which amounts to taking the discrete divergence of $s:E_K\rightarrow \R$. Note that for characteristic functions we have $\partial_1((w,v)^*)=(v-w)^*$.
We denote $\Fop_1(K,\R) = \ker(\partial_1)$. In analogy to the simplicial homology case, we may call the elements in $\ker(\partial_1)$ $1$-cycles. 
Restricting to cycles with non-negative coefficients, we obtain
\(\Fop_1^+(K,\R)=\{s\in \Fop_1(K,\R)\mid s(e)\geq 0,\ \forall e\in E_K\}\).
Elements in $\Fop_1^+(K,\R)$ are called circulations coherently with the max-flow literature \citep{chambers2012homology}.

For later use, define the total incoming and outgoing values at a vertex by
\[
\ff_{in}(s)(v):=\sum_{e\in\In(v)}s(e),
\qquad
\ff_{out}(s)(v):=\sum_{e\in\Out(v)}s(e).
\]
When \(s\) is a circulation these values agree, and we denote their common value by
\(\ff(s)(v)\).
\end{defi}

\begin{prop}[Functoriality]\label{prop:functoriality}
    A quasi-finite directed simplicial map $\alpha:K\rightarrow K'$ descends to a linear map
    \(\Fop_1(\alpha):\Fop_1(K,\R)\rightarrow \Fop_1(K',\R)\).
\end{prop}

\begin{proof}
    For every arc \(e'\) in \(K'\), define
    \(\widetilde{\alpha}_1(s)(e')=\sum_{e\in \alpha^{-1}(e')}s(e)\). Similarly, for
    every vertex \(v'\) in \(K'\), define
    \(\widetilde{\alpha}_0(s)(v')=\sum_{v\in \alpha^{-1}(v')}s(v)\). Note that, for
    these maps to be well defined, we need the quasi-finiteness of \(\alpha\). Thus we use the local finiteness of the graphs to define \(\partial_1\) and the quasi-finiteness of \(\alpha\) to define \(\widetilde{\alpha}_1\) and \(\widetilde{\alpha}_0\).

    Fix \(s\in\DCh_1(K,\R)\) and \(v'\in (K')^0\). Since \(K\) and \(K'\)
    are locally finite and \(\alpha\) is quasi-finite, all the sums below are finite. We have
    \[
    \begin{aligned}
    (\partial_1\widetilde{\alpha}_1s)(v')
    &=
    \sum_{e'\in\In(v')}\sum_{e\in\alpha^{-1}(e')}s(e)
    -
    \sum_{e'\in\Out(v')}\sum_{e\in\alpha^{-1}(e')}s(e).
    \end{aligned}
    \]
    On the other hand,
    \[
    \begin{aligned}
    (\widetilde{\alpha}_0\partial_1s)(v')
    &=
    \sum_{v\in\alpha^{-1}(v')}
    \left(
    \sum_{e\in\In(v)}s(e)-\sum_{e\in\Out(v)}s(e)
    \right).
    \end{aligned}
    \]
    Every arc collapsed by \(\alpha\) to \(v'\) occurs once with each sign in the
    last expression and therefore cancels. The remaining terms are exactly the two
    sums in the first expression. Hence
    \(\partial_1\circ\widetilde{\alpha}_1=\widetilde{\alpha}_0\circ\partial_1\).
    Thus \(\widetilde{\alpha}_1\) sends \(1\)-cycles to \(1\)-cycles, and we define
    \(\Fop_1(\alpha):=\widetilde{\alpha}_1|_{\Fop_1(K,\R)}\).
\end{proof}

\begin{cor}
    If $\alpha$ is also a rigid simplicial map, then $\Fop_1(\alpha)$ is induced by the restriction of $\alpha$ to $1$-cycles. 
\end{cor}

Indeed, rigidity prevents any arc from being collapsed, so the map constructed in
\Cref{prop:functoriality} is precisely the restriction of the simplicial map on chains.

For the readers familiar with the definition of simplicial homology, we point out that, if $K$ is a finite graph, we have the following commutative diagram:

\[
\begin{tikzcd}
\h_1(K,\R)\ar[r,""]&\h_1(K',\R)\\
\Fop_1(K,\R)\ar[u]\ar[r,""]&\Fop_1(K',\R)\ar[u].
\end{tikzcd}
\]
where the vertical arrows, in general, are neither injective nor surjective. 

In the following result, we introduce a way to lift a circulation from a directed graph to another. 

We give the definition of a covering space as it will be used throughout the manuscript.

\begin{defi}
Given a topological space $X$, a covering (space) of $X$ is given by a space $E$ and a continuous map $p:E\rightarrow X$ (also called a covering map) such that for each $x\in X$ there is an open set $U$ so that $p^{-1}(U)$ is the union of disjoint open sets each of which is mapped homeomorphically to $U$ by $p$.
\end{defi}

\begin{prop}[Lifting Circulations]\label{prop:lift}
    Let \(\alpha:K\rightarrow K'\) be a rigid directed simplicial map between locally finite directed graphs such that, for every vertex \(v\in K^0\), the induced maps
    \[
    \alpha:\In(v)\to \In(\alpha(v)),
    \qquad
    \alpha:\Out(v)\to \Out(\alpha(v))
    \]
    are bijections. Then one can lift elements in \(\Fop_1(K',\R)\) to elements in \(\Fop_1(K,\R)\) via a map
    \(\Fop_1(\alpha^*):\Fop_1(K',\R)\rightarrow \Fop_1(K,\R)\).
\end{prop}

\begin{proof}
Given \(s\in \Fop_1(K',\R)\), define \(\Fop_1(\alpha^*)(s)\in \DCh_1(K,\R)\) by
\(\Fop_1(\alpha^*)(s)(e)=s(\alpha(e))\) for every arc \(e\) in \(K\). By hypothesis,
for every vertex \(v\in K^0\), the maps
\(\alpha:\In(v)\to \In(\alpha(v))\) and
\(\alpha:\Out(v)\to\Out(\alpha(v))\) are bijections.
Hence
\begin{align*}
\partial_1(\Fop_1(\alpha^*)(s))(v)
&=
\sum_{(w,v)\in \In(v)}
\Fop_1(\alpha^*)(s)((w,v))
-
\sum_{(v,w)\in \Out(v)}
\Fop_1(\alpha^*)(s)((v,w)) \\
&=
\sum_{(\alpha(w),\alpha(v))\in \In(\alpha(v))}
s((\alpha(w),\alpha(v)))
-
\sum_{(\alpha(v),\alpha(w))\in \Out(\alpha(v))}
s((\alpha(v),\alpha(w))) \\
&=
\partial_1(s)(\alpha(v))
=
0.
\end{align*}
Therefore \(\Fop_1(\alpha^*)(s)\in \Fop_1(K,\R)\), and the result follows.
\end{proof}

Now we specialize these definitions to our scenario.
We consider $f:K\rightarrow Q$ a quasi-finite, directed simplicial map between  finite graphs where the geometric realization of $Q$ is homeomorphic to $\Ss$. We also require that $f$ is rigid, which simplifies some of the upcoming steps and matches the simplicial model constructed in \Cref{sec:reeb_pipeline}.
Moreover, we ask that the directed structure on $Q$ is interpretable w.r.t. the application discussed in the introduction, and thus with the framework detailed in \Cref{sec:reeb}.
That is, we ask that 1) arcs and edges of $Q$ are in bijection, 2)  via the homeomorphism into $\Ss$, the arcs describe either clockwise or counterclockwise circulation around $\Ss$ - as in \Cref{fig:reeb_map}. This amounts to fixing a direction around the circle in which the water is allowed to flow. 
Note that, in this context, $K$ is allowed exactly one direction per edge as, for each edge $\{v,w\}$ in $K$, only one of $(f(v),f(w))$ and $(f(w),f(v))$ can be a directed edge in $Q$.
Lastly, we consider $Q'$ such that its geometric realization is homeomorphic to $\R$ and $p:Q'\rightarrow Q$ which is a directed simplicial map such that the induced map  $\R\rightarrow \Ss$ is the universal covering of $\Ss$.

We have the following immediate results.

\begin{prop}
    With the notation previously introduced,  $\Fop_1(Q,\R)\cong \R \cong \Fop_1(Q',\R)$. In particular, $\Fop_1(p^*)$ is an isomorphism.
\end{prop}

    \begin{proof}
        By construction for any vertex $v$ in $Q$ we have $\In(v)=\{e_v^+\}$ and $\Out(v)=\{e_v^-\}$.
        Thus, for any $s\in \Fop_1(Q,\R)$, 
        $\partial_1(s)(v)= \lambda_{e_v^+}-\lambda_{e_v^-}=0$. Thus $s$ is constant on the arcs of $Q$. The proof for $Q'$ is analogous. Lastly, by construction, we have $\Fop_1(p^*)(s)(e)=s(p(e))$; thus the constant value which identifies $s\in \Fop_1(Q,\R)$ is the same which identifies $\Fop_1(p^*)(s)\in \Fop_1(Q',\R)$. 
    \end{proof}

Since we have $Q\cong \Ss$ and the isomorphism $\Fop_1(Q,\R)\cong \R$ is given by $s\in \Fop_1(Q,\R)$ being constant on the edges of $Q$, we obtain the following corollary which has an important interpretation.

\begin{cor}\label{cor:fix_iso}
Given  $f:K\rightarrow Q$ as above, for every $e,e'$ arcs of $Q$ and for every $s\in\Fop_1(K,\R)$, we have:
\[
\textstyle\Fop_1(f)(s)(e)=\Fop_1(f)(s)(e')=\sum_{(v,w)\in f^{-1}(e')}s((v,w)).
\]

 As a consequence $\Fop_1(f)(s)$ is a constant function and we call $\Fop_1(f)(s)$ the flow of $s$.   
\end{cor}

This follows because \(\Fop_1(f)(s)\) is a circulation on \(Q\), and the preceding
proposition shows that every circulation on \(Q\) is constant on its arcs.

\begin{rmk}
    \Cref{cor:fix_iso} states that for every $s\in\Fop_1(K,\R)$ the map $\Fop_1(f)$ applied on $s$ captures the amount of mass that $s$ transfers between the fibers of $Q$, as in \Cref{fig:sliced_flow}. Note that this also implies that we can have non-zero flow only if $f: K\rightarrow Q$ is surjective.  
\end{rmk}

Further restricting ourselves to non-negative coefficients, we obtain the non-negative valued function
\(\Fop_1^+(f):\Fop_1^+(K,\R)\rightarrow \Fop_1^+(Q,\R)\cong \R_{\geq 0}\).

\begin{figure}
    \centering
    	\includegraphics[width = \textwidth]{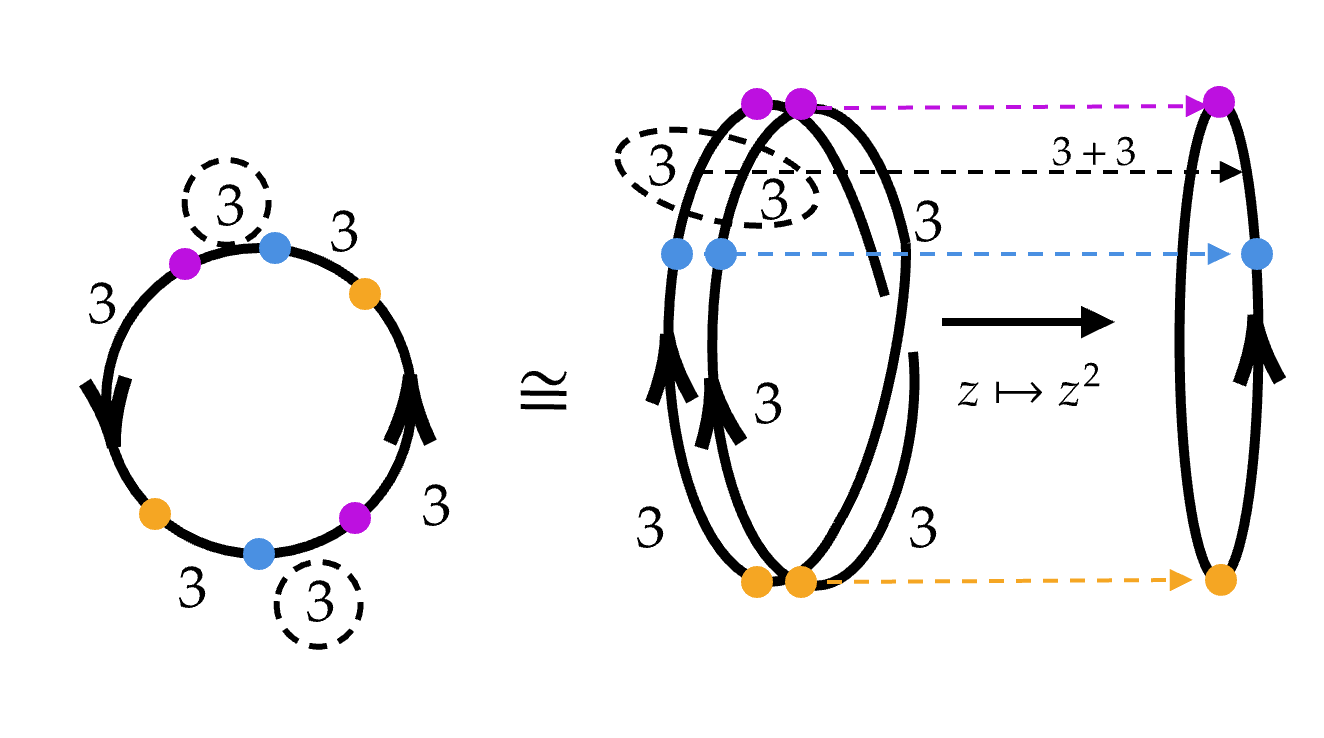}
\caption{An example of circular max-flow. Consider the map of complex numbers $z\mapsto z^2$ restricted to $\Ss$. Discretizing this map appropriately, yields the result showcased in the figure. If we consider the circulation $s$ which assigns $3$ to every edge, this would give a circular flow of $6$ as we have fibers with cardinality $2$ - as shown by the dashed arrows indicating the induced simplicial map. However, the (geometric realization of the) graph is homeomorphic to $\Ss$ and the \virgolette{flow} in the only \virgolette{pipe} available is $3$. }
\label{fig:pullback_flow}
\end{figure}

\subsection{Capacity Constraints and Max-Flow}
\label{sec:max_flow}

The last piece of the puzzle we need to introduce in this section, is a set of maximal capacities on the edges: a function $C:E_K\rightarrow \R_{\geq 0}$ that tells us how much water can flow through a particular arc. Note that, given the identification between arcs and edges in our current setting, we may also write $C:K^1\rightarrow \R_{\geq 0}$, as we may also use interchangeably the words \virgolette{edge} and \virgolette{arc}.

Given a capacity function \(C\) on the arcs of \(K\), we write \(s\leq C\) to mean that
\(s(e)\leq C(e)\) for every \(e\in E_K\). With a slight abuse of notation, we may
also write \(s\leq C\) for the set
\(\{s\in \Fop_1^+(K,\R)\mid s(e)\leq C(e),\ \forall e\in E_K\}\).

We define the (circular) max-flow of $f:K\rightarrow Q$ with capacities $C$ as:
\begin{equation}\label{eq:max_flow}
\textstyle\f(f,C) = \max_{s\leq C} \Fop_1^+(f)(s).    
\end{equation}

As already mentioned, we can think of $\f(f,C)$ as some kind of sectional flow, describing the flow through the slices given by $f^{-1}(e)$, as shown in \Cref{fig:sliced_flow}.

In the following we may need to change the coefficients of circulations and
whenever a capacity-constrained problem is restricted to integer or rational
circulations, we also assume that its capacities are integer-valued or
rational-valued, respectively. With this convention, we introduce the notation:
\begin{equation}
\textstyle\f_\Z(f,C) = \max_{s\in \Fop_1^+(K,\Z), s\leq C} \Fop_1^+(f)(s),   
\end{equation}
and analogously for $\Q$. If the coefficients are not specified, it is implied that 
we are considering $\R$. Note that, by construction, considering $s:E_K\rightarrow \Z$ is in fact equivalent to considering $s:E_K\rightarrow \N$.

\begin{rmk}[Min-Cost Circulations]\label{rmk:circul}
As already pointed out in the abstract and in the introduction, the circular max-flow problem can be framed as a particular case of min-cost circulation problems as defined in \cite{schrijver2003combinatorial}. 
In particular, using the notation of \cite[Ch.~12]{schrijver2003combinatorial}, fix an arc
$e$ of $Q$ and define $k_e:E_K\rightarrow\Q$ by $k_e(e')=-1$ if $f(e')=e$
and $k_e(e')=0$ otherwise. With zero demand and capacity function
$C:E_K\rightarrow\Q_{\geq 0}$, minimizing the total cost is equivalent to maximizing
$\f_\Q(f,C)$; the optimal cost is the negative of the circular max-flow value. The
Cycle Cancelling Algorithm \citep{klein1967primal, tardos1985strongly, schrijver2003combinatorial}
can therefore be used to compute $\f_\Q(f,C)$. If $C$ is integer-valued, the algorithm
produces an optimal circulation with integer coordinates.
We stress once again that, even though our formulation can be obtained via previously known results, our topological point of view unveils a geometric interpretation of the associated minimum-cost circulation problem which is important for the data-analysis interpretation. Plus, it is essential to obtain the results in the rest of the manuscript. 
\end{rmk}

\subsection{Computing Circular max-flow}\label{sec:optimization}

Now we devise an optimization procedure to compute $\f(f,C)$ resorting to linear programming. We propose this alternative approach to max-flow computation as it is easy to implement with any linear solver and further emphasizes some key differences with classical source-target max-flow.
We now give a computational formulation for the finite graphs considered above.

\begin{figure}
\begin{subfigure}{\textwidth} 
    \centering
    	\includegraphics[width = 0.6\textwidth]{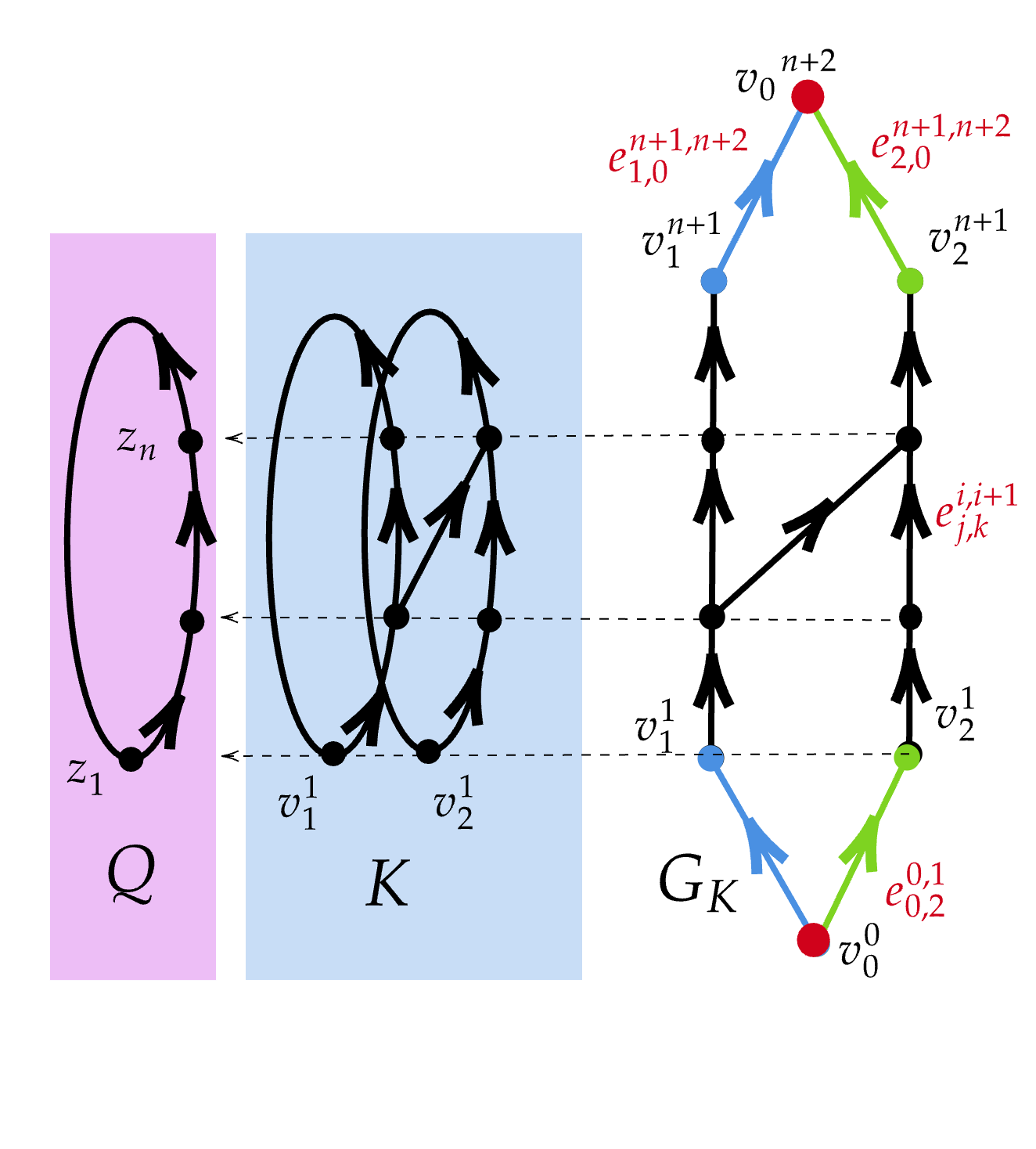}
\end{subfigure}
\caption{Visualization of the augmented graph $G_K$ employed in \Cref{sec:optimization}. We see the external source $v^0_0$ and target $v_0^{n+2}$ and the additional arcs. Light blue and green colors signify the constraints applied on those added arcs: same colors indicate equality constraints on the corresponding flow variables. In red we report some of the real valued variables associated to the arcs in $G_K$ to compute the circular max-flow.}
\label{fig:augmented}
\end{figure}

The main idea is, starting from $f:K\rightarrow Q$,
to \virgolette{open} $K$, obtaining a graph $G_K$ with an external source and an external target, and then solve a source-target max-flow with additional constraints.
\Cref{fig:augmented} may be used for visualizing the upcoming pieces of notation.

First we fix a parametrization of the vertices of $Q$ and then we use it to parametrize the simplices in $K$.
We fix a vertex $z_1\in Q^0$ and then order all the others so that $(z_i,z_{i+1})$ is an arc in $Q$.
 Now we consider $K$ and we name its vertices in the following way: $f^{-1}(z_i)=\{v^i_1,\ldots,v^i_{m_i}\}$. Note that since $f$ is rigid, we cannot have edges between fibers of two non consecutive $z_i$: only edges of the form $\{v^i_j, v^{i+1}_k\}$ or $\{v^n_j, v^{1}_k\}$ are allowed.

Starting from $K$ we \virgolette{open it up} and obtain $G_K$ with two steps. First:
\begin{enumerate}
    \item add a vertex $v^{n+1}_j$ for each vertex $v^1_j$;
    \item for every arc \((v_j^n,v_k^1)\in \Out_K(v_j^n)\), add the arc
\((v_j^n,v_k^{n+1})\) and remove the arc \((v_j^n,v_k^1)\),
\end{enumerate}
 
and then we add an external source and target as follows:

\begin{enumerate}\setcounter{enumi}{2}
    \item we add a vertex $v^0_0$ and a vertex $v^{n+2}_0$;
    \item and add the arcs $(v^0_0,v_j^1)$ and $(v^{n+1}_j,v^{n+2}_0)$ for all $j=1,\ldots,m_1$.
\end{enumerate}

\begin{figure}
\begin{subfigure}{\textwidth} 
        \centering
    	\includegraphics[width = 0.7\textwidth]{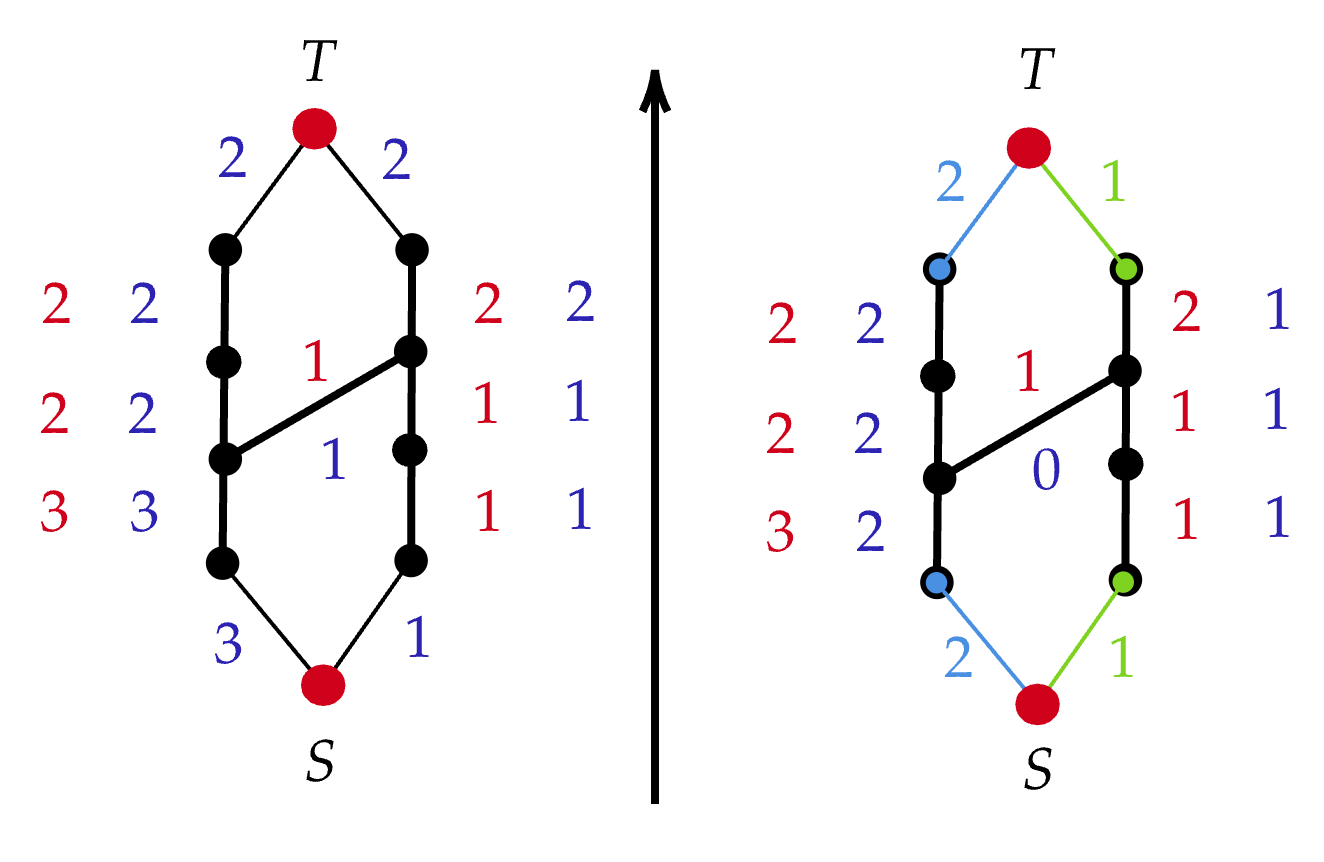}
\end{subfigure}
\caption{An example to compare max-flow computed with the usual source-target approach (left) and with the procedure described in \Cref{sec:optimization}, to ensure the periodicity of the solution. Red values represent capacity constraints and blue values represent circulation coefficients. The additional constraints are displayed, on the right, by arcs of the same color (light blue and green). The coefficient/flow of the arcs leaving the source S, must correspond to the flow of the associated arcs coming into the target T.}
\label{fig:comparison_maxflow}
\end{figure}

Now we bridge between $G_K$ and the linear optimization problem to compute circular max-flow.
For all $i=0,\ldots,n+1$, a variable $e^{i,i+1}_{j,k}$ is associated to each arc $(v^i_j,v^{i+1}_k)\in G_K$ describing the flow on that arc. Then we consider the following constraints:
\begin{enumerate}
    \item (mass conservation constraints) for each $v_j^i$ vertex in $G_K$ with $i=1,\ldots,n+1$:
    \begin{equation}\label{eq:constr_1}
        \sum_{(v^i_j,v^{i+1}_k)\in \Out(v_j^i)} e^{i,i+1}_{j,k} = \sum_{(v^{i-1}_k,v^i_j)\in \In(v_j^i)} e^{i-1,i}_{k,j};
    \end{equation}
    \item (periodicity constraints) for all $j=1,\ldots,m_1$, we ask:
    \begin{equation}\label{eq:constr_2}
       e^{0,1}_{0,j}=e^{n+1,n+2}_{j,0}=\sum_{(v^{n}_k,v^{n+1}_j)\in\In(v^{n+1}_j)} e^{n,n+1}_{k,j}. 
    \end{equation}
    Where the last equality is given by the first constraint applied to $v_j^{n+1}$. 
\end{enumerate}

\begin{figure}
        \centering
    	\includegraphics[width = 0.65\textwidth]{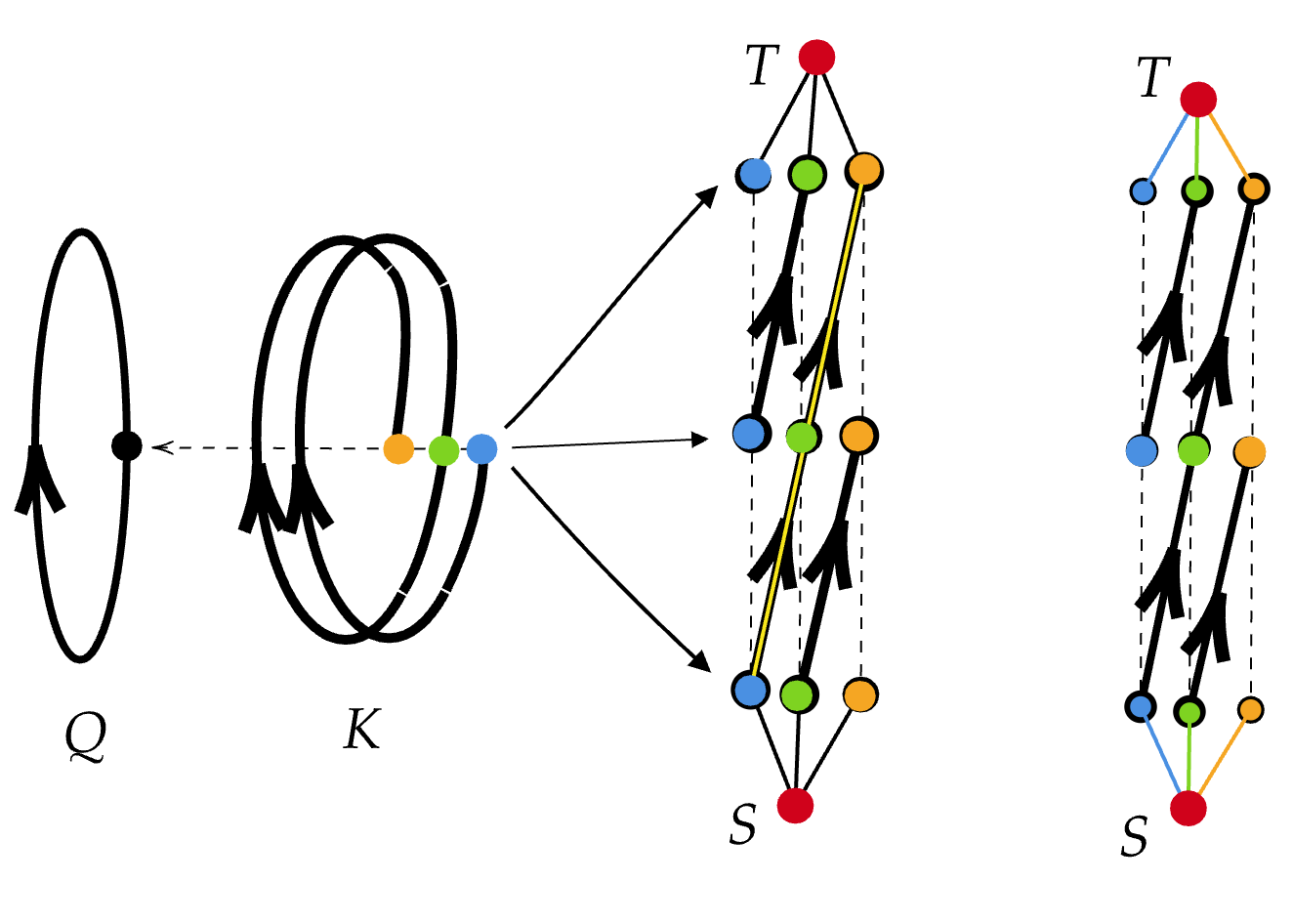}
\caption{An example to compare source-target max-flow and circular max-flow. Consider the map $f:K\rightarrow Q$ on the left. To minimize the number of vertices and edges involved, we use a graph $Q$ which is not a simplicial complex, but the phenomenon described here can clearly be replicated using a subdivision of the graphs involved making everything \virgolette{simplicial}.  \virgolette{Open} the graph $K$ cutting it along the fiber of the only vertex of $Q$ and stack two copies of it before adding the external source and target. And then compute the source-target max-flow (centre). We see a yellow path that goes from source to target ensuring a positive flow. On the rightmost graph, instead, we compute the circular max-flow, adding the additional constraints depicted by the colours of the arcs leaving the source and entering the target. As also clearly visible in the starting graph, the only admissible solution is $0$ and the yellow path is no more admissible: everything entering the orange vertex at the top, must leave the orange vertex at the bottom.}
\label{fig:comparison_stacking}
\end{figure}

\begin{rmk}\label{rmk:constraints}
    The constraints 2. differentiate our definition from computational schemes for classical source-target max-flow in $G_K$. Look at \Cref{fig:comparison_maxflow} and \Cref{fig:comparison_stacking} for example. Note that one can always roll up the graph $G_K$ by identifying $S$ with $T$, but this is different, in general, from identifying $v_i^1$ with $v_i^{n+1}$ for all $i$.
\end{rmk}

We denote by $\Adm(K)$ the set of all admissible values of the variables $e^{i,i+1}_{j,k}$ satisfying these constraints. We will show that $\Adm(K)$ coincides exactly with the set of circulations on $K$.

\begin{prop}\label{prop:equivalence}
  We can build correspondences $s\mapsto e_s$ and $e\mapsto s_e$ giving inverse bijections between $\Fop_1(K,\R)$ and $\Adm(K)$. 
    Moreover, $\Fop_1(f)(s_e)=\sum_{(v^0_0,v^{1}_j)\in \Out_{G_K}(v_0^0)} e^{0,1}_{0,j}$.
\end{prop}

The proof, given in Appendix~\ref{app:circular_flow_proofs}, constructs both
correspondences explicitly and checks the conservation and periodicity constraints.

Thanks to \Cref{prop:equivalence}, given a capacity function $C$, solving
\(\max_{s\leq C}\Fop_1^+(f)(s)\) is equivalent to solving
\begin{equation}\label{eq:opt_problem}
\max_{e\in \Adm(K)}\sum_{(v^0_0,v^{1}_j)\in \Out_{G_K}(v_0^0)} e^{0,1}_{0,j}    
\end{equation}
subject to the linear inequalities granting $0\leq s_e\leq C$ i.e. the corresponding arc variables in \(G_K\) are bounded by the capacity of their rolled-up arc in \(K\).  
Note that we obtain a linear optimization problem with linear constraints.

This computational procedure can also be used to further highlight the differences between the circular flow we defined and usual source-target flows:  \Cref{fig:comparison_maxflow}  shows one example in which  constraints $2.$ are active and fundamental to obtain a solution that can be \virgolette{rolled} back to $K$. As already stressed by \Cref{rmk:constraints}. Moreover, in \Cref{fig:comparison_stacking} we report an example of two stacked copies of $K$ in which we have positive source-target flow but no circular flow.

\subsection{Max-Flow on the Periodic Graph}\label{sec:int_flow}

In this section we study what happens to the circular max-flow of the map $K\rightarrow Q$ when we \virgolette{unroll} it and consider infinitely stacked copies of $K$, obtaining a periodic graph $K'$, with a map $K'\rightarrow Q'$, $\mid Q'\mid \cong \R$.

By \Cref{lem:pull_back}, this stacking is the graph-level counterpart of unrolling a space along one periodic direction and then taking its Reeb representation.

\subsubsection{Pullbacks of Graphs}

We now give a combinatorial description of the graph $K'$ which describes the unrolling of $K$. Topologically, as we will see in later sections, this is the pullback of the geometric realization of $f:K\rightarrow Q$ along the universal covering of $\Ss$. 
However, for the sake of this section, we don't need such general topological approach.

\begin{figure}
    \centering
    	\includegraphics[width = 0.7\textwidth]{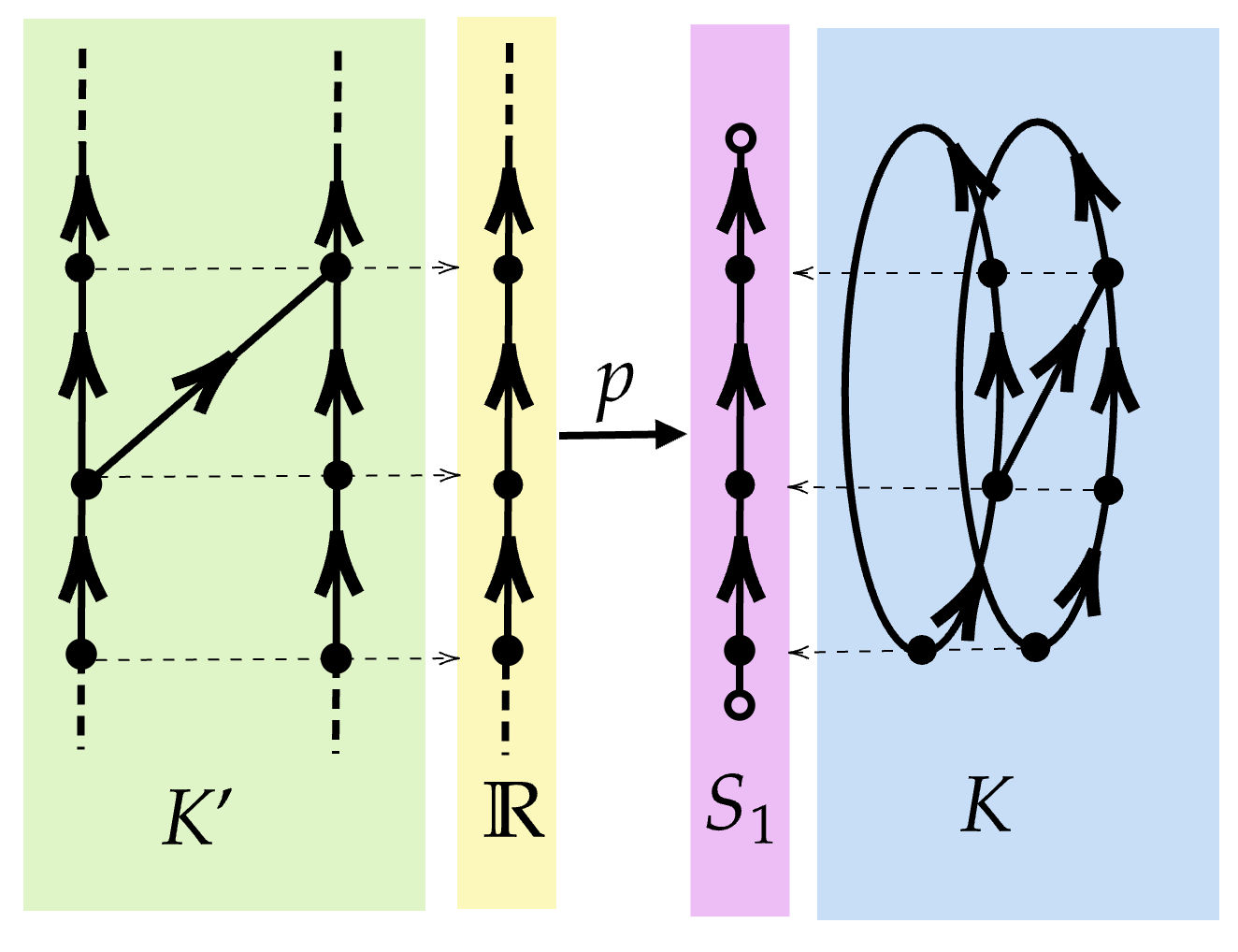}
\caption{Example of the pullback $K'\rightarrow Q'$ of a map $K\rightarrow Q$ via the universal covering of $\Ss$.}
\label{fig:pullback}
\end{figure}

\begin{figure}
    \centering
    	\includegraphics[width = 0.4\textwidth]{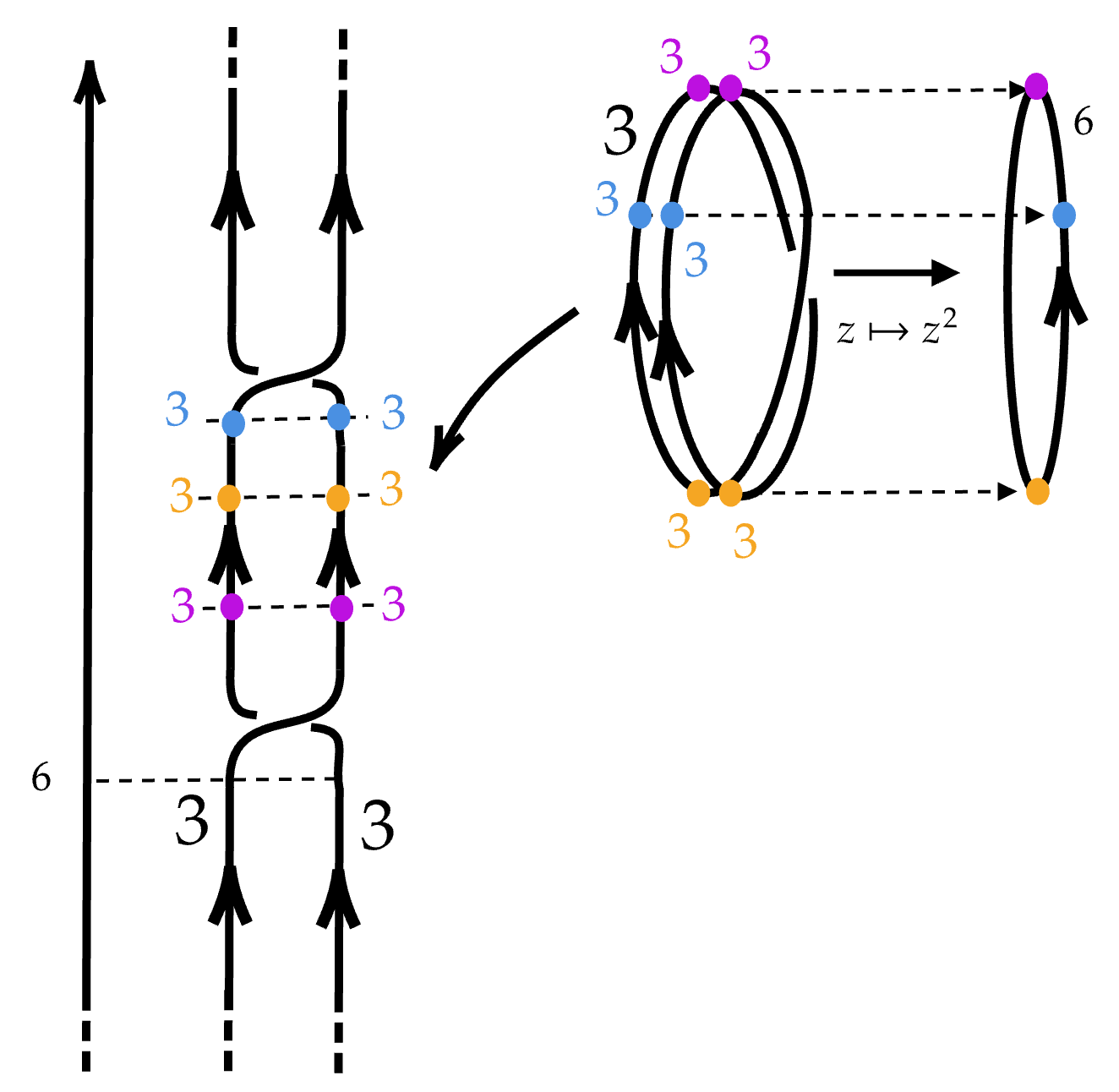}
\caption{Another example of the pullback of $K\rightarrow Q$, involving again the map of complex variable $z\mapsto z^2$ restricted to $\Ss$. We also lift a circulation to $K'$. As already pointed out in \Cref{fig:pullback_flow}, there is only one pipe available in $K$. However, if we take the pullback and \virgolette{unroll} the graph, we see two separate pipes appearing, each with a flow of $3$.}
\label{fig:pullback_circ}
\end{figure}

Consider the maps $f:K\rightarrow Q$ and $p:Q'\rightarrow Q$ as before. 
We now build the (unique, up to isomorphism) graph $K'$ with directed simplicial maps $p':K'\rightarrow K$ and $f':K'\rightarrow Q'$. A visual interpretation of the upcoming steps can be found in \Cref{fig:pullback} and \Cref{fig:pullback_circ}.

First we parametrize the points of $Q'$ coherently with the ones on $Q$. We recall that $Q$ is made by vertices $\{z_i\}_{i=1\ldots,n}$, and $(z_i,z_{i+1})$ and $(z_n,z_1)$ are the arcs in $Q$.
Consider now the map $Q\cong \Ss$ and in particular, $Q\setminus\{z_1\}\cong \Ss\setminus\{x\}$ for some $x\in\Ss$. 
The preimage of $\{x\}$ via the universal cover $\R\rightarrow \Ss$ can clearly be given a total order, since it is a subset of $\R$. Thus we can fix an order preserving bijection between such a fiber and $\Z$. In other words, we can write $p^{-1}(z_1)= \{z_1^r\}_{r\in \Z}$.
Thus, for any other vertex $z_i$ we can parametrize the fiber of $p$ so that  $p^{-1}(z_i)= \{z_i^r\}_{r\in \Z}$ and $z_1^r<z_i^r <z_1^{r+1}$.

We have parametrized the vertices of $K$ so  that $f^{-1}(z_i)=\{v^i_1,\ldots,v^i_{m_i}\}$.
We build the simplicial complex $K'$ by considering vertices  $\{v^{i,r}_1,\ldots,v^{i,r}_{m_i}\}$ for every $i=1,\ldots,n$ and for every $r\in\Z$. Then we consider arcs $(v^{i,r}_k,v^{i+1,r}_j)$ for every arc $(v^{i}_k,v^{i+1}_j)$ in $E_K$ and 
$(v^{n,r}_k,v^{1,r+1}_j)$ for every arc $(v^{n}_k,v^{1}_j)$ in $E_K$. 
Lastly we define the simplicial maps on vertices, and hence induced on edges (and arcs): $f':K'\rightarrow Q'$ and $p':K'\rightarrow K$ as $f'(v^{i,r}_j)=z_i^r$ and $p'(v^{i,r}_j)=v_j^i$. 
In the following we refer to the fibers $f^{-1}(z_i)$ and $(f')^{-1}(z_i^r)$ as layers, while we refer to the preimage in $K'$ of the \virgolette{segment} going from $z_1^r$ to  $z_n^r$, as a deck (see \cite{hatcher}). Roughly speaking, with this notation, layers are indexed by $i$, and decks by $r$.

As already mentioned, one can show that the (geometric realization of) the following diagram is a pullback in the category of topological spaces and, in particular that (the geometric realization of) the map $p'$ is a covering map: 

\[
\begin{tikzcd}
K' \ar[r,"f'"]\ar[d,"p'"]&Q'\ar[d,"p"]\\
K\ar[r,"f"]&Q.
\end{tikzcd}
\]

Moreover, since the maps involved in the previous diagram are all directed simplicial maps, we have the following commutative diagrams of circulations:

\begin{equation}\label{eq:unroll}
\begin{tikzcd}
\Fop_1(K,\R)\ar[r,"\Fop_1(f)"]\ar[d,"\Fop_1(p'^*)"]&\Fop_1(Q,\R)\ar[d,"\Fop_1(p^*)"]\ar[r,"\cong"]&\R\ar[d,"\cong"]\\
\Fop_1(K',\R)\ar[r,"\Fop_1(f')"]&\Fop_1(Q',\R)\ar[r,"\cong"]&\R.
\end{tikzcd}
\end{equation}

If we define the flow of $s\in\Fop_1(K',\R)$ as the value $\Fop_1(f')(s)$, analogously to what was already done for $\Fop_1(K,\R)$, the above diagram implies that lifting circulations along $\Fop_1(p'^*)$  preserves flow.
Similarly, given $C':E_{K'}\rightarrow \R_{\ge 0}$, we define the max-flow constrained by $C'$ as
\(\f(f',C'):=\sup_{s\leq C'} \Fop_1^+(f',\R)(s)\).

We may refer to $K'$ as the unrolled graph or the periodic graph, in contrast with $K$, which is the circular graph.

\begin{defi}
We say that a function $g':E_{K'}\rightarrow \R$ is periodic if
there is $g:E_K\rightarrow \R$ such that 
$g'=g \circ p'$.
In other words, for every arc $(v_j^{i,0},v_k^{i+1,0})$ and every $r\in\Z$, we have
\(g'((v_j^{i,0},v_k^{i+1,0}))=g'((v_j^{i,r},v_k^{i+1,r}))\).
Note that, for every $g:E_K\rightarrow \R$, we can define $g'=g \circ p'$. We may also write $g'=p'^*g$. 
\end{defi}

Note that \Cref{eq:unroll} implies that $\f(f',p'^*C)\geq \f(f,C)$.

We now introduce the deck transformations (see \cite{hatcher}) $T_u:K'\rightarrow K'$ which we use to \virgolette{shift the indexes} of the graph $K'$ upwards or downwards:
\begin{align*}
&T_u(v_j^{i,r}):=
v_j^{i,r+u},\\
&T_u((v_j^{i,r},v_k^{i+1,r})):=
(v_j^{i,r+u},v_k^{i+1,r+u}),\\
& T_u((v_j^{n,r},v_k^{1,r+1}))
:= (v_j^{n,r+u},v_k^{1,r+1+u}),
\end{align*}
with $u\in \Z$.
Then we can define an action of $\Z$ on $\DCh_1(K',\R)$ via $(u.q):=q\circ T_u$, for $q\in\DCh_1(K',\R)$. Similarly, $\Z$ acts on the subsets $\gamma\subset E_{K'}$ via $u.\gamma=T_u(\gamma)$. Given $\gamma\subset E_{K'}$, the union of the orbits of its edges is $p'^{-1}(p'(\gamma))$.

\begin{lem}
Given a periodic function $C':E_{K'}\rightarrow \R_{\geq 0}$ then $\Z$ acts on $\Fop_1^+(K',\R)$. And $q\leq C'$ iff $u.q\leq C'$ for any $u\in \Z$. 
\end{lem}

\begin{proof}

The first claim follows from the fact that for any vertex $v_j^{i,r}$ and every $u\in \Z$, we have $T_u(\In(v_j^{i,r})) = \In(v_j^{i,r+u})$ and $T_u(\Out(v_j^{i,r}))=\Out(v_j^{i,r+u})$.

The last part of the statement holds because of the following equalities. 
Let $q\in \Fop_1^+(K',\R)$ with $q\leq C'$. Then:
\[
u.q((v_j^{i,r},v_k^{i+1,r})) = q((v_j^{i,r+u},v_k^{i+1,r+u})) \leq C'((v_j^{i,r+u},v_k^{i+1,r+u}))=
C'((v_j^{i,r},v_k^{i+1,r})).
\]
The same computation applies to wrap-around arcs $(v^{n,r}_j,v^{1,r+1}_k)$.\end{proof}

\subsubsection{Unrolling Circulations}

Now we prove the main result of this section, which states that
\(\f(f,C)=\f(f',p'^*C)\).
The main obstacle to this would be if some circulations on the unrolled graph did not have corresponding circulations on the original graph.
First we need to setup some additional notation. 
Let $E_i$ be the set of arcs of $K$ that go from layer $i$ to layer $i+1$. Similarly, $E_n$ collects the ones that go from layer $n$ to layer $1$.
Let $E_i^r$ be the set of arcs of $K'$ that go from layer $i$ in deck $r$, to layer $i+1$ in deck $r$. And, lastly, $E^r_n$ collects the arcs that go from layer $n$ in deck $r$, to layer $1$ in deck $r+1$.
Moreover, we recall some classical max-flow definitions.

\begin{defi}
Given a finite directed graph \(G\), let \(v_S\) be a source vertex and \(v_T\) a
target vertex. An ST-circulation is a map \(s:E_G\to\R_{\geq 0}\) such that
\(\partial_1s(v)=0\) for every \(v\in G^0\setminus\{v_S,v_T\}\).

Its ST-flow is
\[
\sum_{e\in\Out(v_S)}s(e)-\sum_{e\in\In(v_S)}s(e)
=
\sum_{e\in\In(v_T)}s(e)-\sum_{e\in\Out(v_T)}s(e).
\]
In the graphs considered below, no arc enters \(v_S\) and no arc leaves \(v_T\), so
this reduces to \(\sum_{e\in\Out(v_S)}s(e)=\sum_{e\in\In(v_T)}s(e)\).
Given capacities \(C:E_G\to\R_{\geq 0}\), an ST-circulation is feasible if
\(s(e)\leq C(e)\) for every \(e\in E_G\). We compute source-target max-flow
$\text{ST-}\f(C)$ by maximizing the ST-flow over feasible ST-circulations. As for
circular flow, we may consider $\text{ST-}\f_\Z(C)$ and $\text{ST-}\f_\Q(C)$ by
changing the coefficient domain, with capacities in the same domain.
\end{defi}

To tackle the main problem of this section we need a technical lemma and a preliminary theorem: 
the lemma guarantees that in any circulation on $K'$ we can find a \virgolette{periodic} structure; while the theorem states that we can map integer circulations from $K'$ to rational circulations on $K$ preserving flow and controlling uniformly the denominators which appear.

\begin{lem}\label{lem:finiteness}
    Fix \(L\in\mathbb N\). There exists \(H_L\in\mathbb N\) such that, for every
\(s:E_{K'}\to\mathbb N\) satisfying
\[
\sum_{e\in E_i^r}s(e)\le L
\]
for every \(i\) and \(r\), for each layer \(i\) and every \(a\in\mathbb Z\), there exist
two decks \(r<r'\) in \(\{a,\ldots,a+H_L\}\) such that
\[
s(e)=s(T_{(r'-r)}e)
\qquad\text{for every }e\in E_i^r.
\]
\end{lem}

  \begin{proof}
        For $l\in \N$, set $(a_i^r)_l=\{g:\{1,\ldots,\# E_i^r\}\rightarrow \{0,1,\ldots,l\}\mid \sum_{k=1}^{\# E_i^r} g(k)=l\}$. 

Define $A^r_
l:=\bigcup_{i=1}^n (a_i^r)_l$ for any choice of $r\in \Z$.
Since $\# E_i^r = \# E_i^{r'}$, then $(a_i^r)_l = (a_i^{r'})_l$ for every $r'\in \Z$.
Thus, we can remove the dependency on the deck $r$ and simply consider $A_l:=\bigcup_{r} A^r_l$.  

We know  that $A_l$ is a finite set. In fact, $(a_i^r)_l$ are finite sets (they are a subset of the set of functions between two finite sets).

Now we set $A_L := \bigcup_{l=0,\ldots,L} A_l$.
    Since $A_L$ is a finite union of finite sets, it is a finite set. 

Consider $s:E_{K'}\rightarrow \N$ with $\sum_{e\in E_i^r} s(e) \leq L$ for every $i$ and $r$, and fix an ordering inside every set $E_i^r$, coherently across different decks.
For any fixed $i$ we can build a function $\varphi_i:\Z\rightarrow A_L$ defined as $r\mapsto (s(e^r_1),\ldots,s(e^r_{\#E_i^r}))$, with a vector $(a_1,\ldots,a_{\# E_i^r})$ being identified with the function $j\mapsto a_j$.
Such a map considers the i-th layer inside each deck, and sends each deck $r$ into a vector which partitions the value $l = \sum_{e\in E_i^r} s(e)$ into the different arcs, as defined by $s$.

Since $A_L$ is finite, this function cannot be injective. In other words, the value of $\varphi_i$ on at least two different decks must coincide. 
More precisely, for every $A\subset\Z$, with  $\#A >\# A_L$, we know that there exist two decks $r,r'\in A$ such that $\varphi_i(r)=\varphi_i(r')$, which means precisely that
\[
s(e)=s(T_{(r'-r)}e)
\qquad\text{for every }e\in E_i^r.
\]

Assume \(r<r'\) and take \(H_L\ge \#A_L\). For any \(a\in\mathbb Z\), the set
\(A=\{a,\ldots,a+H_L\}\) has cardinality \(H_L+1>\#A_L\), and the previous
paragraph gives the desired decks \(r<r'\) in \(\{a,\ldots,a+H_L\}\).
    \end{proof}

\begin{teo}[Rolling Up Circulations]\label{teo:unrolled_flow}
For every $C:E_K\rightarrow \N$, there exist a natural number $H$ and a flow-preserving operator 
\[
\textstyle\varphi_H:\{q\in\Fop_1^+(K',\Z)\mid q\leq p'^*C\} \rightarrow \{s\in\Fop_1^+(K,\Q)\mid s\leq C\},
\]
such that $H\cdot \varphi_H(q)\in \Fop_1^+(K,\Z)$.
\end{teo}

The proof is given in Appendix~\ref{app:circular_flow_proofs}. It averages the decks
between two repeated layer profiles, producing a rational circulation on \(K\) with
a controlled denominator.

Now we are ready to prove the equivalence between the flows.

\begin{teo}[Equivalence of Flows]\label{teo:integer}
    For every $C:E_K\rightarrow \Z_{\ge 0}$ we have:
    \begin{itemize}
        \item $\f(f,C)=\f(f',p'^*C)$;
        \item there exist a finite subgraph $G$ of $K'$ such that, upon adding an external source and an external target, $\text{ST-}\f(p'^*C_{\mid G})=\f(f,C)$.
    \end{itemize}
    In particular, if we consider integer capacities, there always exists a solution with integer coordinates  both for $\f(f,C)$ and $\f(f',p'^*C)$.
\end{teo}

\begin{proof}
        Take $L = \max_i \sum_{e\in E_i} C(e)$, consider $H_L$ as in \Cref{lem:finiteness}, and define $C':=p'^*C$.
    Let $K_L$ be the sub-graph of $K'$ given by $(f')^{-1}([z_1^0,z_1^{H_L+1}])$. From $K_L$ we obtain $G_{K_L}$ adding an external source $S=\{v^{0,0}_1\}$ on the \virgolette{bottom}, connected to $\{v_1^{1,0},\ldots,v_{m_1}^{1,0}\}$ via arcs $(v^{0,0}_1,v_j^{1,0})$, and an external target $T=\{v_1^{0,H_L +1}\}$ on \virgolette{top}, connected to  $\{v_1^{1,H_L +1},\ldots,v_{m_1}^{1,H_L +1}\}$ via arcs $(v_j^{1,H_L +1},v_1^{0,H_L +1})$.
    Similarly to what was done in \Cref{sec:optimization} to obtain $G_K$. We collect the bottom arcs in $E_0^0$ and the top ones in $E_0^{H_L+1}$.

We work with the graph $G_{K_L}$ using the usual source-target max-flow. In particular given $s:E_{G_{K_L}}\rightarrow \R$ ST-circulation, the ST-flow is equal to the total flow through every internal cut
\(E_i^r\), with \(i=1,\ldots,n\) and \(r=0,\ldots,H_L\):
\[
\sum_{e\in E_0^{0}} s(e)
=
\sum_{e\in E_i^r} s(e)
=
\sum_{e\in E_0^{H_L+1}} s(e).
\]

Given $q\in \Fop_1^+(K',\R)$, we can always restrict $q$ on $K_L$
and then 
induce a circulation $\hat{q}$ on $G_{K_L}$ with the same flow as $q$.
This, in particular, shows $\text{ST-}\f(C'_{\mid K_L})\geq \f(f',C')$, where by $C'_{\mid K_L}$ we mean any set of capacity constraints on $G_{K_L}$ which coincides with $C'$ on the arcs of $K_L$ and is integer and bigger than $L$ on the arcs coming from the source or going to the target.

On the other hand, by Ford-Fulkerson \citep{cormen2022introduction}, we know that there exists an ST-circulation $q$ on $G_{K_L}$ with integer non-negative coefficients, which solves $\text{ST-}\f(C'_{\mid K_L})$.
Let $\bar q:E_{K'}\to\mathbb N$ be the extension of $q|_{E_{K_L}}$ by zero outside
$K_L$. For every $i,r$, we have \(\sum_{e\in E_i^r}\bar q(e)\leq L\).
Fix a layer $i$. By \Cref{lem:finiteness} applied to $\bar q$ with $a=0$, there exist two different decks $r,r'\in [0,H_L]$ such that
\(q(e)=q(T_{(r'-r)}e)\) for every \(e\in E_i^r\).
    
    Suppose $r<r'$.
    Since the two profiles on \(E_i^r\) and \(E_i^{r'}\) coincide arc by arc, the
restriction of \(q\) to the slab between these two cuts has matching boundary
data. Hence we may glue copies of this slab periodically along the cuts \(E_i^r\) and \(E_i^{r'}\). In other words, we can periodically extend $q$ above $r'$ and below $r$, obtaining $\Tilde{q}\in \Fop_1^+(K',\Z)$, with the same flow as $q$ and $\Tilde{q}\leq C'$. In this way we have proved that $\text{ST-}\f(C'_{\mid K_L})\leq \f(f',C')$, and the equality follows.
    In particular $\f(f',p'^*C)$ always admits a solution with integer values. 

    We can now consider the operator $\varphi_H$ defined in \Cref{teo:unrolled_flow}.
Let $\widetilde q\in \Fop_1^+(K',\Z)$ be the integer circulation obtained above,
with \(\widetilde q\leq C'=p'^*C\) and
\[
\textstyle\Fop_1(f')(\widetilde q)=\f(f',C').
\]
Applying \Cref{teo:unrolled_flow} to \(\widetilde q\), we obtain
\[
\textstyle\varphi_H(\widetilde q)\in \Fop_1^+(K,\Q),\qquad \varphi_H(\widetilde q)\leq C,
\]
and \(\varphi_H\) preserves flow. Therefore
\[
\textstyle\f(f,C)\geq \Fop_1(f)(\varphi_H(\widetilde q))
   = \Fop_1(f')(\widetilde q)
   = \f(f',C').
\]
Together with the inequality \(\f(f',p'^*C)\geq \f(f,C)\), obtained by lifting
admissible circulations along \(\Fop_1(p'^*)\), this proves
\(
\f(f,C)=\f(f',p'^*C).
\)

Lastly, using \Cref{prop:equivalence} and \Cref{rmk:circul}, since we have integer capacities, we know we can find a solution of $\f(f,C)$ with integer coordinates. 
    \end{proof}

\begin{propfromcor}\label{cor:real_capacities}
For every capacity function \(C:E_K\longrightarrow \mathbb{R}_{\ge 0}\), one has
\(\f(f,C)=\f(f',p'^*C)\).
\end{propfromcor}

The proof is given in Appendix~\ref{app:circular_flow_proofs}. It first uses scaling
for rational capacities and then continuity together with rational approximations
from above and below for arbitrary real capacities.

\section{Max-Flow for Directed Reeb Graphs}\label{sec:reeb_pipeline}

Finally, we assemble together our  contributions. \Cref{sec:reeb} provides the directed Reeb representation and shows that directed tunnel paths are realized by directed paths in the original space, while \Cref{sec:circ_flow} supplies the scalar circular max-flow descriptor. \Cref{sec:mapper_max_flow} connects these two constructions by providing geometric quantities that give natural and interpretable candidates for the capacity constraints of the max-flow problem. It remains only to choose a simplicial model, identify its pullback with the corresponding topological construction, and apply the descriptor.

\subsection{A simplicial model for \(\Rcal(f)\).}

\begin{figure}[!t]
    \centering
    \includegraphics[width = \textwidth]{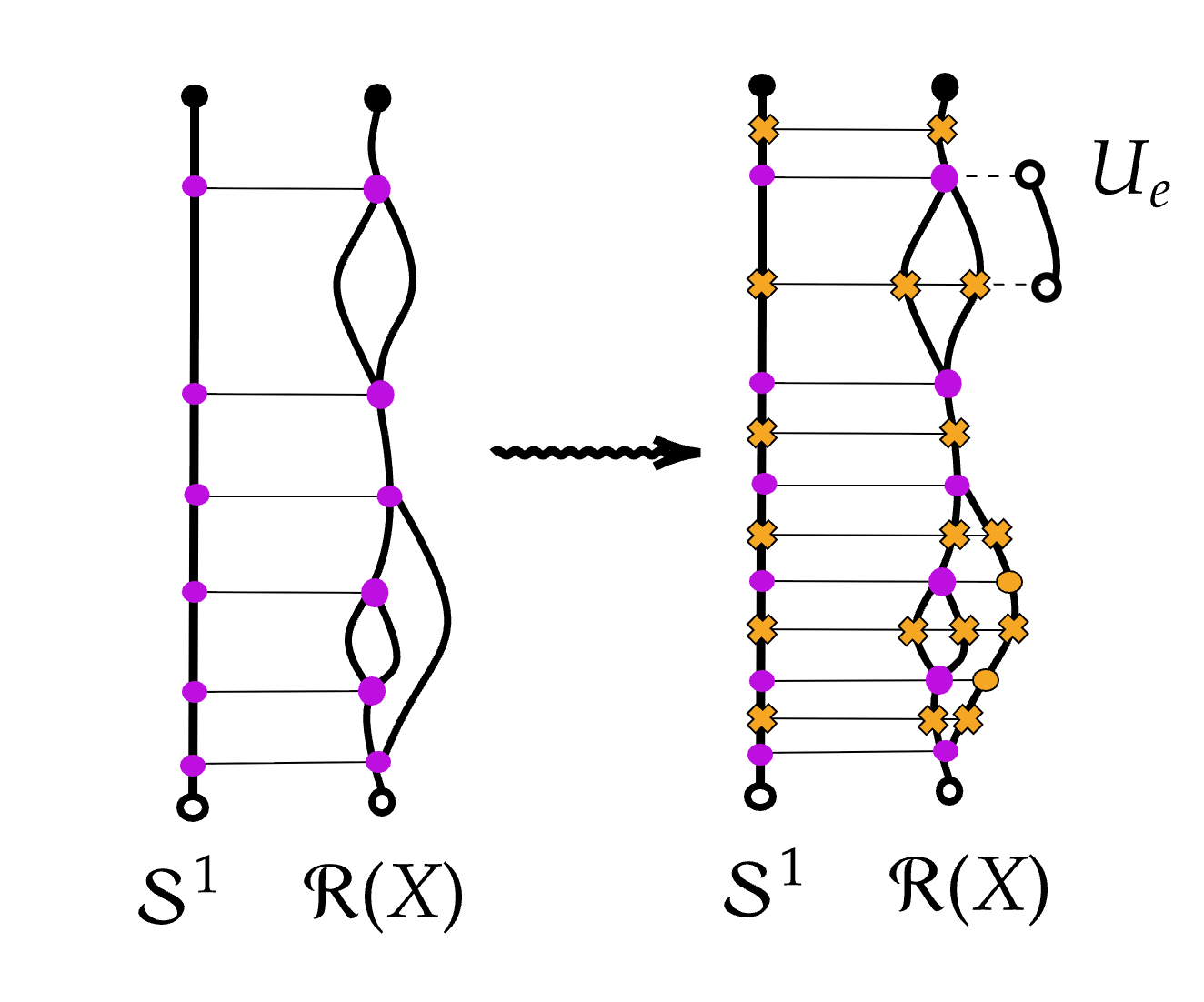}
\caption{Construction of a simplicial model of the Reeb map.
On the left, purple dots mark the topological critical points of
\(\Rcal(f)\) on \(\Rcal(X)\) and their critical values on \(\Ss\).
On the right, orange dots mark the additional preimages of these
critical values. Orange crosses mark the midpoints inserted between consecutive
critical values on \(\Ss\) and all their preimages on \(\Rcal(X)\),
giving compatible first barycentric subdivisions of the circle and
the resulting multigraph stratification.
The dots and crosses on \(\Rcal(X)\) together form the simplicial
singular set \(S_{\Rcal(X)}\).
Horizontal lines connect the displayed Reeb-graph vertices to their
images on \(\Ss\). The arc \(U_e\) illustrates an open edge of the
resulting simplicial model.}
\label{fig:reeb_stratification}
\end{figure}

Our first step is to replace the continuous map \(\Rcal(f):\Rcal(X)\to \Ss\) by a directed
simplicial map \(F:K\to Q\) between graphs (see \Cref{fig:reeb_stratification}).
We adopt the notation of \Cref{sec:circ_flow}, but we write simplicial maps with capital letters.

Choose a finite set \(A\subset\Ss\) containing \(\Crit(\Rcal(f))\) and at least three
points, adding regular values if necessary. The fibres over \(A\) give a finite
multigraph stratification of \(\Rcal(X)\): each component of
\(\Rcal(X)\setminus\Rcal(f)^{-1}(A)\) is an open edge mapped homeomorphically
onto a component of \(\Ss\setminus A\), but distinct edges may have the same
endpoints. Thus this stratification need not yet be a simplicial complex.

To obtain a simplicial model, we take compatible first barycentric subdivisions.
Let \(A^{\mathrm{sd}}\) be obtained from \(A\) by inserting the midpoint of each
complementary arc of \(\Ss\), and define the simplicial singular set
\[
S_{\Rcal(X)}:=\Rcal(f)^{-1}(A^{\mathrm{sd}}).
\]
Here the midpoint of each Reeb edge is chosen using its height map to the
corresponding arc of \(\Ss\). Each pair of parallel edges therefore acquires
distinct intermediate vertices.

We define \(Q\) and \(K\) by:
\begin{itemize}
    \item \(Q\) has vertex set \(A^{\mathrm{sd}}\), and its open edges are the
    path-connected components of \(\Ss\setminus A^{\mathrm{sd}}\).
    \item \(K\) has vertex set \(S_{\Rcal(X)}\), and its open edges are the
    path-connected components of \(\Rcal(X)\setminus S_{\Rcal(X)}\).
\end{itemize}
These are simplicial graphs. By construction,
\(\Rcal(f)(S_{\Rcal(X)})\subseteq A^{\mathrm{sd}}\), and each open edge of \(K\)
is mapped homeomorphically onto an open edge of \(Q\). Using these maps as
edge coordinates, \(\Rcal(f)\) induces a rigid simplicial map
\(F:K\longrightarrow Q\). The chosen directed structure on \(\Ss\) orients the
edges of \(Q\), and pulling it back along \(\Rcal(f)\) orients the edges of \(K\),
making \(F\) a directed simplicial map.

Coarser simplicial models may be chosen when the particular Reeb graph permits
them. After assigning arc capacities below, we verify that subdivisions preserve
the max-flow and that a common refinement makes its value independent of the
chosen simplicial model.

\subsection{Pullback.}
As in \Cref{sec:int_flow}, using the same notation (but capital letters), we unroll the simplicial map \(F\) with a pullback. This yields
a commutative diagram of directed graphs
\begin{equation}\label{eq:simplicial_pullback}
    \begin{tikzcd}
K' \ar[r,"F'"]\ar[d,"P'"]&Q'\ar[d,"P"]\\
K\ar[r,"F"]&Q.
\end{tikzcd}
\end{equation}
The horizontal maps are quasi-finite simplicial maps, and the geometric realization of
Diagram \ref{eq:simplicial_pullback} is homeomorphic to the corresponding pullback diagram in \ref{eq:pullbacks}.

\FloatBarrier

\subsection{Capacities on arcs}\label{sec:arc_capacities}
We now equip the directed simplicial models \(F:K\to Q\) and
\(F':K'\to Q'\) constructed above with arc capacities. Write \(c_X\) for the
pointwise capacity supplied by \Cref{cor:capacity}. For each arc \(e\in E_K\),
let \(U_e\) be the corresponding open edge of
\(\Rcal(X)\setminus S_{\Rcal(X)}\). Since the vertices include every fiber
over a critical value of \(\Rcal(f)\), the function \(c_X\) is continuous on
\(U_e\). Define
\[
C_X:E_K\longrightarrow\R_{\ge0},
\qquad
C_X(e):=\inf_{y\in U_e}c_X(y).
\]
Thus \(C_X(e)\) is the bottleneck capacity of the tunnel represented by
\(e\): a nonnegative flow value assigned to this arc satisfies the capacity
bound at every point of \(U_e\) precisely when it is at most \(C_X(e)\).
We use an infimum because \(U_e\) is open and the bottleneck need not be
attained. Note that $c_X$ is continuous on $U_e$.

For an arc \(e'\in E_{K'}\), let \(U_{e'}\) denote its corresponding open
edge in the unrolled Reeb graph, and set
\[
C_{X\times_{\Ss}\R}:E_{K'}\longrightarrow\R_{\ge0},
\qquad
C_{X\times_{\Ss}\R}(e')
:=\inf_{y'\in U_{e'}}c_{X\times_{\Ss}\R}(y').
\]
By \Cref{sec:capacity_pullback},
\(c_{X\times_{\Ss}\R}=c_X\circ\Rcal(p')\). Moreover,
\(\Rcal(p')\) maps \(U_{e'}\) homeomorphically onto \(U_{P'(e')}\),
so
\[
\begin{aligned}
C_{X\times_{\Ss}\R}(e')
&=\inf_{y'\in U_{e'}}c_X\bigl(\Rcal(p')(y')\bigr)\\
&=\inf_{y\in U_{P'(e')}}c_X(y)
=C_X\bigl(P'(e')\bigr).
\end{aligned}
\]
Consequently, the arc capacities satisfy the pullback relation required
for the flow comparison:
\[
C_{X\times_{\Ss}\R}=C_X\circ P'.
\]

\subsection{Equality of flows.}
We can now apply the circular max-flow of \Cref{sec:circ_flow} to the directed Reeb representation. For the chosen simplicial model \(F:K\to Q\) and the arc capacity \(C_X\) of \Cref{sec:arc_capacities}, define: 
\begin{equation}\label{eq:reeb_flow}
\f(\Rcal(f),c_X):=\f(F,C_X).    
\end{equation}

By \Cref{cor:real_capacities}, the circular max-flow is invariant under the pullback in
Diagram \ref{eq:simplicial_pullback}. Therefore,
\begin{equation}\label{eq:mapper_flow_eq}
\f(\Rcal(f),c_X)\;:=\;\f(F,C_X)
\;=\;\f(F',C_{X\times_{\Ss}\R})
\;=:\;\f(\Rcal(f'),c_{X\times_{\Ss}\R}).
\end{equation}

\begin{rmk}
The value \(\f(\Rcal(f),c_X)\) is independent of the chosen compatible simplicial
model. To see this, suppose an open edge \(U\) is subdivided at a point \(w\)
over a regular value, giving open edges \(U_1\) and \(U_2\). Continuity of \(c_X\)
at \(w\) gives
\[
\inf_{y\in U}c_X(y)
=\min\left\{\inf_{y\in U_1}c_X(y),\inf_{y\in U_2}c_X(y)\right\}.
\]
Indeed, \(c_X(w)\) is a limit from both sides and hence cannot be smaller than
either infimum. Conservation at the new degree-\(2\) vertex forces the same
circulation coefficient on the two new arcs. Feasible circulations before and
after subdivision therefore correspond bijectively and have the same flow.

For a fixed circle model, a compatible model of the Reeb map must include the
entire fibre over each circle vertex. Two choices of compatible simplicial
models admit a common refinement, obtained by taking the union of their circle
vertex sets and all its preimages. Both vertex sets already contain every
critical value, so the additional subdivisions occur over regular values.
Applying the preceding argument proves that both models give the same max-flow.
\end{rmk}

\subsection{Conclusions.}
In conclusion, our combined pipeline is fully captured by the diagrams in \ref{eq:pullbacks} (interpreted in \(\dTop\)),
\Cref{eq:reeb_flow}, and \Cref{eq:mapper_flow_eq}. Taken together, these results show how a directed
structure in the data can be encoded via Reeb graphs, and how this encoding leads to a computable
invariant \(\f(\Rcal(f),c_X)\) that quantifies the maximum flow along a chosen direction in $X$,
with capacities determined by cross-sectional tunnel areas represented on the Reeb graph. Crucially,
this invariant is well defined under periodic boundary conditions and is compatible with modeling
$X$ either directly as a subset of a torus or through its periodic unrolling along the chosen
direction, in line with standard modeling practice for periodic systems. The simulation studies in
\Cref{sec:sim_1,sec:sim_2} illustrate these points in two complementary ways: they show that the
resulting circular max-flow curves discriminate between different spatial organizations of point
cloud backbones, and they highlight the instability of classical source–target max-flow on periodic
domains compared with the intrinsic formulation proposed here.

\section{Discussion}
\label{sec:discussion}

The paper has two main constructions. First, it represents a periodic tunnel system by a directed
Reeb graph. The Reeb quotient compresses a space \(X\) equipped with
\(f:X\to\Ss\) into a one-dimensional network whose edges track level-set components and whose
branching records tunnel reconnections. Pulling back the orientation of \(\Ss\) makes this network
directed, and the path-lifting result shows that its directed tunnel segments correspond to actual
directed paths in \(X\). This gives a compact representation of directional connectivity before
any scalar summary is chosen.

Second, the paper turns that representation into a real-valued descriptor. Fiberwise
Hausdorff-measure integrals provide pointwise cross-sectional capacities; their continuity at
regular and Morse critical values supports the use of edgewise bottlenecks. Circular max-flow then
summarizes the capacity-augmented graph without introducing an arbitrary source, target, or cut.
The resulting optimization is related to minimum-cost circulation and gives the same value on the
circular graph and on its periodically unrolled counterpart.

For the motivating application, \(X\) is the void around a thickened finite atomic backbone in
\((\Ss)^3\), and a coordinate projection selects the transport direction. The first step
describes the void as a periodic tunnel network; the second measures the maximum throughput that
its cross-sections can support. The simulations illustrate both the sensitivity of this descriptor
to spatial organization and the instability created by imposing a source--target cut on a torus.

Several directions remain open, especially concerning the novel max-flow formulation. First, the circular topology complicates the classical max-flow duality: standard
reachability and cut constructions do not directly apply, and establishing an analogue of the
max-flow/min-cut principle may require more general topological tools (cf.\
\cite{hoffman1974generalization, chambers2012homology, krishnan2014flow}). Second, classical loop parametrizations of circulation cones suggest a natural connection
between circular flow, minimum-cost circulations, and knapsack-type formulations. Indeed,
expressing a non-negative circulation in terms of directed cycles turns capacity constraints
on edges into linear packing constraints on cycle coefficients. This perspective points to
possible links with multidimensional knapsack problems and suggests that ideas from the
optimization literature (e.g., \cite{roy1998maximum, laurent1992characterization}) could be
useful for developing heuristics or relaxation guarantees in this setting. We leave a
systematic investigation of these connections to future work.

\FloatBarrier

\section*{Statements and Declarations}

\paragraph{Acknowledgments}
M.P. thanks Mathieu Carrière for hosting him at Centre Inria d'Université Côte d'Azur in 2024, and Steve Oudot for hosting him at Centre Inria de Saclay in 2024; the discussions during these visits contributed to the development of this project.

\paragraph{Funding.}
This work was initiated and primarily supported by the Independent Research Fund Denmark (grant no.\ 1026-00037). Parts of its subsequent development were carried out while M.P. was supported by the Wallenberg AI, Autonomous Systems and Software Program (WASP) and by the SciLifeLab and Wallenberg National Program for Data-Driven Life Science (DDLS) through the project Topological Data Analysis of Functional Genome to Find Covariation Signatures (WASPDDLS22:009), and later by the Fondo Istituzionale per la Ricerca of Università della Svizzera italiana through the project Using Topological Data Analysis to Understand Microglia Shape Variability in Space and Time.

\paragraph{Competing interests.}
The authors declare that they have no competing interests.

\paragraph{Data availability.}
No external datasets were used in this study. All data reported in the simulation studies are
synthetically generated, and the corresponding data-generating processes and parameter choices are
described in the appendix.

\paragraph{Author contributions.}
M. P. conceived the project, developed the methodology, established the main mathematical results, and prepared the original draft. L. F. contributed to the development and refinement of the mathematical framework through discussion and critical feedback, supervised the project, and contributed to the writing, review, and editing of the manuscript. Both authors read and approved the final manuscript.

\paragraph{Use of large language models}

Large language models were used occasionally to improve the clarity and style of the manuscript and to assist in checking the consistency and correctness of selected arguments. All mathematical content and final judgments remain the responsibility of the authors.

\newpage

\FloatBarrier
\appendix

\section{Simulations}

\subsection{Simulation Study 1: Capturing Spatial Structure Around a Backbone}
\label{sec:sim_1}

This first study illustrates, in the periodic setting of the introduction, how the pipeline of
\Cref{sec:reeb_pipeline} summarizes the geometry of a backbone point cloud \(B\) through circular
max-flow. Concretely, we embed \(B\) in the torus \(T=\Ss\times\Ss\times\Ss\) (periodic boundary
conditions), thicken \(B\) by balls of radius \(r\), and study the induced ``tunnel structure'' of
the corresponding void space. This is the geometric situation discussed in the materials-science
motivation and is covered by \Cref{ex:atoms}.

\begin{figure}[htbp]
\centering
\begin{subfigure}{0.47\textwidth}
    \includegraphics[width = 0.9\textwidth]{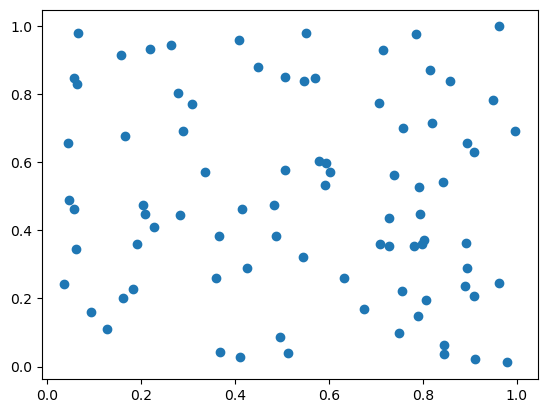}
\caption{The 2D projection of a realization of a homogeneous Poisson process.}
\label{fig:CSR}
\end{subfigure}
\centering
\begin{subfigure}{0.47\textwidth}
    \centering
    \includegraphics[width = \textwidth]{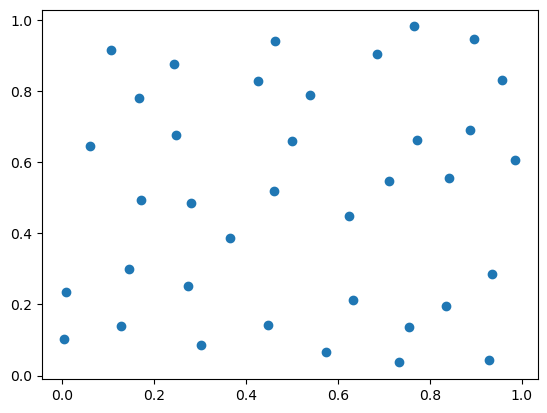}
\caption{The 2D projection of a realization of a Mat\'ern repulsive process.}
\label{fig:rep}
\end{subfigure}

\begin{subfigure}{0.47\textwidth}
    \includegraphics[width = 0.9\textwidth]{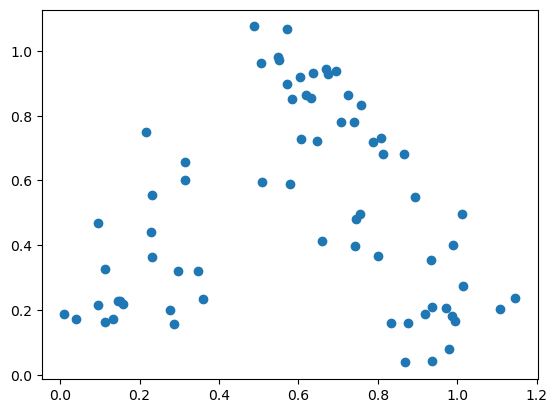}
\caption{The 2D projection of a realization of a Mat\'ern cluster process.}
\label{fig:clust}
\end{subfigure}
\centering
\begin{subfigure}{0.47\textwidth}
    \centering
    \includegraphics[width = \textwidth]{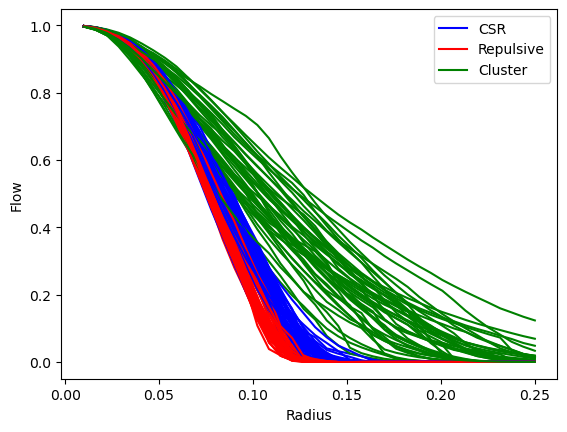}
\caption{Circular max-flow as a function of the thickening radius \(r\), colored by the generating process.}
\label{fig:ris_1}
\end{subfigure}

\caption{Results and plots for Simulation Study 1.}
\label{fig:sim_1}
\end{figure}

\paragraph{Point process models.}
We generate i.i.d.\ realizations of \(B\subset[0,1]^3\) from three classical families of stationary
point processes (PPs) using standard simulation methods \citep{baddeley2015spatial}. The observation
window is \([0,1]^3\) for all models.
\begin{itemize}
    \item \textbf{Homogeneous Poisson (CSR).}
    Points are sampled uniformly and independently in \([0,1]^3\); the parameter is the intensity
    \(\lambda\), i.e.\ the expected number of points in the window. We take \(\lambda=200\)
    (\Cref{fig:CSR}).

    \item \textbf{Mat\'ern repulsive (hard-core thinning).}
    Starting from a homogeneous Poisson process with intensity \(\lambda\), we apply a Mat\'ern-type
    thinning: points that fall within distance \(R>0\) of another point are removed according to a
    random ordering. We use \(\lambda=400\) and \(R=0.1\) (\Cref{fig:rep}).

    \item \textbf{Mat\'ern cluster (Neyman--Scott type).}
    We first sample ``parent'' points by a homogeneous Poisson process of intensity \(\lambda_1\).
    For each parent, we draw a Poisson number of offspring with mean \(\lambda_2\), using
    independent coordinate-wise offsets sampled uniformly from \([0,R]^3\). We take
    \(\lambda_1=20\), \(\lambda_2=10\), and \(R=0.2\) (\Cref{fig:clust}).
\end{itemize}

\paragraph{Thickening and void space.}
Given a realization \(B\), we pass to \(T=(\Ss)^3\) via \(\Ss\cong[0,1]/(0\sim 1)\), fix \(r>0\), and
define the thickened backbone and its complement
\[
\Omega_r := T\setminus \widetilde{B}^{\,r},
\qquad
X_r := \overline{\Omega_r}.
\]
We then run the pipeline of \Cref{sec:reeb_pipeline} on \(X_r\), using the projection onto the last
\(\Ss\)-coordinate as the map \(X_r\to \Ss\) (cf.\ \Cref{ex:atoms}).

\paragraph{Discretization and protocol.}
We discretize \(T\) using a uniform grid with \(101^3\) points and
form a cubical complex. We approximate the Reeb graph using
Mapper, with intervals centered at every other grid value of
the coordinate projection. Each interval covers three consecutive
grid values and overlaps with its neighboring intervals at the
extremes; the same convention is applied periodically at the
endpoints of the coordinate range.

To estimate the cross-sectional capacities, we count the grid
points in each path-connected component of each discretized
level set and multiply this count by \(h^2\), where \(h\) is the
grid spacing. This gives a discrete approximation of the
two-dimensional Hausdorff measure of the corresponding cross
section. Following \Cref{sec:arc_capacities}, we assign each graph edge the
minimum of these estimates along the corresponding tunnel
segment.
For each PP family we generate \(50\) independent backbones. For each realization we
compute the circular max-flow descriptor at \(40\) equally spaced radii \(r\in[0.01,0.25]\).
The PP parameters are chosen so that the expected number of sampled points is comparable across the three families.

\paragraph{Results.}
The resulting curves are shown in \Cref{fig:sim_1}. Circular max-flow separates the three regimes:
for very small radii, the thickenings \(\widetilde{B}^{\,r}\) remain essentially disjoint and the
void spaces \(X_r\) are similar across models, yielding comparable flow values. As \(r\) increases
and the balls begin to overlap, the differences in spatial organization become apparent. Repulsive
processes distribute points more evenly, which closes off corridors earlier and produces smaller
flows; clustered processes concentrate points locally, leaving larger open regions elsewhere, and
thus yield larger flows; Poisson realizations typically lie between these two behaviors.

\subsection{Simulation Study 2: Instability of Source--Target Max-Flow}
\label{sec:sim_2}

This second study illustrates the obstruction highlighted in the introduction: on a torus there is
no canonical inlet/outlet, so any source--target max-flow obtained by ``cutting'' the domain depends
on an arbitrary choice of coordinates. We show that this dependence can be substantial even for
simple geometric signals. Periodicity is kept in the transverse directions.

\begin{figure}[htbp]
\centering
\begin{subfigure}{0.47\textwidth}
    \includegraphics[width = 0.9\textwidth]{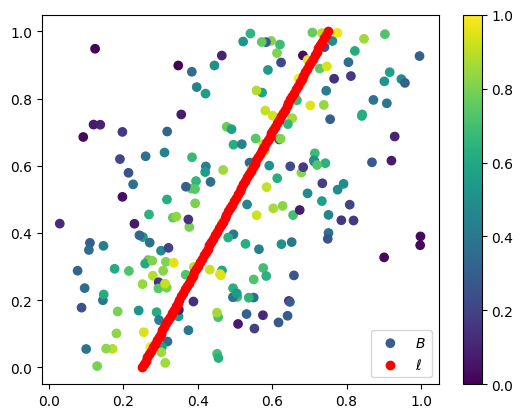}
\caption{A $2D$ projection of a realization of the inhomogeneous PP. Colors indicate thinning acceptance probability; the red segment is \(\ell\).}
\label{fig:inhomo}
\end{subfigure}
\centering
\begin{subfigure}{0.47\textwidth}
    \centering
    \includegraphics[width = 0.9\textwidth]{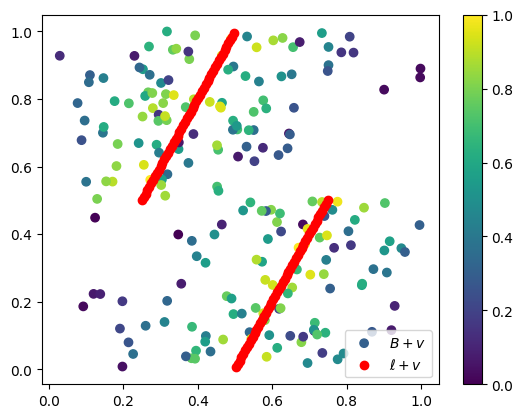}
\caption{The same realization and segment \(\ell\), translated by \(v=(0,0,0.5)\) on the torus and shown in the fixed coordinate box.}
\label{fig:X_v}
\end{subfigure}

\begin{subfigure}{0.47\textwidth}
    \includegraphics[width = 0.9\textwidth]{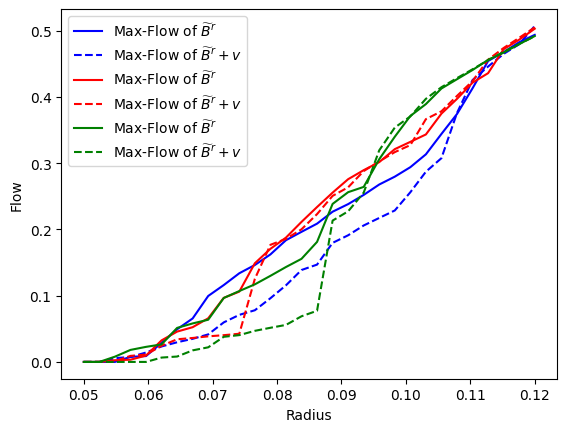}
\caption{Source--target max-flow vs.\ radius \(r\): original \(\widetilde{B}^{\,r}\) (solid) and translated \(\widetilde{B}^{\,r}+v\) (dashed).}
\label{fig:simulation_2_ex}
\end{subfigure}
\centering
\begin{subfigure}{0.47\textwidth}
    \centering
    \includegraphics[width = 0.9\textwidth]{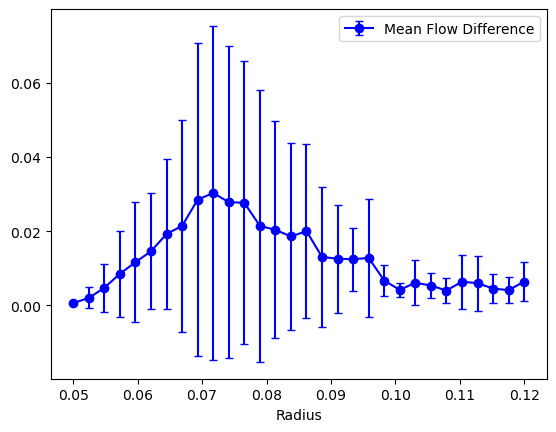}
\caption{Mean absolute difference between source--target flows of \(\widetilde{B}^{\,r}\) and \(\widetilde{B}^{\,r}+v\); the symmetric error-bar magnitude equals the pointwise interquartile range.}
\label{fig:simulation_2_ris}
\end{subfigure}

\caption{Results and plots for Simulation Study 2.}
\label{fig:sim_2}
\end{figure}

\paragraph{A structured inhomogeneous point cloud.}
We generate backbones \(B\subset[0,1]^3\) from an inhomogeneous Poisson process obtained by thinning
a homogeneous Poisson sample of intensity \(\lambda=800\). We represent a diagonal segment by
\(100\) equally spaced points \(\ell_1,\ldots,\ell_{100}\), and retain each candidate point \(x\)
independently with probability
\[
\mathbb{P}(\text{keep }x)
=\max_{1\le j\le100}\exp\!\Bigl(-\frac{\lVert x-\ell_j\rVert^2}{2\sigma^2}\Bigr),
\qquad
\sigma=0.2,
\]
where the discretized diagonal connects \(p_0=(0.25,0.5,0)\) to
\(p_1=(0.75,0.5,1)\) in \([0,1]^3\).
An example realization (projected to \(2D\)) is shown in \Cref{fig:inhomo}.

\paragraph{Thickening and two flow constructions.}
For each realization and each of \(30\) equally spaced radii \(r\in[0.05,0.12]\), we form the thickening
\(\widetilde{B}^{\,r}\subset (\Ss)^3\) exactly as in Simulation~1 (this setup is still covered by
\Cref{ex:atoms}). We then compare:
\begin{itemize}
    \item \textbf{Circular max-flow,} computed intrinsically on the directed Reeb graph produced by
    the pipeline of \Cref{sec:reeb_pipeline};
    \item \textbf{Source--target max-flow,} computed after choosing coordinates on \((\Ss)^3\), ``cutting''
    at height \(1\) (opening the periodic direction), and placing a source at height \(0\) and a
    target at height \(1\) (as in \Cref{fig:augmented}).
\end{itemize}

\paragraph{A translation experiment.}
Let \(v=(0,0,0.5)\) and consider the translated backbone \(B+v\) (projected to \((\Ss)^3\)).
Geometrically, this translation can split the lifted segment \(\ell\) (when represented in the
fixed coordinate box \([0,1]^3\)) into two components, even though the configuration on the torus is
the same up to translation; cf.\ \Cref{fig:X_v}.
Circular max-flow is invariant under such translations, while the source--target construction depends
on where the cut is placed relative to the pattern.

For each realization we compute source--target max-flow curves for \(\widetilde{B}^{\,r}\) and for
\(\widetilde{B}^{\,r}+v\). Representative examples are shown in \Cref{fig:simulation_2_ex}, where the
continuous curve (original) and dashed curve (translated) differ markedly.

\paragraph{Quantifying instability.}
We repeat the experiment \(20\) times independently. For each radius \(r\), we compute the mean
absolute difference between the source--target flows of \(\widetilde{B}^{\,r}\) and
\(\widetilde{B}^{\,r}+v\). The symmetric error-bar magnitude is the pointwise interquartile range
across realizations. The results in
\Cref{fig:simulation_2_ris} show large variability relative to the flow magnitude (notably around
\(r\in[0.07,0.08]\)), confirming that source--target max-flow can be highly sensitive to the
coordinate-dependent choice of how the torus is opened. This is precisely the ambiguity emphasized
in the materials-science discussion: in a periodic setting there is no preferred coordinate system,
hence no preferred source--target placement.

\FloatBarrier
\setcounter{table}{0}
\setcounter{equation}{0}
\section{Proofs from Section 2}\label{app:reeb_pipeline_proofs}
This appendix contains the proof deferred from \Cref{sec:reeb} and an
auxiliary edge-lifting lemma used in it.

\begin{lem}\label{lem:edge_tunnel_lift}
Let \(g:Z\to \R\) be a countably constructible \(\R\)-space, and let
\(\Rcal(g):\Rcal(Z)\to\R\) be its Reeb map.
Let \(I\subset \Rcal(Z)\) be a compact arc contained in the closure of a single open
edge of \(\Rcal(Z)\), and let \(\beta:[0,1]\to I\) be a homeomorphism.

Then there exists a path \(\beta':[0,1]\to Z\) such that
\[
\pi_g(\beta'([0,1]))=I,
\qquad
\pi_g(\beta'(0))=\beta(0),\qquad
\pi_g(\beta'(1))=\beta(1).
\]

Moreover, if \(\Rcal(g)\circ \beta\) is nondecreasing, then \(\beta'\) may be chosen so
that \(g\circ \beta'\) is nondecreasing.
\end{lem}

\begin{proof}
The compact set $\Rcal(g)(I)$ is contained in the interior of a compact interval
$[A,B]\subset\R$. By countable constructibility, the restricted space
\(g^{-1}([A,B])\longrightarrow[A,B]\) admits a finite constructible presentation
with presentation values \(A=t_0<t_1<\cdots<t_N=B\).
Thus $g^{-1}([A,B])$ is homeomorphic, as an $[A,B]$-space, to
\[
\Bigl(\coprod_{i=0}^N(V_i\times\{t_i\})\Bigr)\ \coprod\
\Bigl(\coprod_{i=0}^{N-1}(E_i\times[t_i,t_{i+1}])\Bigr)\Big/\sim,
\]
where the equivalence relation is generated by the attaching maps
\(l_i:E_i\to V_i\) and \(r_i:E_i\to V_{i+1}\).

Write \(a:=\Rcal(g)(\beta(0))\) and \(b:=\Rcal(g)(\beta(1))\).
For the general lifting statement, the case \(a>b\) is obtained from the case
\(a\le b\) by reversing the parameter. For the monotone statement, the additional
assumption that \(\Rcal(g)\circ\beta\) is nondecreasing already implies \(a\le b\).
Thus it is enough to treat the case \(a\le b\).

On the closure of an open edge of an \(\R\)-Reeb graph, the height map is
monotone and restricts to a homeomorphism onto an interval. Hence
\(\Rcal(g)|_I:I\to[a,b]\) is a homeomorphism. Subdivide \([a,b]\) at the
finitely many presentation values \(t_i\) that it meets,
\(a=u_0<u_1<\cdots<u_m=b\).
For every \(k=1,\ldots,m\), the open interval
\((u_{k-1},u_k)\) lies in one \((t_i,t_{i+1})\). The corresponding part of
the Reeb edge is represented by a path-connected component
\(C_k\in\pi_0(E_i)\). Choose \(c_k\in C_k\). The vertical segment
\(\{c_k\}\times[u_{k-1},u_k]\)
with the appropriate attaching-map value at either endpoint that equals
\(t_i\) or \(t_{i+1}\), defines a continuous path \(\gamma_k\) in \(Z\).
Its Reeb image is precisely the subarc
\(I_k:=I\cap\Rcal(g)^{-1}([u_{k-1},u_k])\),
and \(g\circ\gamma_k\) is nondecreasing from \(u_{k-1}\) to \(u_k\).

At an interior subdivision value \(u_k\), the endpoint of \(\gamma_k\) and
the initial point of \(\gamma_{k+1}\) represent the same point of the Reeb
graph. They therefore lie in the same path-connected component of
\(g^{-1}(u_k)\), and may be joined by a path \(\eta_k\) contained in that
component. Along \(\eta_k\), both \(g\) and \(\pi_g\) are constant.
Concatenating \(\gamma_1,\eta_1,\gamma_2,\ldots,\eta_{m-1},\gamma_m\)
and reparametrizing gives a path \(\beta':[0,1]\to Z\) whose Reeb image is
\(I\), whose endpoints project to \(\beta(0)\) and \(\beta(1)\), and for
which \(g\circ\beta'\) is nondecreasing. Reversing this path treats the case
\(a>b\), while the nondecreasing case is exactly the construction above.
\end{proof}

\subsection*{Proof of \Cref{prop:path_conn_regular}}
\begin{proof}
Consider the pullback of \(f\) along the universal covering
\(p:\R\to\Ss\),
\[
\begin{tikzcd}
X\times_{\Ss}\R \ar[r,"f'"]\ar[d,"p'"] & \R\ar[d,"p"]\\
X\ar[r,"f"] & \Ss .
\end{tikzcd}
\]
Since \(f\) is constructible, 
\(
f':X\times_{\Ss}\R\longrightarrow\R
\)
is countably constructible by definition. Moreover, by
\Cref{lem:pull_back} there is a canonical homeomorphism of
\(\R\)-spaces
\[
\Theta:
\Rcal(X\times_{\Ss}\R)
\xrightarrow{\sim}
\Rcal(X)\times_{\Ss}\R.
\]

Set
\(
h:=\Rcal(f)\circ\alpha:[0,1]\to\Ss.
\)
By the path-lifting property of the covering \(p:\R\to\Ss\), choose
a lift
\[
\widetilde h:[0,1]\to\R,
\qquad
p\circ\widetilde h=h.
\]
The map
\[
t\longmapsto (\alpha(t),\widetilde h(t))
\]
takes values in \(\Rcal(X)\times_{\Ss}\R\). We therefore define
\[
\widetilde\alpha:[0,1]\longrightarrow
\Rcal(X\times_{\Ss}\R),
\qquad
\widetilde\alpha(t)
:=
\Theta^{-1}\bigl(\alpha(t),\widetilde h(t)\bigr).
\]
Since \(\alpha\) is an embedding, so is \(\widetilde\alpha\).
Furthermore,
\[
\Rcal(p')\circ\widetilde\alpha=\alpha,
\qquad
\Rcal(f')\circ\widetilde\alpha=\widetilde h.
\]

The image \(\widetilde\alpha([0,1])\) is compact. Since \(f'\) is
countably constructible, over every compact interval in \(\R\) its
Reeb graph has a finite constructible graph stratification. Hence
there exists a subdivision
\[
0=t_0<t_1<\cdots<t_n=1
\]
such that, for every \(i=1,\ldots,n\), the compact arc
\[
I_i
:=
\widetilde\alpha([t_{i-1},t_i])
\]
is contained in the closure of a single open edge of
\(\Rcal(X\times_{\Ss}\R)\).

For each \(i\), apply \Cref{lem:edge_tunnel_lift} to the
countably constructible \(\R\)-space \(f'\) and to the restriction
of \(\widetilde\alpha\) to \([t_{i-1},t_i]\), after the obvious affine
reparametrization of this interval onto \([0,1]\). Reparametrizing
back, we obtain a path
\[
\gamma_i:[t_{i-1},t_i]\to X\times_{\Ss}\R
\]
such that
\[
\pi_{f'}\bigl(\gamma_i([t_{i-1},t_i])\bigr)=I_i,
\]
and
\[
\pi_{f'}(\gamma_i(t_{i-1}))
=
\widetilde\alpha(t_{i-1}),
\qquad
\pi_{f'}(\gamma_i(t_i))
=
\widetilde\alpha(t_i).
\]

For every interior breakpoint \(t_i\), the points
\[
\gamma_i(t_i),
\qquad
\gamma_{i+1}(t_i)
\]
belong to the same fiber
\(
\pi_{f'}^{-1}\bigl(\widetilde\alpha(t_i)\bigr).
\)
By definition of the Reeb relation this fiber is path-connected, so
there exists a path
\[
\eta_i:[0,1]\to
\pi_{f'}^{-1}\bigl(\widetilde\alpha(t_i)\bigr)
\]
joining \(\gamma_i(t_i)\) to \(\gamma_{i+1}(t_i)\).

Concatenating
\(
\gamma_1,\eta_1,\gamma_2,\eta_2,\ldots,
\eta_{n-1},\gamma_n
\)
and reparametrizing gives a continuous path
\[
\widetilde\alpha':[0,1]\to X\times_{\Ss}\R.
\]
Since each \(\gamma_i\) has Reeb image \(I_i\), while each
\(\eta_i\) is contained in a single fiber of \(\pi_{f'}\), we have
\[
\pi_{f'}\bigl(\widetilde\alpha'([0,1])\bigr)
=
\widetilde\alpha([0,1]),
\]
and
\[
\pi_{f'}(\widetilde\alpha'(0))
=
\widetilde\alpha(0),
\qquad
\pi_{f'}(\widetilde\alpha'(1))
=
\widetilde\alpha(1).
\]

Now define
\[
\alpha':=p'\circ\widetilde\alpha':[0,1]\to X.
\]
Using the commutativity of the Reeb quotient diagram,
\[
\pi_f\circ p'
=
\Rcal(p')\circ\pi_{f'},
\]
we obtain
\[
\begin{aligned}
\pi_f\bigl(\alpha'([0,1])\bigr)
&=
\Rcal(p')\!
\left(
\pi_{f'}\bigl(\widetilde\alpha'([0,1])\bigr)
\right)\\
&=
\Rcal(p')\bigl(\widetilde\alpha([0,1])\bigr)
=
\alpha([0,1]).
\end{aligned}
\]
The same identity at the endpoints gives
\[
\pi_f(\alpha'(0))=\alpha(0),
\qquad
\pi_f(\alpha'(1))=\alpha(1).
\]

It remains to prove the directed statement. We assume without loss of
generality that the universal covering
\(p:\R\to\Ss\) is parametrized so that
\(\Dcal^p(d\Ss)\) is the directed structure on \(\R\)
given by nondecreasing paths. Suppose that
\(\alpha\in d\Rcal(X)\). Since the directed structure on
\(\Rcal(X)\) is pulled back from \(d\Ss\),
\[
h=\Rcal(f)\circ\alpha\in d\Ss.
\]
Since \(p\circ\widetilde h=h\), we have
\[
\widetilde h\in\Dcal^p(d\Ss).
\]
This is precisely the directed structure
on \(\R\) given by nondecreasing paths. Thus
\(\widetilde h\) is nondecreasing.

Since
\[
\Rcal(f')\circ\widetilde\alpha=\widetilde h,
\]
the restriction of \(\Rcal(f')\circ\widetilde\alpha\) to every
\([t_{i-1},t_i]\) is nondecreasing. The directed part of
\Cref{lem:edge_tunnel_lift} therefore allows each \(\gamma_i\) to be
chosen so that
\(
f'\circ\gamma_i
\)
is nondecreasing.

Moreover, along each connecting path \(\eta_i\), the map
\(\pi_{f'}\) is constant. Since
\(
f'=\Rcal(f')\circ\pi_{f'},
\)
the map \(f'\circ\eta_i\) is constant, with value
\(\widetilde h(t_i)\). The successive pieces therefore concatenate
with matching values, and, since \(\widetilde h\) is nondecreasing,
the concatenated path may be parametrized so that
\(
f'\circ\widetilde\alpha'
\)
is nondecreasing on \([0,1]\).

Finally, from the pullback square,
\[
f\circ p'=p\circ f',
\]
and hence
\[
f\circ\alpha'
=
p\circ f'\circ\widetilde\alpha'
\in d\Ss.
\]
By the definition \(dX=\Dcal^f(d\Ss)\), this proves
\(
\alpha'\in dX.
\)
\end{proof}

\FloatBarrier

\FloatBarrier
\setcounter{table}{0}
\setcounter{equation}{0}
\section{Proofs from Section 3}\label{app:capacity_proofs}
This appendix contains the long technical proofs deferred from
\Cref{sec:mapper_max_flow}.

\subsection*{Boundary regularity in \Cref{ex:atoms}}
For either choice of \(\Omega\), its boundary is contained in \(d_B^{-1}(r)\).
The cases \(r<0\) and \(r=0\) give an empty boundary and a boundary
contained in the finite set \(B\), respectively, so assume \(r>0\).
Choose one lift of each point of \(B\) to \(\R^3\), and fix the compact
fundamental domain \(Q=[0,2\pi]^3\). Only finitely many radius-\(r\) balls
centered at lattice translates of these lifts meet \(Q\). The lift of
\(\partial\Omega\) inside \(Q\) is therefore contained in a finite union of
Euclidean spheres. By compactness, each sphere can be covered by finitely many sufficiently
small open Euclidean balls centered on it. Choose these balls so that
the covering projection \(\R^3\to T\) is a diffeomorphism onto its image
on a neighborhood of each closed ball. Intersecting the sphere with
these closed balls gives compact spherical caps with smooth boundary.
Their images in \(T\) are therefore compact smooth embedded surfaces
with smooth boundary, and together they cover \(\partial\Omega\),
as required by \Cref{assumptions}.

Now let \(f\) be a coordinate projection. In the universal cover, every level
of \(f\) is a countable union of parallel coordinate planes. The intersection
of any one of the finitely many spheres above with such a plane is empty, a
point, or a circle, and hence is \(\Hcal^2\)-null. It follows that the
intersection of every two-dimensional boundary piece with \(f^{-1}(t)\) is
\(\Hcal^2\)-null. The critical set of \(f\) restricted to that piece at level
\(t\) is a subset of this intersection, so the a.e. boundary regularity
condition in \Cref{assumptions} holds.

\subsection*{Constructibility in \Cref{ex:atoms}}
Let \(\pi:\R^3\to T\) be the universal covering projection,
\(Q=[0,2\pi]^3\), and choose \(\widetilde B\subset\R^3\) such that the restriction
\(\pi|_{\widetilde B}:\widetilde B\to B\) is bijective. For every \(x\in Q\), the squared Riemannian distance in \(T\)
from \(\pi(x)\) to \(B\) is given by
\[
d_B(\pi(x))^2
=
\min_{\tilde b\in\widetilde B,\,m\in M}
\sum_{j=1}^3(x_j-\tilde b_j-m_j)^2,
\]
where \(M\subset(2\pi\Z)^3\) is finite: only finitely many lattice
translates can realize the minimum on the compact set \(Q\).

Consider first \(\Omega=\{d_B<r\}\), with \(r\geq0\). Then
\[
Q\cap\pi^{-1}(\Omega)
=
\bigcup_{\tilde b\in\widetilde B,\,m\in M}
\left\{
x\in Q:
\sum_{j=1}^3(x_j-\tilde b_j-m_j)^2<r^2
\right\}.
\]
This is a finite union of sets defined by polynomial inequalities,
and hence is semialgebraic. For \(\Omega=\{d_B>r\}\), the same
conclusion follows by replacing the union with an intersection
and \(<r^2\) with \(>r^2\). The cases \(r<0\) give the empty set
or all of \(T\).

Since \(\pi\) is real analytic, the image of this relatively
compact semialgebraic set is subanalytic in \(T\). Thus
\(\Omega\), and also its closure \(X\), are subanalytic
\citep[Section~3, p.~16]{bierstone1988semianalytic}.

We now use the pullback
\(
f':X\times_{\Ss}\R\longrightarrow\R.
\)
Its restriction to
\[
Z:=(f')^{-1}([0,2\pi])
=
\{(x,t)\in X\times[0,2\pi]:f(x)=p(t)\}
\]
is a continuous real-valued subanalytic function on a compact
subanalytic set, because both the coordinate projection \(f\)
and the covering map \(p\) are real analytic.
By \citet[Theorem~1]{shiota1989piecewise}, there exist a finite
simplicial complex \(K\) and a homeomorphism \(h:|K|\to Z\)
such that \(g:=f'\circ h\) is linear on every simplex.

A piecewise-linear real-valued map on a compact polyhedron is
constructible \citep[Example~2.2(ii)]{de2016categorified}.
Consequently, \(f'|_Z\) is constructible in the sense of
\Cref{def:de_constr}. 

Let \(S_0\subset[0,2\pi]\) be a finite
critical set for this presentation, including the endpoints
\(0\) and \(2\pi\).
The deck transformations \((x,t)\mapsto(x,t+2\pi k)\) transport
this presentation to each interval
\([2\pi k,2\pi(k+1)]\). These presentations join along their
common endpoint fibres, with critical values in the closed
discrete set
\[
S=\{s+2\pi k:s\in S_0,\ k\in\Z\}.
\]
Every compact interval meets only finitely many such pieces;
restricting and joining their presentations gives the finite
presentation required by \Cref{def:cont_constr}. Thus \(f'\)
is countably constructible, and \(f:X\to\Ss\) is constructible.

\subsection*{Proof of \Cref{prop:hausd_continuous}}
\begin{proof}

\begin{figure}[htbp]
    \centering
    	\includegraphics[width = 0.8\textwidth]{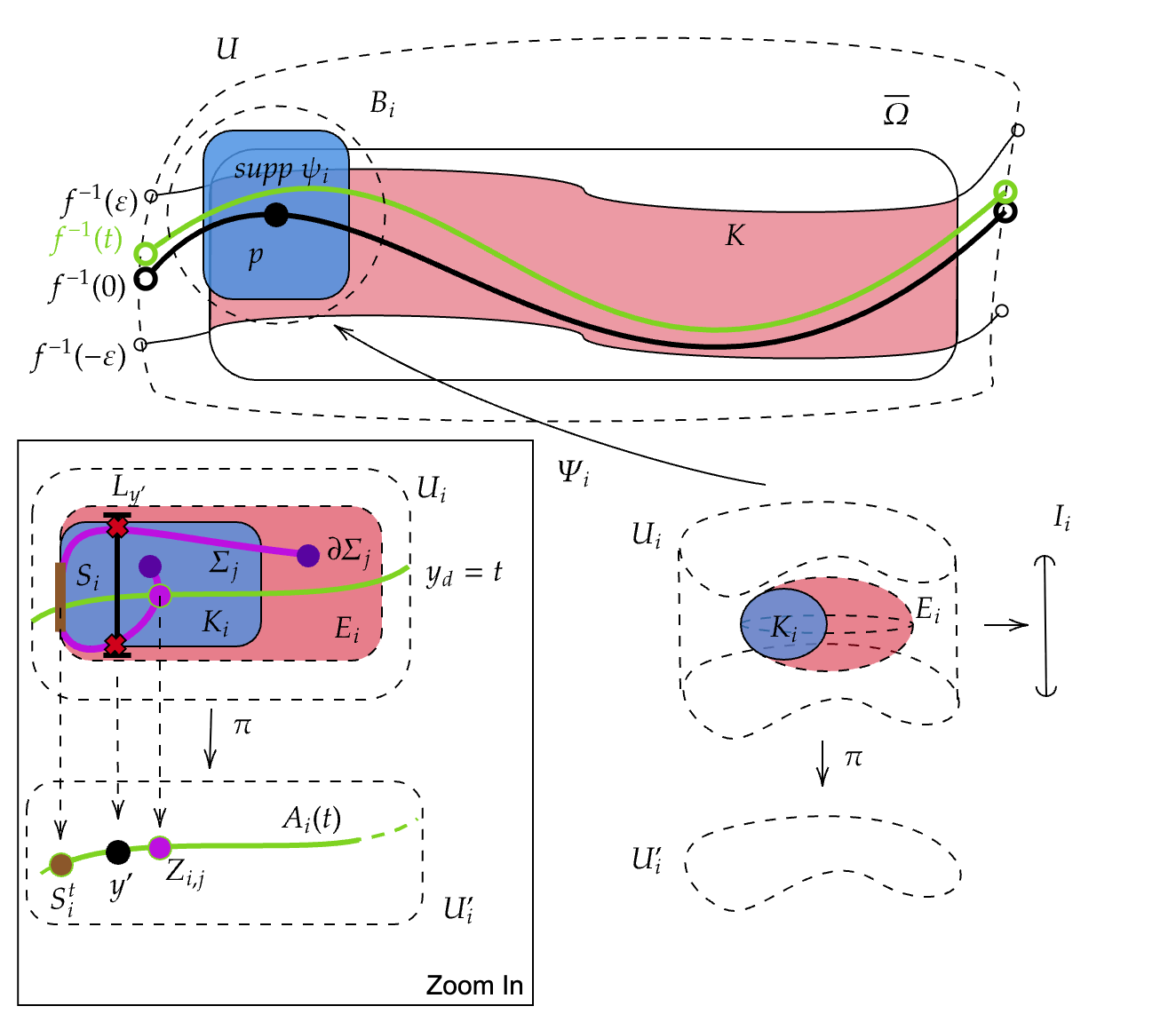}
\caption{A depiction of the local charts involved in the proof of \Cref{prop:hausd_continuous}.}
\label{fig:local_charts}
\end{figure} 

\Cref{fig:local_charts} can help in following the localization arguments pursued in the proof.

Without loss of generality, assume $t_0=0$ by replacing $f$ with $f-t_0$.
Thus $\nabla f(x)\neq 0$ for all $x\in\overline\Omega\cap f^{-1}(0)$.
We prove $c(t)\to c(0)$ as $t\to 0$.

Set
\[
S_0 := \overline\Omega\cap f^{-1}(0),
\]
which is compact. If \(S_0=\emptyset\), compactness of \(\overline\Omega\)
implies that every slice with level sufficiently close to zero is empty,
so the conclusion is immediate. We may therefore assume \(S_0\ne\emptyset\).

\medskip
\noindent\textbf{Step 1: Submersion charts around $S_0$ and a compact strip.}\par
\nopagebreak[4]

For each $p\in S_0$, the submersion theorem yields an open neighborhood
$B_p\subset U$, an open set $U_p\subset\mathbb{R}^d$, and a $C^1$ diffeomorphism
\[
\Psi_p:U_p\to B_p
\]
such that, in coordinates $y=(y',y_d)$ on $U_p$,
\[
f(\Psi_p(y',y_d)) = y_d.
\]
Shrinking the sets if necessary, assume $U_p=U_p'\times I_p$ with $U_p'\subset\mathbb{R}^{d-1}$
open and bounded and $I_p\subset\mathbb{R}$ an open interval containing $0$.

By compactness of $S_0$, choose $p_1,\dots,p_M\in S_0$ such that
\[
S_0 \subset \bigcup_{i=1}^M B_{p_i}.
\]
Set $B_i:=B_{p_i}$, $\Psi_i:=\Psi_{p_i}$, and write $U_i:=U_i'\times I_i$.

Set \(B:=\bigcup_{i=1}^M B_i\subset U\). Since \(\overline\Omega\setminus B\)
is compact and disjoint from \(S_0\), we can choose \(\varepsilon>0\)
such that
\[
[-\varepsilon,\varepsilon]\subset\bigcap_{i=1}^M I_i,
\qquad
K:=\overline\Omega\cap f^{-1}([-\varepsilon,\varepsilon])\subset B.
\]
Indeed, if \(\overline\Omega\setminus B\ne\emptyset\), it suffices also to require
\(\varepsilon<\min_{\overline\Omega\setminus B}|f|\).
Write \(J:=(-\varepsilon,\varepsilon)\). The set \(K\) is compact,
and \(\Omega\cap f^{-1}(t)\subset K\) for every \(t\in J\).

For each $i$, define the pulled-back domain
\[
E_i := \Psi_i^{-1}(\Omega\cap B_i)\subset U_i'\times I_i.
\]

\medskip
\noindent\textbf{Step 2: finite partition of unity around $K$.}\par
\nopagebreak[4]

Since $\{B_i\}_{i=1}^M$ is a finite open cover of the compact set $K$, there exist
functions $\psi_i\in C_c^1(B_i)$ with values in $[0,1]$ such that
\[
\sum_{i=1}^M \psi_i \equiv 1 \quad\text{on }K.
\]

For every $t\in J$, the inclusion $\Omega\cap f^{-1}(t)\subset K$ gives
\[
c(t)
= \int_{\Omega\cap f^{-1}(t)} \mu\,d\mathcal{H}^{d-1}
= \sum_{i=1}^M \int_{\Omega\cap B_i\cap f^{-1}(t)} \mu(x)\,\psi_i(x)\,
  d\mathcal{H}^{d-1}(x)
=: \sum_{i=1}^M c_i(t).
\]
It suffices to prove that each $c_i$ is continuous at $0$.

\medskip
\noindent\textbf{Step 3: Expression of $c_i(t)$ and uniform bounds.}\par
\nopagebreak[4]

Fix $i$. In the coordinates $y=(y',y_d)\in U_i'\times I_i$ we have
$f(\Psi_i(y)) = y_d$, and for $t\in J$ the slice
$\Omega\cap B_i\cap f^{-1}(t)$ corresponds to
\[
\{(y',y_d)\in E_i : y_d=t\},
\qquad \text{ whose projection is }
A_i(t):=\{y'\in U_i':(y',t)\in E_i\}.
\]
Since \(\overline\Omega\cap B_i\) is closed in \(B_i\), we extend
\(\mu|_{\overline\Omega\cap B_i}\) continuously to \(B_i\) using Tietze Extension Theorem \cite[Theorem~35.1]{munkres2000topology}, still denoted by \(\mu\).
By the localization claim below, the choice of extension is irrelevant.

By the area formula for $C^1$ maps, the pullback of $\mathcal{H}^{d-1}$ on $\{y_d=t\}$ via $y'\mapsto \Psi_i(y',t)$ is
$J_i(y',t)\,dy'$, where
\[
J_i(y',t)
:= \sqrt{\det\!\big((D_{y'}\Psi_i(y',t))^\top D_{y'}\Psi_i(y',t)\big)}
\]
is the \((d-1)\)-dimensional Jacobian of \(\Psi_i\)
along the level set hypersurface.
Thus
\[
c_i(t)=\int_{A_i(t)} \Phi_i(t,y')\,dy',
\qquad
\Phi_i(t,y'):=\mu(\Psi_i(y',t))\,\psi_i(\Psi_i(y',t))\,J_i(y',t).
\]
Define, for $t\in J$ and $y'\in U_i'$,
\[
F_i(y',t):=\mathbf{1}_{A_i(t)}(y')\,\Phi_i(t,y'),
\qquad\text{so that}\qquad
c_i(t)=\int_{U_i'} F_i(y',t)\,dy'.
\]

Define the compact set
\[
K_i := \Psi_i^{-1}\bigl(K\cap \operatorname{supp}\psi_i\bigr)\subset U_i'\times I_i.
\]

\smallskip
\noindent Claim. If $F_i(y',t)\neq 0$ then $(y',t)\in K_i$.

\smallskip
\noindent Proof of the claim.
If \(F_i(y',t)\ne0\), then \(x:=\Psi_i(y',t)\in\Omega\cap B_i\),
\(\mu(x)\ne0\), and \(\psi_i(x)\ne0\). Since \(t\in J\),
we have \(|f(x)|<\varepsilon\), so \(x\in K\).
Thus \(x\in K\cap\supp\psi_i\), or
equivalently \((y',t)\in K_i\).
\hfill$\square$

Since $\Phi_i$ is continuous on $U_i'\times I_i$, it is bounded on the compact set $K_i$:
there exists $C_i>0$ such that $|\Phi_i|\le C_i$ on $K_i$. Using the claim, we obtain
\[
|F_i(y',t)|\le C_i\qquad\text{for all }(y',t)\in U_i'\times J.
\]
Because $U_i'$ is bounded, $C_i\mathbf{1}_{U_i'}\in L^1(U_i')$ dominates
$\{F_i(\,\cdot\,,t)\}_{t\in J}$.

\medskip
\noindent\textbf{Step 4: Localizing the boundary in the chart.}\par
\nopagebreak[4]

Fix $i$ and set
\[
S_i := \Psi_i^{-1}\bigl(\partial\Omega\cap K\cap \operatorname{supp}\psi_i\bigr)\subset U_i'\times I_i,
\]
which is compact.

\smallskip
\noindent Claim. $S_i =\partial E_i\cap K_i$.

\smallskip
\noindent Proof of the claim.
Let \(y\in \partial E_i\cap K_i\) and set \(x:=\Psi_i(y)\).
Since \(y\in K_i\), we have \(x\in K\cap\supp\psi_i\subset B_i\), hence \(x\notin \partial B_i\).
Moreover \(y\in\partial E_i\) means \(x\in \partial(\Omega\cap B_i)\).
Because \(x\notin\partial B_i\), this forces \(x\in \partial\Omega\cap B_i\).
Together with \(x\in K\cap\supp\psi_i\), we get
\(x\in \partial\Omega\cap K\cap\supp\psi_i\), i.e. \(y\in S_i\).
Thus \(\partial E_i\cap K_i\subset S_i\).

Conversely, if \(y\in S_i\) then \(x=\Psi_i(y)\in \partial\Omega\cap K\cap\supp\psi_i\subset B_i\),
so \(x\notin\partial B_i\). Hence \(x\in \partial(\Omega\cap B_i)\), i.e. \(y\in\partial E_i\),
and also \(x\in K\cap\supp\psi_i\) gives \(y\in K_i\). Therefore \(y\in \partial E_i\cap K_i\),
so \(S_i\subset \partial E_i\cap K_i\).
\hfill\(\square\)

By hypothesis, $\partial\Omega\subset \bigcup_{j=1}^N \Sigma_j$, so
\[
S_i \subset \bigcup_{j=1}^N \Psi_i^{-1}\bigl(\Sigma_j\cap K\cap \operatorname{supp}\psi_i\bigr).
\]

Define
\[
\Sigma_{i,j}:=\Psi_i^{-1}(\Sigma_j\cap B_i)\subset U_i'\times I_i,
\qquad
\widetilde\Sigma_{i,j}:=\Sigma_{i,j}\cap K_i
=\Psi_i^{-1}\bigl(\Sigma_j\cap K\cap\supp\psi_i\bigr),
\]
and set $\Sigma_{i,j}^\circ:=\Sigma_{i,j}\setminus\partial\Sigma_{i,j}$.

\medskip
\noindent\textbf{Step 5: The boundary slice at level $0$.}\par
\nopagebreak[4]

Define the slice
\[
S_i^{0}:=\{y'\in U_i' : (y',0)\in S_i\}.
\]
We claim that $\mathcal{L}^{d-1}(S_i^{0})=0$.

Using the decomposition above,
\[
S_i^{0}
\subset
\Bigl(\ \bigcup_{\substack{1\le j\le N\\ \dim\Sigma_j=d-1}}
\pi\bigl(\widetilde\Sigma_{i,j}\cap\{y_d=0\}\bigr)\Bigr)
\ \cup\
\Bigl(\ \bigcup_{\substack{1\le j\le N\\ \dim\Sigma_j\le d-2}}
\pi\bigl(\widetilde\Sigma_{i,j}\bigr)\Bigr),
\]
where $\pi:U_i'\times I_i\to U_i'$ is the projection.

The second union has $\mathcal{L}^{d-1}$-measure zero because each $\widetilde\Sigma_{i,j}$ in the union
has Hausdorff dimension at most $d-2$ (hence $\mathcal{H}^{d-1}$-measure zero in $\mathbb{R}^d$),
and $\pi$ is Lipschitz.

Fix now $j$ with $\dim\Sigma_j=d-1$. Consider the $C^1$ function
\[
g_{i,j}:=\bigl(y\mapsto y_d\bigr)\big|_{\Sigma_{i,j}^\circ}:\Sigma_{i,j}^\circ\to\mathbb{R}.
\]
Since $f(\Psi_i(y))=y_d$, we have $g_{i,j}=(f|_{\Sigma_j^\circ}\circ \Psi_i)\big|_{\Sigma_{i,j}^\circ}$.
By assumption \textup{(ii)}, the inclusion $K\subset\overline\Omega$, and the fact that
$\Psi_i$ is a $C^1$ diffeomorphism,
the set of critical points of $g_{i,j}$ at level $0$,
\[
C_{i,j}
:=
\bigl\{y\in \widetilde\Sigma_{i,j}\cap \Sigma_{i,j}^\circ\cap\{y_d=0\}:\ dg_{i,j}(y)=0\bigr\},
\]
has $\mathcal{H}^{d-1}$-measure zero. Since $\pi$ is Lipschitz, it maps $\mathcal{H}^{d-1}$-null sets to
$\mathcal{H}^{d-1}$-null sets, and in $\mathbb{R}^{d-1}$ the measures $\mathcal{H}^{d-1}$ and $\mathcal{L}^{d-1}$
coincide up to a constant. Therefore $\mathcal{L}^{d-1}\bigl(\pi(C_{i,j})\bigr)=0$.

Set
\[
R_{i,j}:=
\bigl(\widetilde\Sigma_{i,j}\cap \Sigma_{i,j}^\circ\cap\{y_d=0\}\bigr)\setminus C_{i,j}.
\]
On $R_{i,j}$ we have $dg_{i,j}\neq 0$. For each $y\in R_{i,j}$, by continuity of $dg_{i,j}$
there exists an open neighborhood $O_y\subset \Sigma_{i,j}^\circ$ such that $dg_{i,j}\neq 0$ on $O_y$.
Let
\[
O_{i,j}:=\bigcup_{y\in R_{i,j}} O_y,
\]
an open neighborhood of $R_{i,j}$ in $\Sigma_{i,j}^\circ$ satisfying $dg_{i,j}\neq 0$ on $O_{i,j}$.
Therefore $0$ is a regular value of $g_{i,j}|_{O_{i,j}}$, and
\[
M_{i,j}:=O_{i,j}\cap\{y_d=0\}=(g_{i,j}|_{O_{i,j}})^{-1}(0)
\]
is a $C^1$ submanifold of $\Sigma_{i,j}^\circ$ of dimension $d-2$. Hence $\mathcal{H}^{d-1}(M_{i,j})=0$
and, as above, $\mathcal{L}^{d-1}\bigl(\pi(M_{i,j})\bigr)=0$. Since $R_{i,j}\subset M_{i,j}$,
we obtain $\mathcal{L}^{d-1}\bigl(\pi(R_{i,j})\bigr)=0$.

Since
\[
\widetilde\Sigma_{i,j}\cap \Sigma_{i,j}^\circ\cap\{y_d=0\}=R_{i,j}\cup C_{i,j},
\]
we conclude
\[
\mathcal{L}^{d-1}\Bigl(\pi\bigl(\widetilde\Sigma_{i,j}\cap \Sigma_{i,j}^\circ\cap\{y_d=0\}\bigr)\Bigr)=0.
\]

Finally, \(\Sigma_j\cap B_i\) is a relatively open subset of
\(\Sigma_j\). Thus the manifold boundary of \(\Sigma_{i,j}\) is the
preimage of \(\partial\Sigma_j\cap B_i\), and
\begin{equation}\label{eq:boundary}
\widetilde\Sigma_{i,j}\setminus \Sigma_{i,j}^\circ
\subset \partial\Sigma_{i,j}\cap K_i
\subset \Psi_i^{-1}(\partial \Sigma_j\cap B_i).
\end{equation}
Putting the pieces together, since $\Psi_i^{-1}(\partial \Sigma_j\cap B_i)$ has dimension at most $d-2$,
\[
\mathcal{L}^{d-1}\bigl(\pi(\widetilde\Sigma_{i,j}\cap\{y_d=0\})\bigr)=0.
\]
Altogether, $\mathcal{L}^{d-1}(S_i^{0})=0$.

\medskip
\noindent\textbf{Step 6: Pointwise convergence.}\par
\nopagebreak[4]

Fix $y'\in U_i'\setminus S_i^0$. If
$\Phi_i(0,y')=0$, then $\Phi_i(t,y')\to 0$ by continuity, and
\[
|F_i(y',t)|\leq |\Phi_i(t,y')|\longrightarrow 0=F_i(y',0).
\]

Suppose instead that $\Phi_i(0,y')\neq 0$. We claim that
$(y',0)\notin\partial E_i$. Indeed, if $(y',0)\in\partial E_i$, then, since
$\Psi_i(y',0)\in B_i$ and $B_i$ is open, we would have
\(x:=\Psi_i(y',0)\in\partial\Omega\). Since \(\Phi_i(0,y')\ne0\),
\(\psi_i(x)\ne0\). Since \(x\in\overline\Omega\) and \(f(x)=0\),
we have \(x\in K\cap\supp\psi_i\).
It would follow that $(y',0)\in S_i$, contrary to $y'\notin S_i^0$.

Therefore there exists \(\delta>0\) such that
\((y',t)\) lies either in \(E_i\) for every \(|t|<\delta\), or in its
complement for every \(|t|<\delta\). Hence, recalling that
\(
\mathbf{1}_{A_i(t)}(y')=\mathbf{1}_{E_i}(y',t),
\)
we have
\[
\mathbf{1}_{A_i(t)}(y')
=
\mathbf{1}_{A_i(0)}(y')
\qquad\text{for all sufficiently small }|t|.
\]
Since
\[
F_i(y',t)
=
\mathbf{1}_{A_i(t)}(y')\,\Phi_i(t,y')
\]
and \(t\mapsto\Phi_i(t,y')\) is continuous, it follows that
\[
F_i(y',t)\longrightarrow
\mathbf{1}_{A_i(0)}(y')\,\Phi_i(0,y')
=
F_i(y',0).
\]

\medskip
\noindent\textbf{Step 7: Dominated convergence.}\par
\nopagebreak[4]

We have pointwise convergence for a.e.\ $y'\in U_i'$ and a uniform $L^1$-dominating bound.
Thus dominated convergence gives
\[
\lim_{t\to 0} c_i(t)
= \lim_{t\to 0}\int_{U_i'} F_i(y',t)\,dy'
= \int_{U_i'} F_i(y',0)\,dy'
= c_i(0).
\]
Hence each $c_i$ is continuous at $0$, and since $c=\sum_{i=1}^M c_i$ on $J$,
$c$ is continuous at $0$.
\end{proof}

\subsection*{Proof of \Cref{cor:hausd_localized}}

\begin{proof}
We adapt the proof of \Cref{prop:hausd_continuous}, again taking \(t_0=0\).
Set \(S_0:=K\cap f^{-1}(0)\). If \(S_0=\emptyset\), compactness of \(K\)
and the support condition give \(c(t)=0\) for \(t\) sufficiently close to zero.
Otherwise, assumption \textup{(i)} and compactness give finitely many
submersion charts \(\Psi_i:U_i'\times I_i\to B_i\subset W\) covering \(S_0\),
as in Step~1 of the proof of \Cref{prop:hausd_continuous}. Set \(B:=\bigcup_i B_i\).
Since \(K\setminus B\) is compact and contains no zero of \(f\), choose
\(\varepsilon>0\) so that
\[
[-\varepsilon,\varepsilon]\subset I\cap\bigcap_i I_i,
\qquad
L:=K\cap f^{-1}([-\varepsilon,\varepsilon])\subset B.
\]
Write \(J:=(-\varepsilon,\varepsilon)\). The support condition implies
\[
\operatorname{supp}_{\overline\Omega}(\mu)\cap f^{-1}(J)\subset L.
\]

Apply Steps~2--7 of the proof of \Cref{prop:hausd_continuous}
with \(L\) in place of the compact strip \(K\), and with
\(L_i:=\Psi_i^{-1}(L\cap\supp\psi_i)\) in place of \(K_i\).
Indeed, the support condition ensures that the partition of
unity around \(L\) still gives \(c(t)=\sum_i c_i(t)\) for
\(t\in J\), and that every nonzero integrand is supported in
the corresponding \(L_i\). Since \(L\subset K\subset W\), the boundary sets involved
are contained in \(\bigcup_{j=1}^N\Sigma_j\), and
assumption \textup{(ii)} gives the required measure-zero
condition at level zero. The same dominated-convergence argument therefore
proves continuity at zero.
\end{proof}

\subsection*{Proof of \Cref{cor:manifold}}
\begin{proof}
Replacing $f$ by $f-t_0$, we may assume $t_0=0$. Thus $df_x\neq 0$ on
$\overline{\Omega}\cap f^{-1}(0)$, and we prove $c(t)\to c(0)$ as $t\to 0$.

\medskip
\noindent\textbf{Step 1: Submersion charts and a compact strip.}
Set
\[
S_0:=\overline{\Omega}\cap f^{-1}(0),
\]
which is compact. If \(S_0=\emptyset\), compactness of
\(\overline\Omega\) gives \(c(t)=0\) for \(t\) near zero, so we are done.
Otherwise, by the submersion theorem, for every $p\in S_0$ there exist an open neighborhood
$B_p\subset U$, an open set $U_p\subset\mathbb{R}^d$, and a $C^1$ diffeomorphism
$\Phi_p:U_p\to B_p$ such that, in local coordinates $y$ on $U_p$,
\[
(f\circ \Phi_p)(y)=y_d.
\]

By compactness of $S_0$, choose finitely many points $p_1,\dots,p_N\in S_0$ such that
\[
S_0\subset \bigcup_{i=1}^N B_i,\qquad B_i:=B_{p_i},\qquad \Phi_i:=\Phi_{p_i},\qquad U_i:=U_{p_i}.
\]
Shrinking $U_i$ if necessary, assume $\overline{U_i}$ is compact and $U_i$ has the product form
$U_i=U_i'\times I_i$ where $I_i\subset\mathbb{R}$ is an open interval containing $0$ and the last
coordinate is $y_d\in I_i$.

Choose $\varepsilon>0$ such that $[-\varepsilon,\varepsilon]\subset I_i$ for all $i$, and define
\[
K:=\overline{\Omega}\cap f^{-1}([-\varepsilon,\varepsilon]).
\]
Then $K$ is compact, $\Omega\cap f^{-1}(t)\subset K$ for $|t|<\varepsilon$, and (possibly shrinking $\varepsilon$)
we may assume $K\subset \bigcup_{i=1}^N B_i$.

\medskip
\noindent\textbf{Step 2: Partition of unity.}
Choose nonnegative \(C^1\) functions \(\psi_i\) with compact support
in \(B_i\), extended by zero to \(T\), whose sum is one on a neighborhood
of \(K\). For \(|t|<\varepsilon\),
\[
c(t)=\sum_{i=1}^N c_i(t),
\qquad
c_i(t):=\int_{\Omega\cap f^{-1}(t)\cap B_i}\mu(x)\psi_i(x)\,d\mathcal{H}^{d-1}(x).
\]
It suffices to show each $c_i$ is continuous at $0$.

\medskip
\noindent\textbf{Step 3: Reduction to the localization corollary.}
Fix \(i\) and set
\[
\Omega_i:=\Phi_i^{-1}(\Omega\cap B_i)\subset U_i.
\]
On \(U_i\), the function \(f\circ\Phi_i\) is the coordinate projection
\(y\mapsto y_d\). Denote its extension to all of \(\mathbb R^d\) by
\(f_i(y):=y_d\); in particular, \(df_i\ne0\) everywhere.

Let \(J_i(y)\) be the positive continuous \((d-1)\)-dimensional
Jacobian of the restriction of \(\Phi_i\) to the horizontal slice through
\(y\), computed using the Riemannian metric on \(T\). On
\(\overline{\Omega_i}\cap U_i\), define
\[
\widehat\mu_i(y):=\mu(\Phi_i(y))\,\psi_i(\Phi_i(y))\,J_i(y).
\]
Since \(\supp\psi_i\) is compact in \(B_i\), this weight vanishes
outside a compact subset of \(U_i\) and extends by zero to a continuous
nonnegative function on \(\overline{\Omega_i}\).
The area formula gives, for \(|t|<\varepsilon\),
\begin{equation}\label{eq:ci-local}
c_i(t)=\int_{\Omega_i\cap f_i^{-1}(t)}\widehat\mu_i(y)\,
d\mathcal H^{d-1}(y).
\end{equation}

The set
\[
K_i:=\Phi_i^{-1}\bigl(K\cap\supp\psi_i\bigr)
\subset\overline{\Omega_i}\cap U_i
\]
is compact, and
\[
\operatorname{supp}_{\overline{\Omega_i}}(\widehat\mu_i)
\cap f_i^{-1}((-\varepsilon,\varepsilon))\subset K_i.
\]
If \(K_i=\emptyset\), then \(c_i(t)=0\) for \(|t|<\varepsilon\).
Otherwise, choose an open neighborhood \(W_i\) of \(K_i\) with compact
closure in \(U_i\). Inside \(U_i\), the boundary of \(\Omega_i\) is
\(\Phi_i^{-1}(\partial\Omega\cap B_i)\), so
\[
\partial\Omega_i\cap W_i
\subset\bigcup_{j=1}^J\Phi_i^{-1}(\Sigma_j\cap\Phi_i(W_i)).
\]
For each \(j\), the compact set
\(\Sigma_j\cap\Phi_i(\overline{W_i})\) can be covered by
finitely many compact \(C^1\) submanifold pieces of
\(\Sigma_j\), each of dimension \(\dim\Sigma_j\), contained
in \(B_i\), and with \(C^1\) boundary. Since each piece
has the same dimension as \(\Sigma_j\), its manifold
interior lies in \(\Sigma_j^\circ\), and their tangent
spaces agree there.

The preimages of these pieces under \(\Phi_i\) are compact
\(C^1\) submanifolds with \(C^1\) boundary covering
\(\partial\Omega_i\cap W_i\). For each \(j\) with
\(\dim\Sigma_j=d-1\), the identity \(f_i=f\circ\Phi_i\)
and the chain rule imply that the critical points of
\(f_i\) on the manifold interiors of the corresponding
preimage pieces, at level zero and inside \(K_i\),
are contained in
\[
\Phi_i^{-1}\Bigl(
\bigl\{x\in\Sigma_j^\circ\cap B_i\cap\overline\Omega
\cap f^{-1}(0):
d(f|_{\Sigma_j^\circ})(x)=0\bigr\}
\Bigr).
\]
These sets have \(\mathcal H^{d-1}\)-measure zero by assumption
\textup{(ii)}, and the \(C^1\) coordinate changes preserve this nullity.
All hypotheses of \Cref{cor:hausd_localized} therefore hold for
\(\Omega_i\), \(f_i\), \(\widehat\mu_i\), and the localization set
\(K_i\). The corollary implies that \(c_i\) is continuous at zero.

\medskip
\noindent\textbf{Step 4: Conclusion.}
Each \(c_i\) is finite on an interval around \(0\), by the local integral
representation above, and is continuous at \(0\). There are finitely many
terms, so on a common interval their sum \(c\) is finite and continuous at
\(0\), hence at \(t_0\).
\end{proof}

\subsection*{Proof of \Cref{cor:hausd}}
\begin{proof}
Fix $t_0$ so that the hypotheses of \Cref{cor:manifold} hold. Write $\mu=\mu^+-\mu^-$, where
\[
\mu^+(x):=\max\{\mu(x),0\},\qquad
\mu^-(x):=\max\{-\mu(x),0\}.
\]
Then $\mu^\pm:\overline{\Omega}\to[0,\infty)$ are continuous.
Applying \Cref{cor:manifold} to $\mu^+$ and $\mu^-$ separately, we obtain
that the functions
\[
t \longmapsto \int_{\Omega\cap f^{-1}(t)} \mu^+(x)\,d\mathcal{H}^{d-1}(x)
\quad\text{and}\quad
t \longmapsto \int_{\Omega\cap f^{-1}(t)} \mu^-(x)\,d\mathcal{H}^{d-1}(x)
\]
are finite on a common neighborhood of \(t_0\) and continuous at
\(t_0\). On that neighborhood, both parts are integrable, so by linearity
\[
\int_{\Omega\cap f^{-1}(t)} \mu(x)\,d\mathcal{H}^{d-1}(x)
=
\int_{\Omega\cap f^{-1}(t)} \mu^+(x)\,d\mathcal{H}^{d-1}(x)
-
\int_{\Omega\cap f^{-1}(t)} \mu^-(x)\,d\mathcal{H}^{d-1}(x),
\]
so the left-hand side is the difference of two functions continuous at $t_0$,
hence is continuous at $t_0$.
Since $t_0$ was arbitrary, the claim follows.
\end{proof}

\subsection*{Proof of \Cref{prop:morse_cont_man}}
\begin{proof}
Replacing \(f\) by \(f-t_0\), we may assume \(t_0=0\). Let
\(\varepsilon_0>0\) be as in assumption \textup{(b)}.

\medskip
\noindent\textbf{Step 1: remove the Morse points and treat the regular part.}
By the \(C^1\) Morse lemma for \(C^2\) functions
\citep{bromberg2001extension}, for each critical point \(p_\ell\) we may choose a
\(C^1\) diffeomorphism
\[
\Phi_\ell:W_\ell\to \Phi_\ell(W_\ell)\subset \Omega,
\qquad
\Phi_\ell(0)=p_\ell,
\]
where \(W_\ell\subset\R^d\) is open, such that, writing
\(z=(u,v)\in\R^{\lambda_\ell}\times\R^{d-\lambda_\ell}\),
\[
f(\Phi_\ell(u,v))
=
-\|u\|^2+\|v\|^2
=:Q_\ell(u,v).
\]
Since the points \(p_\ell\) lie in the interior of \(\Omega\), are isolated critical
points, and are finite in number, we may choose \(\eta_\ell>0\) so small that
\[
\overline{B_{\eta_\ell}(0)}\subset W_\ell,
\qquad
\Phi_\ell(\overline{B_{\eta_\ell}(0)})\subset \Omega,
\]
the compact sets \(\Phi_\ell(\overline{B_{\eta_\ell}(0)})\) are pairwise disjoint, and
\(p_\ell\) is the only critical point of \(f\) in
\(\Phi_\ell(\overline{B_{\eta_\ell}(0)})\). Set
\[
V_\ell:=\Phi_\ell(B_{\eta_\ell}(0)).
\]
Set
\[
B:=\bigcup_{\ell=1}^L V_\ell,
\qquad
\Omega_{\rm reg}:=\Omega\setminus\overline B,
\]
with the convention that \(B=\emptyset\) if \(L=0\). Then \(\Omega_{\rm reg}\) is open,
the sets \(\overline{V_\ell}\) are pairwise disjoint and contained in \(\Omega\), and each
\(\partial V_\ell\) is a compact \(C^1\) hypersurface in \(T\), since
\(\Phi_\ell\) is a \(C^1\) diffeomorphism. Our goal is now to apply \Cref{cor:manifold} on 
\(\Omega_{\rm reg}\).

By assumption \textup{(b)}, \(f\) has no critical points in
\[
\overline\Omega\cap f^{-1}((-\varepsilon_0,\varepsilon_0))
\]
other than \(p_1,\dots,p_L\). These points lie in the open sets \(V_\ell\), and hence
do not belong to \(\overline{\Omega_{\rm reg}}\). Choose
\(0<\varepsilon<\varepsilon_0\). Then \(f\) has no critical points on
\[
\overline{\Omega_{\rm reg}}\cap f^{-1}([-\varepsilon,\varepsilon]).
\]
In particular, \(df_x\neq0\) for every
\[
x\in\overline{\Omega_{\rm reg}}\cap f^{-1}(0).
\]

Moreover,
\[
\partial\Omega_{\rm reg}
\subset
\partial\Omega\cup\bigcup_{\ell=1}^L\partial V_\ell.
\]
The pieces inherited from \(\partial\Omega\) satisfy the boundary decomposition hypothesis.
Moreover, if \(j\) is such that \(\dim\Sigma_j=d-1\), then by assumption \textup{(c)}
the set
\[
\bigl\{x\in \Sigma_j^\circ\cap \overline{\Omega}\cap f^{-1}(0):
d(f|_{\Sigma_j^\circ})(x)=0\bigr\}
\]
has \(\Hcal^{d-1}\)-measure zero. Intersecting this set with
\(\overline{\Omega_{\rm reg}}\subset\overline\Omega\), we obtain the required a.e.\
boundary regularity for the inherited boundary pieces of \(\Omega_{\rm reg}\).

We now check the
new boundary pieces. In the Morse coordinates above, \(\partial V_\ell\) corresponds to
\[
\partial B_{\eta_\ell}(0)=\{z:\|z\|^2=\eta_\ell^2\},
\]
and \(f\) corresponds to
\[
Q_\ell(u,v)=-\|u\|^2+\|v\|^2.
\]
On
\[
\{Q_\ell=0\}\cap\partial B_{\eta_\ell}(0),
\]
the differentials \(dQ_\ell\) and \(d\|z\|^2\) are linearly independent. Indeed, if
\(\lambda_\ell=0\) or \(\lambda_\ell=d\), the intersection is empty. Otherwise, at every point \((u,v)\) of this intersection,
both vectors \(u\) and \(v\) are nonzero, and the relation
\[
a\,dQ_\ell+b\,d\|z\|^2=0
\]
implies
\[
a\,\nabla Q_\ell+b\,\nabla\|z\|^2
=
a(-2u,2v)+b(2u,2v)
=
2\bigl((-a+b)u,(a+b)v\bigr)
=0,
\]
where the gradients are taken in the Euclidean coordinates
\(z=(u,v)\). Equating the \(u\)- and \(v\)-components to zero
therefore gives
\[
(-a+b)u=0,
\qquad
(a+b)v=0.
\]
Hence \(a=b\) and \(a=-b\), so \(a=b=0\). 

By the chain rule, for every \(w\in T_z\partial B_{\eta_\ell}(0)\),
we have
\(d(Q_\ell|_{\partial B_{\eta_\ell}(0)})_z(w)
=\langle\nabla Q_\ell(z),w\rangle\). Hence, at a critical point of
\(Q_\ell|_{\partial B_{\eta_\ell}(0)}\), the gradient
\(\nabla Q_\ell\) would be orthogonal to the tangent space
of the sphere. The linear independence just proved
excludes this at every point of
\(\{Q_\ell=0\}\cap\partial B_{\eta_\ell}(0)\).

Since \(\Phi_\ell\) restricts to a \(C^1\) diffeomorphism
from \(\partial B_{\eta_\ell}(0)\) onto \(\partial V_\ell\)
and \(f\circ\Phi_\ell=Q_\ell\), the chain rule therefore
implies that \(f|_{(\partial V_\ell)^\circ}\) has no
critical points at level zero.

Consequently the new hypersurface pieces \(\partial V_\ell\)
satisfy the required boundary regularity at level \(0\).

Therefore \(\Omega_{\rm reg}\) satisfies the hypotheses of \Cref{cor:manifold}
at \(t_0=0\). Hence
\[
t\longmapsto
\int_{\Omega_{\rm reg}\cap f^{-1}(t)}\mu\,d\Hcal^{d-1}
\]
is continuous at \(0\).

Since the compact sets \(\overline{V_\ell}\) are pairwise disjoint and contained in
\(\Omega\), we have
\[
\overline B=\bigcup_{\ell=1}^L\overline{V_\ell}
\]
and hence the disjoint decomposition
\[
\Omega
=
\Omega_{\rm reg}
\sqcup
\bigsqcup_{\ell=1}^L \overline{V_\ell}.
\]
Therefore, for \(t\) close to \(0\),
\[
c(t)
=
\int_{\Omega_{\rm reg}\cap f^{-1}(t)}\mu\,d\Hcal^{d-1}
+
\sum_{\ell=1}^L c_\ell(t),
\qquad
c_\ell(t):=
\int_{\overline{V_\ell}\cap f^{-1}(t)}\mu\,d\Hcal^{d-1}.
\]
It remains to prove that each \(c_\ell\) is continuous at \(0\).

\medskip
\noindent\textbf{Step 2: local terms in Morse coordinates.}
Fix \(\ell\). To simplify notation, write
\[
\Phi:=\Phi_\ell,\qquad
\eta:=\eta_\ell,\qquad
\lambda:=\lambda_\ell,\qquad
m:=d-\lambda,\qquad
Q:=Q_\ell,
\]
and set
\[
K:=\overline{B_\eta(0)}\subset\R^d.
\]
Then
\[
\overline{V_\ell}=\Phi(K),
\qquad
f\circ\Phi=Q.
\]
We prove that
\[
c_\ell(t)=\int_{\Phi(K)\cap f^{-1}(t)}\mu\,d\Hcal^{d-1}
\]
is continuous at \(0\).

\medskip
\noindent\textbf{Step 3: saddle case \(\lambda\in\{1,\dots,d-1\}\).}
Assume \(1\le\lambda\le d-1\). We treat the two one-sided limits separately.

\smallskip
\noindent Case \(t\ge0\).
Set
\[
r(t,u):=\sqrt{t+\|u\|^2},
\qquad
\Psi_t(u,\omega):=(u,r(t,u)\omega),
\qquad
\omega\in\mathbb S^{m-1}.
\]
For \(t\ge0\), define
\[
D_t^+:=\{(u,\omega):\Psi_t(u,\omega)\in K\}.
\]
Since \(K=\overline{B_\eta(0)}\), we have
\[
D_t^+
=
\Bigl\{(u,\omega):t+2\|u\|^2\le \eta^2\Bigr\}
=
\overline{B_{\sqrt{(\eta^2-t)/2}}^\lambda(0)}
\times\mathbb S^{m-1}.
\]
Fix \(t_1>0\) small enough that \(t_1<\eta^2\), and set
\[
D^+:=
\overline{B_{\eta/\sqrt2}^\lambda(0)}
\times\mathbb S^{m-1}.
\]
Then \(D_t^+\subset D^+\) for every \(t\in[0,t_1]\).

For \(t>0\), the map \(\Phi\circ\Psi_t\) is a \(C^1\) parametrization of
\(\Phi(K)\cap f^{-1}(t)\). Let
\(
\mathcal J_t^+(u,\omega)
\)
denote its \((d-1)\)-dimensional Jacobian, computed with respect to the Riemannian metric
on \(T\). By the area formula,
\[
c_\ell(t)
=
\int_{D_t^+}
\mu(\Phi(\Psi_t(u,\omega)))\,
\mathcal J_t^+(u,\omega)\,
du\,d\Hcal^{m-1}(\omega).
\]
Extending the integrand by zero from \(D_t^+\) to \(D^+\), we obtain
\[
c_\ell(t)
=
\int_{D^+}
\mathbf 1_{D_t^+}(u,\omega)\,
\mu(\Phi(\Psi_t(u,\omega)))\,
\mathcal J_t^+(u,\omega)\,
du\,d\Hcal^{m-1}(\omega).
\]
Here the entire integrand is understood to be zero outside \(D_t^+\).
We use the same convention for \(D_t^-\subset D^-\) below.

At \(t=0\), the parametrization \(\Phi\circ\Psi_0\) is \(C^1\) and injective away from
\(\{0\}\times\mathbb S^{m-1}\), whose product measure is zero, and its image is
\[
(\Phi(K)\cap f^{-1}(0))\setminus\{p_\ell\}.
\]
Since \(d\ge2\), the point \(\{p_\ell\}\) has \(\Hcal^{d-1}\)-measure zero. Hence the same
area-formula expression holds for \(c_\ell(0)\), after defining \(\mathcal J_0^+\) on
\(D^+\setminus(\{0\}\times\mathbb S^{m-1})\) and setting it arbitrarily on the null set
\(\{0\}\times\mathbb S^{m-1}\).

For \(u\neq0\), as \(t\downarrow0\), we have
\[
\Psi_t(u,\omega)\to\Psi_0(u,\omega),
\]
and the differentials \(D(\Phi\circ\Psi_t)(u,\omega)\) converge to
\(D(\Phi\circ\Psi_0)(u,\omega)\). Hence
\[
\mathcal J_t^+(u,\omega)\to\mathcal J_0^+(u,\omega)
\]
for every \(u\neq0\). Moreover,
\[
\mathbf 1_{D_t^+}(u,\omega)
=
\mathbf 1_{\{t+2\|u\|^2\le \eta^2\}}
\to
\mathbf 1_{\{2\|u\|^2\le \eta^2\}}
\]
for every \(u\) except those satisfying \(2\|u\|^2=\eta^2\), a \(du\)-null sphere in
\(\R^\lambda\). Since \(\mu\circ\Phi\) is continuous, we also have
\[
\mu(\Phi(\Psi_t(u,\omega)))\to\mu(\Phi(\Psi_0(u,\omega)))
\]
for every \(u\neq0\). Therefore the full integrand converges pointwise for a.e.\
\((u,\omega)\in D^+\).

It remains to dominate the integrands. The maps \(\Phi\) and \(D\Phi\) are bounded on the
compact set \(K\). For \(t>0\), and for \(t=0\) away from the null set
\(\{0\}\times\mathbb S^{m-1}\), the derivatives of \(\Psi_t\) are uniformly bounded:
indeed
\[
\partial_{u_i}r(t,u)=\frac{u_i}{\sqrt{t+\|u\|^2}}
\]
has absolute value at most \(1\), and \(r(t,u)\le\eta\) on \(D_t^+\). Since the Riemannian
metric is continuous and \(\Phi(K)\) is compact, the Jacobians \(\mathcal J_t^+\) are
uniformly bounded on the relevant domain, up to the null set at \(t=0\). Since \(\mu\) is
bounded on \(\overline\Omega\), the integrands are dominated by an integrable constant on
the finite-measure set \(D^+\). Dominated convergence gives
\[
\lim_{t\downarrow0}c_\ell(t)=c_\ell(0).
\]

\smallskip
\noindent Case \(t\le0\).
Now \(Q(u,v)=t\) is equivalent to
\[
\|u\|^2=\|v\|^2+(-t).
\]
Let \(\omega'\in\mathbb S^{\lambda-1}\), and define
\[
\widetilde\Psi_t(v,\omega')
:=
\bigl(\sqrt{-t+\|v\|^2}\,\omega',v\bigr).
\]
For \(t\le0\), set
\[
D_t^-:=\{(v,\omega'):\widetilde\Psi_t(v,\omega')\in K\}.
\]
As before,
\[
D_t^-
=
\Bigl\{(v,\omega'):(-t)+2\|v\|^2\le\eta^2\Bigr\}
\subset
D^-:=
\overline{B_{\eta/\sqrt2}^{m}(0)}\times\mathbb S^{\lambda-1}
\]
for \(t\in[-t_1,0]\), with \(t_1>0\) sufficiently small.

Repeating the same argument with the parametrizations \(\Phi\circ\widetilde\Psi_t\), their
\((d-1)\)-dimensional Jacobians, and dominated convergence on the fixed finite-measure set
\(D^-\), we obtain
\[
\lim_{t\uparrow0}c_\ell(t)=c_\ell(0).
\]
Combining the two one-sided limits proves that \(c_\ell\) is continuous at \(0\) in the
saddle case.

\medskip
\noindent\textbf{Step 4: local minima and local maxima.}
If \(\lambda=0\), then \(Q(v)=\|v\|^2\). Hence \(\{Q=t\}=\emptyset\) for \(t<0\), and
\(\{Q=0\}=\{0\}\). Since \(d\ge2\), we have
\[
\Hcal^{d-1}(\{p_\ell\})=0,
\]
so \(c_\ell(t)=0\) for \(t<0\) and \(c_\ell(0)=0\). For \(0<t<\eta^2\), the set
\(\{Q=t\}\cap K\) is \(\partial B_{\sqrt t}(0)\). Since \(\Phi\) is Lipschitz on \(K\) as a map from the Euclidean metric to the
Riemannian distance on \(T\),
and since \(\mu\) is bounded, there is a constant \(C>0\) such that
\[
0\le c_\ell(t)
\le
C\,\Hcal^{d-1}\bigl(\partial B_{\sqrt t}(0)\bigr)
=
C'\,t^{(d-1)/2}.
\]
Thus \(c_\ell(t)\to0=c_\ell(0)\) as \(t\to0^+\).

If \(\lambda=d\), then \(Q(u)=-\|u\|^2\), and the same argument, with the signs reversed,
shows that \(c_\ell(t)=0\) for \(t>0\), \(c_\ell(0)=0\), and \(c_\ell(t)\to0\) as
\(t\to0^-\).

Therefore \(c_\ell\) is continuous at \(0\) also in the minimum and maximum cases.

\medskip
\noindent\textbf{Step 5: conclusion.}
The regular part contribution is continuous at \(0\) by Step~1, and each local term
\(c_\ell\) is continuous at \(0\) by Steps~2--4. Since the sum is finite, \(c\) is
continuous at \(0\), and hence at \(t_0\).
\end{proof}

\FloatBarrier
\setcounter{table}{0}
\setcounter{equation}{0}
\section{Proofs from Section 4}\label{app:circular_flow_proofs}
This appendix contains the proofs deferred from \Cref{sec:circ_flow}.

\subsection*{Proof of \Cref{prop:equivalence}}
\begin{proof}
We use the notation \(\ff_{in}\), \(\ff_{out}\), and \(\ff\) introduced in
\Cref{sec:circulations}.

    Consider $e=\{e^{i,i+1}_{j,k}\}\in \Adm(K)$.
The values of $e$ clearly define a $1$-chain on $G_K$, but we want to \virgolette{roll} such chain back on $K$ and so we define $s_e:E_K\rightarrow \R$ as follows:

\begin{itemize}
\item for \(i=1,\ldots,n-1\), set
\[
s_e((v^i_j,v^{i+1}_k))=e^{i,i+1}_{j,k}
\]
for all \(j,k\) such that \((v^i_j,v^{i+1}_k)\in \Out_{K}(v^i_j)\);
    \item to go from \virgolette{top to bottom}, instead, set
\[
s_e((v^n_j,v^1_k))=e^{n,n+1}_{j,k}
\]
for all \(j,k\) such that \((v^n_j,v^1_k)\in \Out_K(v^n_j)\).
\end{itemize}

Viceversa, given $s\in\Fop_1(K,\R)$, we can define $e_s \in \Adm(K)$ as:
\begin{itemize}
    \item for \(i=1,\ldots,n-1\), set
\[
(e_s)^{i,i+1}_{j,k}=s((v^i_j,v^{i+1}_k))
\]
for all \(j,k\) such that \((v^i_j,v^{i+1}_k)\in \Out_K(v^i_j)\);
    \item set
\[
(e_s)^{n,n+1}_{j,k}=s((v^n_j,v^1_k))
\]
for all \(j,k\) such that \((v^n_j,v^1_k)\in \Out_K(v^n_j)\);
    \item lastly, for all $j=1,\ldots,m_1$, define $(e_s)^{0,1}_{0,j}=(e_s)^{n+1,n+2}_{j,0}=\ff(s)(v^1_j)$.
\end{itemize}

Now we prove $e_s \in \Adm(K)$ and that these two correspondences are inverses to each other.

\begin{itemize}
      \item We start by proving that, given $s\in\Fop_1(K,\R)$, $e_s \in \Adm(K)$. We need to prove that $e_s$ satisfies constraints 1. and 2.:
        \begin{itemize}
            \item[(1.a)] for each $v_j^i$ vertex in $G_K$ with $i=2,\ldots,n-1$ we have:
            \begin{align*}
                &\sum_{(v^i_j,v^{i+1}_k)\in \Out(v_j^i)} (e_s)^{i,i+1}_{j,k} = 
                \sum_{(v^i_j,v^{i+1}_k)\in \Out(v_j^i)} s((v^i_j,v^{i+1}_k)) = \\
                &\sum_{(v^{i-1}_k,v^i_j)\in \In(v_j^i)} s((v^{i-1}_k,v^i_j)) =
                \sum_{(v^{i-1}_k,v^i_j)\in \In(v_j^i)} (e_s)^{i-1,i}_{k,j}
            \end{align*}
 \item[(1. b)] For \(v_j^1\), we have
\[
\begin{aligned}
\sum_{(v^1_j,v^2_k)\in \Out_{G_K}(v^1_j)}
(e_s)^{1,2}_{j,k}
&=
\sum_{(v^1_j,v^2_k)\in \Out_K(v^1_j)}
s((v^1_j,v^2_k)) \\
&=\ff(s)(v_j^1) \\
&=(e_s)^{0,1}_{0,j} \\
&=
\sum_{(v^0_0,v^1_j)\in \In_{G_K}(v^1_j)}
(e_s)^{0,1}_{0,j}.
\end{aligned}
\]
\item[(1. c)] We now consider the same scenario as above but with $i=n$. For each vertex $v_j^{n}$ in $G_K$ we have:
            \begin{align}
                &\sum_{(v^n_j,v^{n+1}_k)\in \Out_{G_K}(v_j^{n})} (e_s)^{n,n+1}_{j,k} = 
                \sum_{(v^n_j,v^{1}_k)\in \Out_K(v_j^{n})} s((v^n_j,v^{1}_k)) = \\
                &\sum_{(v^{n-1}_k,v^n_j)\in \In_K(v^n_j)} s((v^{n-1}_k,v^{n}_j)) =
                \sum_{(v^{n-1}_k,v^n_j)\in \In_{G_K}(v^n_j)} (e_s)^{n-1,n}_{k,j}.
            \end{align}
            \item[(1.d)] Finally, for $v_j^{n+1}$, we have
\[
\begin{aligned}
\sum_{(v_j^{n+1},v_0^{n+2})\in \Out_{G_K}(v_j^{n+1})}
(e_s)^{n+1,n+2}_{j,0}
&=(e_s)^{n+1,n+2}_{j,0}\\
&=\ff(s)(v_j^1)\\
&=
\sum_{(v^n_k,v^1_j)\in \In_K(v_j^1)}
s((v^n_k,v^1_j))\\
&=
\sum_{(v^n_k,v^{n+1}_j)\in \In_{G_K}(v_j^{n+1})}
(e_s)^{n,n+1}_{k,j}.
\end{aligned}
\]
            \item[(2.)] We are left to check that $(e_s)^{0,1}_{0,j}=(e_s)^{n+1,n+2}_{j,0}$ for every $j$. But this is true by construction. 
            \end{itemize}
            Thus $e_s\in\Adm(K)$.
    \item Now we prove that $s_e\in\Fop_1(K,\R)$.
          We need to check that for every vertex $v$ in $K$, $\partial(s_e)(v)=0$. For any vertex $v_j^i$ with $i\neq 1,n$ this is easily implied by constraint 1. defining $\Adm(K)$. Thus, consider $v^n_j$ for $j\in\{1,\ldots,m_n\}$. 
            \begin{align*}
            \textstyle\ff_{out}(s_e)(v^n_j) =
            &\displaystyle\sum_{(v^n_j,v^1_k)\in \Out_K(v^n_j)} s_e((v^n_j,v^1_k))=\sum_{(v^n_j,v^{n+1}_k)\in \Out_{G_K}(v^n_j)} e^{n,n+1}_{j,k} = \\
            &\sum_{(v^{n-1}_k,v^n_j)\in \In_{G_K}(v^n_j)} e^{n-1,n}_{k,j} =
            \sum_{(v^{n-1}_k,v^n_j)\in \In_K(v^n_j)} s_e((v^{n-1}_k,v^n_j)) =\textstyle\ff_{in}(s_e)(v^n_j).
            \end{align*}

            Similarly, take $v^1_j$ for some $j\in\{1,\ldots,m_1\}$. 
            \begin{align*}
            \textstyle\ff_{out}(s_e)(v^1_j) =&\displaystyle\sum_{(v^1_j,v^2_k)\in \Out_K(v^1_j)} s_e((v^1_j,v^2_k))=\sum_{(v^1_j,v^2_k)\in \Out_{G_K}(v^1_j)} e^{1,2}_{j,k} = \\
            &\sum_{(v^{0}_0,v^1_j)\in \In_{G_K}(v^1_j)} e^{0,1}_{0,j} = e^{0,1}_{0,j} = e^{n+1,n+2}_{j,0}= \\
            &\sum_{(v^{n}_k,v^{n+1}_j)\in\In_{G_K}(v^{n+1}_j)} e^{n,n+1}_{k,j} = \sum_{(v^{n}_k,v^{1}_j)\in\In_K(v^{1}_j)} s_e((v^{n}_k,v^{1}_j)) = 
            \textstyle\ff_{in}(s_e)(v^1_j).
            \end{align*}

    \item Now we prove $e_{s_e}=e$. 
        \begin{itemize}
            \item for \(i=1,\ldots,n-1\), and for all \(j,k\) such that
\((v^i_j,v^{i+1}_k)\in \Out_{G_K}(v^i_j)\), we have
\[
(e_{s_e})^{i,i+1}_{j,k}
=
s_e((v^i_j,v^{i+1}_k))
=
e^{i,i+1}_{j,k}.
\]
            \item for all \(j,k\) such that \((v^n_j,v^1_k)\in \Out_K(v^n_j)\) instead, we have: $(e_{s_e})^{n,n+1}_{j,k}=s_e((v^n_j,v^{1}_k))=e^{n,n+1}_{j,k}$;
            \item lastly, consider any $j=1,\ldots,m_1$, then:
            \begin{align*}
                &(e_{s_e})^{0,1}_{0,j}=(e_{s_e})^{n+1,n+2}_{j,0}=\ff(s_e)(v^1_j) =\\
                &\sum_{(v^{n}_k,v^{1}_j)\in\In_K(v^{1}_j)} s_e((v^{n}_k,v^{1}_j)) = \sum_{(v^{n}_k,v^{n+1}_j)\in\In_{G_K}(v^{n+1}_j)} e^{n,n+1}_{k,j} =e^{0,1}_{0,j}=e^{n+1,n+2}_{j,0}.
            \end{align*}

        \end{itemize}

    \item Now we need to prove $s_{e_s}=s$. We have $s_{e_s}((v^i_j,v^{i+1}_k))=(e_s)^{i,i+1}_{j,k}=s((v^i_j,v^{i+1}_k))$ for all $i=1,\ldots,n-1$. We just need to check edges of the form $(v_j^n,v^1_k)$. We have: 
    $s_{e_s}((v^n_j,v^{1}_k))=(e_s)^{n,n+1}_{j,k}=s((v^n_j,v^{1}_k))$. So we are done.
    
    \item In the points above, we have already shown that:
    \[
    (e_{s_e})^{0,1}_{0,j}=\textstyle\ff(s_e)(v^1_j)
    \]
    and so:
    \[
    \textstyle\sum_{(v^0_0,v^{1}_j)\in \Out_{G_K}(v_0^0)} e^{0,1}_{0,j} = \sum_{j=1}^{m_1} \ff(s_e)(v^1_j)=\Fop_1(f)(s_e).
    \]
  \end{itemize}
\end{proof}

\subsection*{Proof of \Cref{teo:unrolled_flow}}
\begin{proof}
For the following proof, for simplicity, we write $F$ for $\Fop_1^+(f)$ and $F'$ for $\Fop_1^+(f')$.

Note that, by \Cref{cor:fix_iso}, for every $s\in \Fop_1^+(K,\R)$ and for every $i$, we have:
\[
\textstyle F(s)=\sum_{j=1}^{m_i}\ff_{in}(s)(v^i_j)=\sum_{j=1}^{m_i}\ff_{out}(s)(v^i_j).
\]

Set $L = \max_i \sum_{e\in E_i} C(e)$, and consider $H_L$ as in \Cref{lem:finiteness}.
Fix the first layer, i.e. $i=1$ and take \(q\in \Fop_1^+(K',\N)\) with \(q\le p'^*C\). We know that there exist $r<r'$ such that  
\[
    q((v_j^{1,r},v_k^{2,r}))=
    q((v_j^{1,r'},v_k^{2,r'})),
    \]
    for every arc \((v_j^{1,r},v_k^{2,r})\in E_1^r\).
To simplify the calculations, we translate deck indices and suppose \(r=0\);
then \(N:=r'=r'-r\).

Since \(q\) is a circulation, and since the profiles on \(E_1^0\) and \(E_1^N\)
agree, we have
\[
\ff_{in}(q)(v^{1,0}_k)=\ff_{out}(q)(v^{1,0}_k)
=\ff_{out}(q)(v^{1,N}_k)=\ff_{in}(q)(v^{1,N}_k)
\]
for every \(k=1,\ldots,m_1\).
We now define $s_N:E_K\rightarrow \R$ as the average of $q$ over the first $N$ decks. For every $i=1,\ldots,n-1$ we set:
\begin{align*}
s_N((v_j^{i},v_k^{i+1})):=\frac{1}{N} \sum_{u=0}^{N-1} u.q((v_j^{i,0},v_k^{i+1,0})), 
\end{align*}
and, analogously,
\begin{align*}
s_N((v_j^{n},v_k^{1})):=\frac{1}{N} \sum_{u=0}^{N-1} u.q((v_j^{n,0},v_k^{1,1})). 
\end{align*}
We need to check $s_N\in \Fop_1^+(K,\Q)$ and $s_N\leq C$. Since clearly $s_N(e)\in \Q$, we show $\partial s_N\equiv 0$ and $s_N\leq C$. 

First we consider $i\neq 1,n$.

\begin{align*}
&\textstyle\ff_{in}(s_N)(v^i_j)=
\displaystyle \sum_{(v_k^{i-1},v_j^{i}) 
\in \In(v^i_j)} s_N((v_k^{i-1},v_j^{i})) =
\sum_{(v_k^{i-1},v_j^{i}) \in \In(v^i_j)} \frac{1}{N}\sum_{u=0}^{N-1} u.q((v_k^{i-1,0},v_j^{i,0}))=\\
&\sum_{(v_k^{i-1,0},v_j^{i,0}) \in 
\In(v^{i,0}_j)} \frac{1}{N}\sum_{u=0}^{N-1} q((v_k^{i-1,u},v_j^{i,u}))=
 \frac{1}{N}\sum_{u=0}^{N-1}
 \sum_{(v_k^{i-1,u},v_j^{i,u}) \in 
\In(v^{i,u}_j)}
 q((v_k^{i-1,u},v_j^{i,u}))= \\
&\frac{1}{N}\sum_{u=0}^{N-1} \textstyle\ff_{in}(q)(v^{i,u}_j) = \displaystyle\frac{1}{N}\sum_{u=0}^{N-1}
\textstyle\ff_{out}(q)(v^{i,u}_j) = 
\displaystyle\frac{1}{N}\sum_{u=0}^{N-1}
 \sum_{(v_j^{i,u},v_k^{i+1,u}) \in 
\Out(v^{i,u}_j)}
 q((v_j^{i,u},v_k^{i+1,u}))= \\
& \sum_{(v_j^{i},v_k^{i+1}) \in \Out(v^i_j)} \frac{1}{N}\sum_{u=0}^{N-1} u.q((v_j^{i,0},v_k^{i+1,0}))=\textstyle\ff_{out}(s_N)(v^i_j).
\end{align*}

The case $i=n$ is analogous, we just need to replace $i$ with $n$ and the last line with:
\begin{align*}
&\displaystyle\frac{1}{N}\sum_{u=0}^{N-1}
 \sum_{(v_j^{n,u},v_k^{1,u+1}) \in 
\Out(v^{n,u}_j)}
 q((v_j^{n,u},v_k^{1,u+1}))= \\
 & \sum_{(v_j^{n},v_k^{1}) \in \Out(v^n_j)} \frac{1}{N}\sum_{u=0}^{N-1} u.q((v_j^{n,0},v_k^{1,1}))=\textstyle\ff_{out}(s_N)(v^n_j).
\end{align*}

Consider now $i=1$.

\begin{align*}
&\textstyle\ff_{in}(s_N)(v^1_j)=
\displaystyle \sum_{(v_k^{n},v_j^{1}) 
\in \In(v^1_j)} s_N((v_k^{n},v_j^{1})) =
\sum_{(v_k^{n},v_j^{1}) \in \In(v^1_j)} \frac{1}{N}\sum_{u=0}^{N-1} u.q((v_k^{n,0},v_j^{1,1}))=\\
&\sum_{(v_k^{n,0},v_j^{1,1}) \in 
\In(v^{1,1}_j)} \frac{1}{N}\sum_{u=0}^{N-1} q((v_k^{n,u},v_j^{1,u+1}))=
 \frac{1}{N}\sum_{u=0}^{N-1}
 \sum_{(v_k^{n,u},v_j^{1,u+1}) \in 
\In(v^{1,u+1}_j)}
 q((v_k^{n,u},v_j^{1,u+1}))= \\
&\displaystyle\frac{1}{N}
 \sum_{u=0}^{N-1} \textstyle\ff_{in}(q)(v^{1,u+1}_j) = 
\displaystyle\frac{1}{N} \left(
 \sum_{u=0}^{N-2} \textstyle\ff_{in}(q)(v^{1,u+1}_j)
 +\textstyle\ff_{in}(q)(v^{1,N}_j) \right)=\\
&\displaystyle\frac{1}{N} \left(
 \sum_{u=0}^{N-2} \textstyle\ff_{in}(q)(v^{1,u+1}_j)
 +\textstyle\ff_{in}(q)(v^{1,0}_j) \right) = 
 \displaystyle\frac{1}{N} 
 \sum_{u=0}^{N-1} \textstyle\ff_{in}(q)(v^{1,u}_j)=\\
&\displaystyle\frac{1}{N}\sum_{u=0}^{N-1}
\textstyle\ff_{out}(q)(v^{1,u}_j) =
\displaystyle\frac{1}{N}\sum_{u=0}^{N-1}
 \sum_{(v_j^{1,u},v_k^{2,u}) \in 
\Out(v^{1,u}_j)}
 q((v_j^{1,u},v_k^{2,u}))= \\
& \sum_{(v_j^{1},v_k^{2}) \in \Out(v^1_j)} \frac{1}{N}\sum_{u=0}^{N-1} u.q((v_j^{1,0},v_k^{2,0}))=\textstyle\ff_{out}(s_N)(v^1_j).
\end{align*}

Now we show $s_N\leq C$:
\[
s_N((v_k^{i-1},v_j^{i})) =
 \frac{1}{N}\sum_{u=0}^{N-1} u.q((v_k^{i-1,0},v_j^{i,0})) \leq \frac{1}{N}\sum_{u=0}^{N-1} p'^*C((v_k^{i-1,u},v_j^{i,u})) = C((v_k^{i-1},v_j^{i})).
\]
The same argument applies to wrap-around edges \((v_j^n,v_k^1)\), using the
lifts \((v_j^{n,u},v_k^{1,u+1})\).

To conclude this step of the proof, for any $i=1,\ldots,n$ we have:
\begin{align*}
&F(s_N)=\sum_{j=1}^{m_i}\ff(s_N)(v_j^i) = \sum_{j=1}^{m_i}\frac{1}{N}\sum_{u=0}^{N-1}
\textstyle\ff(q)(v^{i,u}_j)=
\displaystyle\frac{1}{N}\sum_{u=0}^{N-1}\sum_{j=1}^{m_i}
\textstyle\ff(q)(v^{i,u}_j)= \\
&\displaystyle\frac{1}{N}\sum_{u=0}^{N-1}
\textstyle F'(q) = 
\textstyle F'(q).
\end{align*}

To remove the assumption $r=0$ one needs to set 
\begin{align*}
s_N((v_j^{i},v_k^{i+1})):=\frac{1}{N} \sum_{u=0}^{N-1} u.q((v_j^{i,r},v_k^{i+1,r})), 
\end{align*}
and 
\begin{align*}
s_N((v_j^{n},v_k^{1})):=\frac{1}{N} \sum_{u=0}^{N-1} u.q((v_j^{n,r},v_k^{1,r+1})), 
\end{align*}
 and repeat analogous calculations.

So far, for every \(q\in\Fop_1^+(K',\N)\) with \(q\le p'^*C\), we have built
\(s_N\in\Fop_1^+(K,\Q)\) with \(s_N\le C\).

As a last step in the proof, we only need to find $H$ independent on $q$ and $s_N$ such that  $H\cdot s_N\in \Fop_1^+(K,\Z)$. Since $\mid r-r'\mid \leq H_L$, if we consider $H=H_L!=\prod_{i=0}^{H_L-1} (H_L-i)$, we always have $H\cdot s_N\in \Fop_1^+(K,\Z)$.

\end{proof}

\subsection*{Proof of \Cref{cor:real_capacities}}
\begin{proof}
We assume \(E_K\ne\varnothing\); otherwise both max-flow values are zero
and the statement is trivial.
By \Cref{teo:integer}, the equality holds for every integer-valued capacity function.
Since \(K\) is finite, multiplying by a common denominator shows that the same equality
also holds for every rational-valued capacity function
\[
C:E_K\to \mathbb{Q}_{\ge 0}.
\]

We now pass to arbitrary real-valued capacities.

If \(C_1,C_2:E_K\to \mathbb{R}_{\ge 0}\) satisfy \(C_1\le C_2\) arc-wise, then
\[
\textstyle
\{s\in \Fop_1^+(K,\mathbb R)\mid s\le C_1\}
\subset
\{s\in \Fop_1^+(K,\mathbb R)\mid s\le C_2\},
\]
hence
\[
\f(f,C_1)\le \f(f,C_2).
\]
Similarly,
\[
\f(f',p'^*C_1)\le \f(f',p'^*C_2).
\]

First of all we prove that \(C\mapsto \f(f,C)\) is continuous.
Since \(K\) is finite, the set \(E_K\) of directed edges of \(K\) is finite. Fix capacities
\(C:E_K\to \mathbb R_{\ge 0}\) and a sequence \(C_n:E_K\to \mathbb R_{\ge 0}\) such that \(C_n\to C\) arc-wise. We claim that
\[
\f(f,C_n)\to \f(f,C).
\]

Set
\[
M:=\max_{e\in E_K} C(e)+1.
\]
Up to discarding finitely many indices, we may assume \(C_n(e)\le M\) for every
\(e\in E_K\) and every \(n\). For each \(n\), choose
\(s_n\in \Fop_1^+(K,\mathbb R)\) such that
\[
\textstyle
s_n\le C_n
\qquad\text{and}\qquad
\Fop_1(f)(s_n)=\f(f,C_n).
\]
Such a maximizer exists: on the finite-dimensional space \(\mathbb R^{E_K}\), the set
\[
\textstyle
\{s\in \Fop_1^+(K,\mathbb R)\mid s\le C_n\}
\]
is non-empty and compact, and \(\Fop_1(f)\) is linear.

We first prove upper semicontinuity. Let \((n_k)_k\) be a subsequence such that
\[
\lim_{k\to\infty}\f(f,C_{n_k})
=
\limsup_{n\to\infty}\f(f,C_n).
\]
Since \(0\le s_{n_k}(e)\le M\) for every \(e\in E_K\), the sequence
\((s_{n_k})_k\) lies in the compact cube \([0,M]^{E_K}\). Passing to a
further subsequence, not relabeled, there exists
\(
s_\infty\in[0,M]^{E_K}
\)
such that \(s_{n_k}\to s_\infty\) arc-wise. We claim that \(s_\infty\)
is admissible for \(\f(f,C)\).

Indeed, since
\(
\Fop_1(K,\mathbb R)=\ker(\partial_1)
\)
and \(\partial_1\) is linear, \(\Fop_1(K,\mathbb R)\) is closed in
\(\mathbb R^{E_K}\). Hence
\(
s_\infty\in\Fop_1^+(K,\mathbb R).
\)
Moreover, for every \(e\in E_K\),
\[
s_\infty(e)
=
\lim_{k\to\infty}s_{n_k}(e)
\le
\lim_{k\to\infty}C_{n_k}(e)
=
C(e).
\]
Thus \(s_\infty\le C\), and \(s_\infty\) is admissible for
\(\f(f,C)\). Therefore,
\[
\begin{aligned}
\limsup_{n\to\infty}\f(f,C_n)
&=
\lim_{k\to\infty}\f(f,C_{n_k})\\
&=
\lim_{k\to\infty}\Fop_1(f)(s_{n_k})\\
&=
\Fop_1(f)(s_\infty)
\le
\f(f,C).
\end{aligned}
\]
This proves upper semicontinuity.

We now prove lower semicontinuity, namely
\[
\f(f,C)\le \liminf_{n\to\infty}\f(f,C_n).
\] 
Fix \(\varepsilon>0\) and choose
\(s\in \Fop_1^+(K,\mathbb R)\) with \(s\le C\) such that
\[
\textstyle
\Fop_1(f)(s)\ge \f(f,C)-\varepsilon.
\]
If \(s=0\), then \(\f(f,C)\le \varepsilon\) and \(\liminf_{n\to\infty}\f(f,C_n)\ge 0\ge \f(f,C)-\varepsilon.\) So the claim holds. Otherwise define
\[
\alpha_n
:=
\min\left(
1,\,
\min_{\{e\in E_K:\,s(e)>0\}}
\frac{C_n(e)}{s(e)}
\right)
\]
and set
\[
s_n':=\alpha_n s.
\]
Since \(C_n(e)\ge 0\) for every \(e\), we have \(\alpha_n\ge 0\), and hence
\(s_n'\in \Fop_1^+(K,\mathbb R)\). Moreover, if \(s(e)=0\), then
\[
s_n'(e)=0\le C_n(e),
\]
while, if \(s(e)>0\), the definition of \(\alpha_n\) gives
\[
s_n'(e)=\alpha_n s(e)\le C_n(e).
\]
Thus \(s_n'\le C_n\). Therefore
\[
\textstyle
\f(f,C_n)
\ge
\Fop_1(f)(s_n')
=
\alpha_n\Fop_1(f)(s).
\]

For every arc \(e\) with \(s(e)>0\), since \(C_n(e)\to C(e)\) and \(s(e)\le C(e)\), we have
\[
\frac{C_n(e)}{s(e)}
\longrightarrow
\frac{C(e)}{s(e)}
\ge 1.
\]
Since \(E_K\) is finite, it follows that \(\alpha_n\to 1\). 
Taking \(n\to\infty\), we obtain
\[
\textstyle
\liminf_{n\to\infty}\f(f,C_n)
\ge
\Fop_1(f)(s)
\ge
\f(f,C)-\varepsilon.
\]
Since \(\varepsilon>0\) is arbitrary,
\[
\liminf_{n\to\infty}\f(f,C_n)\ge \f(f,C).
\]
Together with the upper semicontinuity proved above, this shows
\(
\f(f,C_n)\to \f(f,C).
\)

We now combine this continuity result with the density of \(\mathbb Q\) in \(\mathbb R\). Take an arbitrary \(C:E_K\to \mathbb R_{\ge 0}\). We can choose
rational-valued capacities
\[
C_n^- ,\, C_n^+ :E_K\to \mathbb Q_{\ge 0}
\]
such that, for every arc \(e\in E_K\),
\[
C_n^-(e)\le C(e)\le C_n^+(e),
\qquad
C_n^-(e)\uparrow C(e),
\qquad
C_n^+(e)\downarrow C(e).
\]
For example, one may take
\[
C_n^-(e)=\frac{\lfloor 2^nC(e)\rfloor}{2^n},
\qquad
C_n^+(e)=\frac{\lceil 2^nC(e)\rceil}{2^n}.
\]

Since \(C_n^\pm\) are rational-valued, the rational case gives
\[
\f(f,C_n^-)=\f(f',p'^*C_n^-),
\qquad
\f(f,C_n^+)=\f(f',p'^*C_n^+)
\]
for every \(n\).

By the monotonicity we showed at the beginning of the proof,
\[
\f(f',p'^*C_n^-)\le \f(f',p'^*C)\le \f(f',p'^*C_n^+),
\]
hence
\[
\f(f,C_n^-)\le \f(f',p'^*C)\le \f(f,C_n^+).
\]
Now \(C_n^\pm\to C\) arc-wise, so by the previous part of the proof,
\[
\f(f,C_n^-)\to \f(f,C),
\qquad
\f(f,C_n^+)\to \f(f,C).
\]
Passing to the limit in the previous inequalities yields
\[
\f(f',p'^*C)=\f(f,C).
\]
This proves the claim.
\end{proof}

\vskip 0.2in
\bibliography{references}

\end{document}